\documentclass[11pt,a4paper,reqno]{article}

\usepackage[margin=2.5cm]{geometry}
\usepackage[T1]{fontenc}
\usepackage[utf8]{inputenc}
\usepackage{alphabeta}

\usepackage[varqu,varl]{inconsolata} % sans serif typewriter
\usepackage{amssymb}
\usepackage{amsmath}
\usepackage{amsthm}
\usepackage{color}
\usepackage[dvipsnames]{xcolor}
\usepackage{graphicx}
\usepackage{subfigure}
\usepackage[pagewise]{lineno}%turn on/off by \linenumbers below
\usepackage{bm}
\usepackage{cancel}
\usepackage{enumerate}
\usepackage{ifthen}
\usepackage{listings}
\usepackage{mathtools}
\usepackage{paralist}
\usepackage{wrapfig}
\usepackage{url}
\usepackage{xcolor}
\usepackage{xspace}
\usepackage{blkarray}
\usepackage{hyperref}
\usepackage{bbm}
\usepackage[norefs,nocites]{refcheck}
\usepackage{pdflscape}
\usepackage{rotating, caption}
\usepackage{pgfplots}
\pgfplotsset{compat=1.11}
\usetikzlibrary{arrows,matrix}
\usetikzlibrary{shapes.geometric}

\usepackage[siunitx,RPvoltages]{circuitikz}
\usepackage{tikz}
\usetikzlibrary{arrows}

\newcommand\nc{\newcommand}
\nc\rnc{\renewcommand}

\nc\TW{\textwidth}

\makeatletter
\def\rf#1{(\@rf#1,.)}
\def\@rf#1,{\ref{eq:#1}\@ifnextchar . {\@endrf}{, \@rf}}
\def\@endrf.{}
\makeatother

\nc\alref[1]{Algorithm~\ref{al:#1}\xspace}
\nc\apref[1]{Appendix~\ref{sc:#1}\xspace}
\nc\coref[1]{Corollary~\ref{co:#1}\xspace}
\nc\dfref[1]{Definition~\ref{df:#1}\xspace}
\nc\exref[1]{Example~\ref{ex:#1}\xspace}
\nc\fgref[1]{Figure~\ref{fg:#1}\xspace}
\nc\lmref[1]{Lemma~\ref{lm:#1}\xspace}
\nc\scref[1]{Section~\ref{sc:#1}\xspace}
\nc\ssref[1]{Subsection~\ref{ss:#1}\xspace}
\nc\tbref[1]{Table~\ref{tb:#1}\xspace}
\nc\thref[1]{Theorem~\ref{th:#1}\xspace}

\nc\aprf[1]{App.~\ref{sc:#1}\xspace}
\nc\fgrf[1]{Fig.~\ref{fg:#1}\xspace}
\nc\scrf[1]{\S\ref{sc:#1}\xspace}
\nc\ssrf[1]{\S\ref{ss:#1}\xspace}

\nc\mc[3]{\multicolumn{#1}{#2}{#3}}
\rnc\_[1]{\mathbf{#1}}
\rnc\`{^{-1}} %inverse of a matrix
\rnc\:{\,{:}\,} %for range like 1:n; ordinary : has too much space
\nc\x{\times}
\rnc\*{{\times}} %suppresses space round it, I can't find what \* originally means
\nc\ds{\displaystyle}
\nc\ninf{{-\!\infty}}
\nc\R{\mathbb{R}}
\nc\Ap{\widetilde{A}}

\rnc{\d}{\mathrm{d}}
\nc\val{\mathop{\textup{Val}}}
\nc\sij{\sigma_{ij}}
\nc\set[2]{\{\,#1\mid\mbox{#2}\,\}}
\nc\eset[1]{\{\,\mbox{#1}\,\}}
\nc\elst[1]{[\,\mbox{#1}\,]}

\nc\dbd[2]{\frac{\partial #1}{\partial #2}} %displayed generic partial derivative
\nc\tdbd[2]{\partial #1/\partial #2} %generic partial derivative for running text

\nc\matlab{{\sc Matlab}\xspace}
\nc\matlabtm{\matlab{}\texttrademark\xspace}
\nc\SA{structural analysis\xspace}
\nc\Smethod{$\Sigma$-method\xspace}
\nc\sigmx{signature matrix\xspace}
\nc\sysJ{system Jacobian\xspace}
\nc\ccmatrix{Kron matrix\xspace} %or cycle-cutset matrix

\nc\para[1]{\par\medskip\noindent{\bfseries\slshape#1.}}
\nc\HI[1]{{\color{red}#1}}
\nc\HB[1]{{\color{blue}#1}}

\rnc\.[2]{\ifthenelse{#1=1}{\dot{#2}}{\ifthenelse{#1=2}{\ddot{#2}}{\ifthenelse{#1=3}{\dddot{#2}}{??}}}}
\nc\?[1]{\mathsf{\MakeUppercase{#1}}}
\nc\ddt{\dfrac{\d}{\d t}}
\nc\DOF{\textsc{dof}\xspace}
\nc\e[1]{\texttt{e#1}} %"exponent" as in 1\e-6
\nc\imag{\,\textup{i}} %imaginary unit
\nc\phip{\.1{\phi}}
\nc\qp{\.1{\_q}}
\nc\xp{\.1{x}}
\nc\s{\text{s}}
\nc\inter{\shortintertext}
\nc\mx[1]{\begin{matrix}#1\end{matrix}}
\nc\bmx[1]{\begin{bmatrix}#1\end{bmatrix}}
\nc\pmx[1]{\begin{pmatrix}#1\end{pmatrix}}
\nc\smallbmx[1]{\left[\begin{smallmatrix}#1\end{smallmatrix}\right]}
\nc\ol[1]{\overline{#1}}
\nc\ul[1]{\underline{#1}}

\nc\cpp{\texttt{C++}\xspace}
\nc\daets{\textsc{Daets}\xspace}
\nc\dymola{\textsc{Dymola}\xspace}
\nc\maplesim{\textsc{MapleSim}\xspace}
\nc\modelica{\textsc{Modelica}\xspace}
\nc\twentysim{\textsc{20sim}\xspace}
\nc\dof{degrees of freedom\xspace}

\newtheorem{theorem}{Theorem}[section]
\newtheorem{corollary}[theorem]{Corollary}
\newtheorem{lemma}[theorem]{Lemma}
\newtheorem{definition}[theorem]{Definition}
\newtheorem{remark}[theorem]{Remark}
\newtheorem{examp}{Example}
\newenvironment{example}{\begin{examp}\small\sf}{\end{examp}}

\newtheorem{alg}{Algorithm}[section]
\newenvironment{algo}[1]%Argument gives name of algorithm: mandatory at present!
{%
\smallskip%
\begin{alg}[#1]\sf~\\}
{\smallskip
\end{alg}}

\nc\Asmp[1]{{\bf (A#1)}\xspace}
\nc\cpH{compact \pH}
\nc\evl{eigenvalue\xspace}
\nc\evls{eigenvalues\xspace}
\nc\evc{eigenvector\xspace}
\nc\evcs{eigenvectors\xspace}
\nc\Mdl[2]{{\bf M#1-#2}\xspace}
\nc\PD{positive definite\xspace}
\nc\pH{port-Hamiltonian\xspace}
\nc\PH{Port-Hamiltonian\xspace}

\def\<#1,#2>{\langle\,#1,\,#2\,\rangle} %functional-analysis style
\def\(#1,#2){#1\mathop{\cdot}#2} %applied math style, overrides standard TeX \(
\nc\0{\pmb{0}}
\nc\der[2]{#1^{(#2)}} %higher derivatives
\nc\ddbd[3]{\frac{\partial^2#1}{\partial#2\partial#3}} %2nd partial derivative
\nc\T{{\scalebox{.6}{$\top$}}} %transpose operation
\nc\+{^\T} %or {^\T} %transpose including the "superscript"
\nc\diag{\mathop{\textrm{diag}}}
\nc\emp{\emptyset}
\nc\lam{\lambda}
\nc\1{\mathbbm{1}} %vector of all ones
\nc\dcup{\stackrel{\pmb{\cdot}}{\smash\cup}}
\nc\nv{{n}} %no. of vertices (nodes)
\nc\nb{{b}} %no. of branches (edges)
\nc\nalpha{{\ol\alpha}} %negated bitstring
\nc\nbeta{{\ol\beta}}
\nc\nT{{\nv{-}1}}
\nc\nN{{\nb{-}\nv{+}1}}
\rnc\^[1]{\widehat{#1}}
\nc\Ham{\mathcal{H}}
\nc\calT{\mathcal{T}}
\nc\tee{\mathpzc{t}}
\rnc\v{\nu}
\rnc\i{\iota}
\nc\idvD{\^x} %for combined vector [\i_d;\v_D]
\nc\I{\_I}%{\mathbb{I}} %identity mx, since ordinary I is current source
\nc\Tivl{\mathbb{T}} %interval where t lives
\nc\XX{\mathbb{X}} %region where state variable x lives
\rnc\qp{\dot{q}} %has another definition in macros.tex
\nc\tA{\widetilde A}
\nc\bA{\bar A}
\nc\ta{\widetilde a}
\nc\To{{\,\rightarrow\,}}
\nc\zp{\dot{z}}

\rnc\phi{\varphi}

\nc\iStor{\mathcal{S}} %all energy storage
\nc\iPort{\mathcal{P}} %all external ports
\nc\Dir{\mathcal{D}} % Dirac structure for energy routing
\nc\iT{\texttt{T}}
\nc{\iN}{{\iT^*}}
\nc\iTt{\tilde{\texttt{T}}}
\nc{\iNt}{{\iTt^*}}
\nc\init[1]{\alpha_1(#1)}
\nc\ter[1]{\alpha_2(#1)}
\nc\G{G} %{\mathfrak{G}} %the graph

\nc\edges[1]{{\hspace{.2pt}\texttt{#1}}} %hspace keeps subscripts away from previous letter
\nc\Edges[1]{{\hspace{.2pt}\pmb{\edges{#1}}}}
\nc\iE{\edges{E}} %all edges of the graph
\nc\iS{\edges{S}} %subset of edges of the graph
\nc\iCc{\Edges{C}} %"Cc", all capacitors
\nc\ic{\edges{c}} %Cc intersect T
\nc\iC{\edges{C}} %Cc intersect N
\nc\id{\edges{d}} %in ss:resistexplicit, Rr = d union D
\nc\iD{\edges{D}} %
\nc\iLl{\Edges{L}} %"Ll", all inductors
\nc\il{\edges{l}} %Ll intersect T
\nc\iL{\edges{L}} %Ll intersect N
\nc\iRr{\Edges{R}} %"Rr", all resistors
\nc\iGg{\Edges{G}} %"Gg", all conductors
\nc\iDd{\Edges{D}} %"Dd", all dissipative elements
\nc\ir{\edges{r}} %Rr intersect T
\nc\iR{\edges{R}} %Rr intersect N
\nc\ig{\edges{g}} %Gr intersect T
\nc\iG{\edges{G}} %Gr intersect N
\nc\iVv{\Edges{V}} %all voltage sources (\subseteq T)
\nc\iv{\edges{v}} %same as \iVv
\nc\iV{\edges{V}}
\nc\iIi{\Edges{I}} %all current sources (\subseteq N)
\nc\iI{\edges{I}} %same as \iIi
\nc\iXx{\Edges{X}} %"Xx", in Model 2, all edges except external ports
\nc\ix{\edges{x}} %Xx intersect T
\nc\iX{\edges{X}} %Xx intersect N

\nc\ZS{\mathcal{Z}} %state space of energy storage
\nc\Res{\mathcal{R}} %resistive structure
\nc\Flo{\mathpzc{F}} %a space of flows
\nc\Eff{\mathpzc{E}} %a space of efforts
\nc\flo{\mathpzc{f}\kern -1pt} %an individual flow
\nc\eff{\mathpzc{e}\kern 0.5pt} %an individual effort

\rnc\$[2]{F_{\edges{#1}\edges{#2}}} %for the F submatrices, same font as \ic etc
\rnc\s[1]{\mbox{\footnotesize$#1$}}
\nc\st[1]{\mbox{\scriptsize #1}}

\definecolor{pale-gray}{gray}{0.85}

\DeclareMathAlphabet{\mathpzc}{OT1}{pzc}{m}{it}
\nc\cCl{\_C}
\nc\cG{\scalebox{.85}{$\mathpzc{G}$}}
\nc\cIl{\scalebox{1.2}{$\mathpzc{I}$}}
\nc\cLl{\_L}
\nc\cDl{\_D}
\nc\cGl{\_G}
\nc\cV{\scalebox{.85}{$\mathpzc{V}$}}
\nc\cVl{\scalebox{1.2}{$\mathpzc{V}$}}

\nc\cC{\mathcal{C}}%for Jacobian of capacior relation
\nc\cL{\mathcal{L}}%for Jacobian of inductor relation
\nc\cD{\mathcal{D}}%for Jacobian of dissipator relation

\nc\Span     {\mathop{\rm span}\nolimits}
\nc\rank     {\mathop{\rm rank}\nolimits}
\nc\Loop     {\mathop{\rm cycle}\nolimits} %\cycle is a plain TeX command!
\nc\Cutset   {\mathop{\rm cutset}\nolimits}

\definecolor{myred}{rgb}{0.50, 0. 0.}
\definecolor{keywordscolor}{rgb}{0,0,0.6}

\definecolor{solarized@base03}{HTML}{002B36}
\definecolor{solarized@base02}{HTML}{073642}
\definecolor{solarized@base01}{HTML}{586e75}
\definecolor{solarized@base00}{HTML}{657b83}
\definecolor{solarized@base0}{HTML}{839496}
\definecolor{solarized@base1}{HTML}{93a1a1}
\definecolor{solarized@base2}{HTML}{EEE8D5}
\definecolor{solarized@base3}{HTML}{FDF6E3}

\definecolor{solarized@yellow}{HTML}{B58900}
\definecolor{solarized@orange}{HTML}{CB4B16}
\definecolor{solarized@red}{HTML}{DC322F}
\definecolor{solarized@magenta}{HTML}{D33682}
\definecolor{solarized@violet}{HTML}{6C71C4}
\definecolor{solarized@blue}{HTML}{268BD2}
\definecolor{solarized@darkblue}{HTML}{16537E}
\definecolor{solarized@cyan}{HTML}{2AA198}
\definecolor{solarized@green}{HTML}{859900}

 \nc\li[1]{\protect\lstinline[basicstyle=\ttfamily\color{solarized@blue}]{#1}}

\rnc{\tee}{e}
\rnc{\ell}{e^*}

\nc\Nc[1]{{\color{magenta}#1}\xspace}
\nc\Ns[1]{{\color{green}#1}}
\nc\Nsb[1]{\fbox{\Ns{#1}}}
\nc\Ncl[1]{{\sout{\Nc{#1}}}}
\nc\Nr[2]{\Ncl{#1} {\Ns{#2}}}
\nc\Ncs[2]{\Ncl{#1} \Ns{#2}}

\begin {document}

\thispagestyle{empty}
%    ===== Author, Title, Date =====
\date{\today}
\title{An Energy-based, always Index \texorpdfstring{$\le1$}{<=1} and Structurally Amenable Electrical Circuit Model
\footnotemark[1]}
\author{
Nedialko Nedialkov\footnotemark[2]
\and John D. Pryce\footnotemark[3]
\and Lena Scholz\footnotemark[4]
}

\maketitle

{

\footnotetext[1]{LS acknowledges the support of the Berlin Mathematics Research Center MATH+ by Project AA4-5.\\
NN acknowledges the support of the Natural Sciences and Engineering Research Council of Canada (NSERC), FRN  RGPIN-2019-07054.}

\footnotetext[2]{Department of Computing and Software, McMaster University, Hamilton, Ontario, CN, \texttt{nedialk@mcmaster.ca}.}

\footnotetext[3]{School of Mathematics, Cardiff University, Senghennydd Road, Cardiff, UK, \texttt{PryceJD1@cardiff.ac.uk}.}

\footnotetext[4]{Institute for Mathematics, Technische Universit\"at Berlin, Berlin, Germany, \texttt{lena.scholz@tu-berlin.de}.}
}

\noindent
{\bf Abstract:}\\
Combining three themes: port-Hamiltonian energy-based modelling, structural analysis as used in the circuit world, and structural analysis of general differential-algebraic equations, we form a new model for electrical circuits, the compact port-Hamiltonian equations.
They have remarkable simplicity and symmetry, and always have index at most 1 and other good numerical properties.
The method has been implemented in \textsc{Matlab}.
We give proofs and numerical results.

\vspace{1cm}

\vspace{1cm}

\vfill

{\bf Note:} This report is an expanded version of a manuscript submitted to the SIAM Journal on Scientific Computing (SISC).

\vfill

\newpage
\setcounter{page}{1}

\section {Introduction}

{\em Differential-algebraic equation (DAE)}\/ systems are ubiquitous models of the time evolution of scientific and engineering systems, for good reasons.
A complex system is usually multi-physical, and at the mathematical level can be modelled conveniently as various elementary parts interacting with each other subject to basic physical laws, or an idealisation thereof.
For example, a mechanical system treated as ideal rigid bodies, springs and dampers interacting subject to Newton's laws and the constraints of joints, sliding contact etc.
Or an electrical circuit treated as ideal resistor, capacitor and inductor  elements interacting subject to Kirchhoff's laws and the constraints of the elements' constitutive relations.
Or an electro-mechanical system mixing the two.
Whatever the physical domains, the complete model becomes for some $n$ a set of $n$ equations $f_i=0$ in $n$ variables $x_j$ and some time derivatives of these---that is, a DAE.%
\footnote{This includes the possibility of a hybrid system that can switch between different such DAE representations.}
This paradigm underlies the object-oriented \modelica language, and popular environments such as \dymola, \maplesim and \twentysim.
Continuum models, described by PDEs, often become DAEs after discretisation.
As energy is the main reference quantity for the coupling of different subsystems in multi-physical systems that behave as an energy-transferring network, {\em energy-based network modelling} is beneficial for these kinds of systems.

Successful numerics must find the DAE's {\em causality}---which variables and derivatives influence which in the equations.
An ODE's causality is obvious (the $\xp_j$ are explicit functions of the $x_j$); a DAE's usually not, owing to hidden dependencies.
How far it differs from an ODE's, hence how hard numerical solution is (assuming sufficient smoothness etc.), is measured by the DAE's index.
We use the {\em differentiation index} \cite[Def.~4]{CamG95}: the most times any $f_i=0$ needs to be differentiated with respect to $t$ and the result adjoined to the original equations, so that the enlarged set can be solved to give an ODE in the $x_j$.
One can solve DAEs of index $\le1$ as is, by backward differentiation formula (BDF) codes such as DASSL, SUNDIALS and \matlab's \li{ode15i}; those of index $>1$ are harder and typically need pre-treatment by {\em index reduction} techniques.
\smallskip

In this wide field we contribute to one corner---electrical circuits---by combining three themes,
\begin{compactenum}[(1)]
  \item \pH energy-based modelling, see e.g.~\cite{vdSch13,vdSchJ14},
  \item \SA as used in the circuit world ({\em circuit-SA}, for short), see e.g.~\cite{BarSS11,EstT00,Ria08},
  \item \SA of general DAEs ({\em DAE-SA}, for short) see e.g.~\cite{Pan88,Pry01},
\end{compactenum}
in a way that seems not to have been done before.
The resulting {\em compact \pH (CpH)} circuit DAE is remarkably simple, and always has index $\le1$ and other good numerical properties.
In this introduction we first say enough about themes (1) and (2) to let us describe the CpH model and state the main \thref{mainthm}.
Theme (3) is about the method of proof, not the model, so we then give enough detail of it to outline the proof of \thref{mainthm}.

\paragraph{The \pH approach}
We use the theory of dissipative {\em \pH DAEs (pHDAEs)}, which generalises classical \pH system theory.
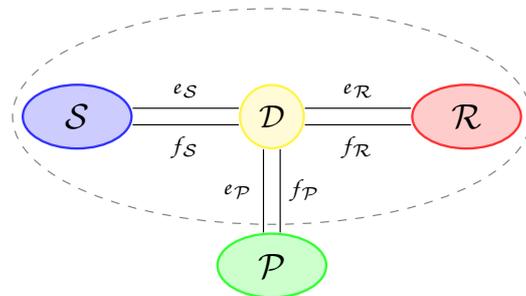
\begin{wrapfigure}[11]{r}{0.5\TW}
\tikzstyle{storage}=[ellipse, thick, minimum width=1.3cm, align=center, draw=blue!80, fill=blue!20]
% The dissipation part is represented by a red ellipse.
\tikzstyle{dissipation}=[ellipse, thick, minimum width=1.3cm, align=center, draw=red!80, fill=red!20]                                    
% The Dirac structure is represented by a yellow circle.
\tikzstyle{dirac}=[circle, thick, minimum size=0.3cm, draw=yellow!80, fill=yellow!20]
%External ports are represented by invisible circle.
\tikzstyle{ports}=[ellipse, thick, minimum width=1.3cm, align=center, draw=green!80, fill=green!20]
\scalebox{1.1}{
\begin{tikzpicture}[>=latex, every node/.append style={midway},ampersand replacement=\&]  
% The elements are placed using a matrix:
  \matrix[row sep=1cm,column sep=1cm]{
         \node (S) [storage] {\large $\iStor$}; \&
         \node (D) [dirac] {$\Dir$};       \&
         \node (R) [dissipation] {\large $\Res$};  \\
      \& \node (P) [ports] {\large $\iPort$}; \& \\
 };
%Draw the interconnections 
\begin{scope}[xshift=40mm, yshift=40mm]
    \draw ([yshift=1mm]S.east) -- ([yshift=1mm]D.west) node[style={font=\scriptsize},above] {$\eff_\iStor$};
    \draw ([yshift=-1mm]S.east) -- ([yshift=-1mm]D.west)
    node[style={font=\scriptsize}, below] {$\flo_\iStor$};
    \draw ([yshift=1mm]D.east) -- ([yshift=1mm]R.west) node[style={font=\scriptsize},above] {$\eff_\Res$};
    \draw ([yshift=-1mm]D.east) -- ([yshift=-1mm]R.west) node[style={font=\scriptsize}, below] {$\flo_\Res$};
    \draw ([xshift=1mm]D.south) -- ([xshift=1mm]P.north) node[style={font=\scriptsize},right] {$\flo_\iPort$};
    \draw ([xshift=-1mm]D.south) -- ([xshift=-1mm]P.north) node[style={font=\scriptsize}, left] {$\eff_\iPort$};
\end{scope}        
%Draw ellipse around the pH system
\draw[gray,dashed] (D) ellipse (3.1cm and 1.3cm);
\end{tikzpicture}
}\vspace{-2ex}
\caption{Port-Hamiltonian system (from \cite{vdSchJ14}). 
\label{fg:pHgenl}}
\end{wrapfigure}
The basic idea, see \fgref{pHgenl}, is to split the system into energy-storing elements grouped in an object $\iStor$, energy-dissipating  elements grouped in an object $\Res$, control inputs and outputs (or external ports) grouped in $\iPort$, and energy-routing elements denoted by $\Dir$ connecting the parts.
The energy storing part is represented by a {\em Hamiltonian} $\Ham$ giving total system energy.
State variables are chosen so that $\Ham$ is an {\em algebraic} function of them.
The energy-routing structure $\Dir$ is linked to other parts via pairs $(\flo,\eff)$ of vectors of flow and effort variables called {\em ports}.
Whatever the physical domain, the product of $\flo$ and $\eff$ (computationally a dot product $\(\flo,\eff)$) has the dimension of power.
One can regard $\Dir$ as a set of equations encoding the network topology; the fact that $\Dir$ conserves energy is encoded in it being a {\em Dirac structure}, see \apref{pHdyn}.
For a circuit, $\Dir$ is given by Kirchhoff's laws.
Such energy-based modelling has many advantages, e.g.: it accounts for the physical interpretation of its variables; it is well suited for modular network modelling since a power-conserving interconnection of \pH systems is again \pH; \pH models are numerically robust with respect to disturbances; they allow for robust simulation and control methods, and easy model refinement (encoded in the structure is the preservation of energy, passivity or stability).
For more details see \cite{BeaMXZ18,MehM19,vdSch04,vdSch13,vdSchJ14}.
The \pH approach makes clear that it is flexible which parts of a system belong to the external part $\iPort$ and which to the dissipative part $\Res$.
This choice is part of the modelling paradigm used, or relates to the goals that are envisaged.

\paragraph{Circuit \SA}
The circuit is treated as a connected undirected graph $G$ with 2-terminal components (resistors/conductors, capacitors, inductors, voltage and current sources) placed on edges, joined at nodes.
Edges carry orientations as a current and voltage housekeeping convention; reversing them does not change physical behaviour.
Kirchhoff's laws say the suitably oriented sum of voltage drops round any cycle of $G$ is zero (KVL) and sum of currents across any cutset of $G$ is zero (KCL).
(Spanning) trees in $G$ are key to circuit analysis.
Each tree specifies a set of a {\em fundamental cycles} and a set of {\em fundamental cutsets}.
It suffices to impose the KVL equations for a tree's fundamental cycles and the KCL equations for a tree's fundamental cutsets, since by linear combination these generate all possible KVL and KCL equations: see e.g.~\cite{BalB69,biggs1993algebraic,ChuDK87,SesR61}.

\paragraph{The CpH DAE}
Some physical and modelling assumptions are required, formalised by \Asmp1 to \Asmp6 in \ssref{assump}. 
That the circuit is {\em well posed}, \Asmp1--\Asmp2, is needed for
ensuring that one can choose a {\em normal tree}: \Asmp3 and \dfref{normalT}. The {\em constitutive relations} that define element behaviour must be {\em passive}: elements cannot create energy from nothing.
Mathematically this relates to positive definiteness of certain Hessians and Jacobians, \Asmp4--\Asmp6.
\scref{CpHmodels} describes three versions of the CpH approach (``models'') differing in the generality of the dissipative structures and external ports they handle. 
The simplest of them---Model 2 in \ssref{model2}---can now be described in a nutshell as follows:
\begin{quote}
Choose a DAE variable for each non-source circuit element as: charge $q$ for a capacitor; flux linkage $\phi$ for an inductor; voltage $\v$ or current $\i$ at will for a resistor.
For a voltage source, take voltage as a control input and current as output; for a current source, vice versa.
Choose an arbitrary normal tree $\iT$.
Then the CpH DAE is formed by the constitutive relations, together with the KVL equations round each fundamental cycle of $\iT$ and the KCL equations across each fundamental cutset.
\end{quote}

We write the KCL and KVL equations as $\i_\iT = F^\T\i_\iN$ and $\v_\iN = -F\v_\iT$, see \rf{KL}, where $\i_\iT,\i_\iN$ denote the vectors of currents on tree and cotree edges respectively, similarly $\v_\iT,\v_\iN$, and $F$ is the {\em \ccmatrix} (also {\em cycle-cutset matrix}) of the tree $\iT$ (used already in the 1950s by Bashkow \cite{Bashkow57} based on Kron's tensor analysis for electrical circuits \cite{Kro39}, see also \cite{Ama62,Bra62,Har63}).
The result describes the circuit by a first-order DAE of control input-output form 
\begin{subequations}
\begin{align}
  0 &= f(t,x,\xp,u(t)) &\text{\em (state equation)},\label{eq:IOforma}\\
  y &= g(t,x,\xp,u(t)) &\text{\em (output equation)}.\label{eq:IOformb}
\end{align}
\end{subequations}
Here vectors $x,f$ have size equal to the number of non-source edges of $G$, and $x$ consists of the $q,\phi,\v,\i$ variables of these edges;
while $u,y,g$ have size equal to the number of source edges, and $u(t)$ holds the voltage sources $V(t)$ and current sources $I(t)$ of these edges, which are assumed to be given functions; $y$ holds the corresponding currents on $V$ edges and voltages on $I$ edges.

\begin{theorem}\label{th:mainthm}
If assumptions \Asmp1--\Asmp6 hold, the system \rf{IOforma,IOformb} is SA-amenable and \pH.
For given input function $u(t)$ the state equation \rf{IOforma} is a DAE of index $\le1$.
\end{theorem}

\paragraph{DAE \SA}
A DAE is {\em SA-amenable} if its causality can be found from the pattern of what variables and derivatives occur in what equations.
This can be done by the Pantelides graph method \cite{Pan88}, %Mattsson \& S\"oderlind dummy derivatives \cite{MatS93}, 
or the Pryce \Smethod \cite{Pry01} used here.
For first-order DAEs such as \rf{IOforma} they are equivalent.
In the \Smethod the DAE is SA-amenable if and only if a certain matrix, the {\em \sysJ} $\_J$, is nonsingular at some sample point.
This forms a {\em success check} of the DAE-SA, after which the method can compute a {\em structural index}.
This can overestimate the differentiation index, but here we prove separately that for a CpH DAE structural and differentiation index are always equal.
Examples of SA-unamenable DAEs can be found in \cite{ReiMB00,Sch18,SchS16b}.  There is much recent work on ways to transform such DAEs to equivalent amenable ones, e.g.\ \cite{IwaT18, TanNP17}.
These are useful approaches but apply to an individual DAE.
As most circuit DAEs are generated automatically within a modelling environment, it is desirable to find an {\em SA-amenable method}, i.e.,\ a modelling approach that generates a DAE known {\sl a priori} to be SA-amenable.

\begin{proof}[Sketch proof of \thref{mainthm}]
Assumptions \Asmp1, \Asmp2 make it possible to choose a normal tree, \Asmp3.
This puts some blocks of zeros in $F$, which improve the sparsity of the \sysJ $\_J$ sufficiently to give it a block-triangular structure, with three square diagonal blocks related to the groups of capacitor edges, inductor edges and dissipator (dissipative element) edges respectively.
Then \Asmp4, \Asmp5, \Asmp6, which relate to the constitutive relations, are just what is needed to prove that each diagonal block is nonsingular, whence $\_J$ is nonsingular, proving SA-amenability.
That index $\le1$ is a by-product.
Now, \Smethod theory shows just which equations need differentiating to convert the DAE to an implicit ODE---namely, those associated with off-tree capacitors and on-tree inductors.
If there are some of these, the DAE has index 1, else index 0.
That \rf{IOforma,IOformb} has \pH structure is a direct result of how it was constructed.
\end{proof}

The remainder of the text is organised as follows. 
In \scref{prelim} we give basic linear algebra results for circuit graphs, define the \ccmatrix, and state our circuit assumptions. 
A contribution we make to \SA theory, in the circuit sense, is to show that a normal tree can be computed by converting the (suitably reordered) circuit incidence matrix to Reduced Row Echelon Form, see \thref{normalT},
as well as by Kruskal's algorithm, see \ssrf{computeTandF}, which tests have shown to be computationally more efficient.
The {\ccmatrix}, and useful counts, e.g.\ of certain cycles and cutsets, appear as a byproduct.
In \scref{CpH} we present the \cpH method for electrical circuits
and the general CpH circuit model in form of a pHDAE. Three specific CpH models are presented in \scref{CpHmodels}. 
We show \thref{mainthm} applies to them.
The relations of the CpH circuit model to other circuit formulations is briefly discussed in \scref{otherforms}. 
In \scref{Numerics} we discuss numerical implementation, with computational examples.
\scref{conclusions} gives some concluding remarks and outlines work in progress.

%%%%%%%%%%%%%%%%%%%%%%%%%%%%%%%%%%%%%%%%%%%%%%%%%%%%%%%%%%%%%%%%%%%%%%%%%%%%%%%%%%
%%% PRELIMINARIES
%%%%%%%%%%%%%%%%%%%%%%%%%%%%%%%%%%%%%%%%%%%%%%%%%%%%%%%%%%%%%%%%%%%%%%%%%%%%%%%%%
\section{Preliminaries}\label{sc:prelim}

\nc\setn[1]{\{\,#1\,\}}
\nc{\cutsym}{\chi}

\subsection{Graphs and their Linear Algebra}

\subsubsection{Graph basics}
Where a term's meaning differs from that which seems commonest in computer science it is {\em *starred} at its first occurrence.
These differences are summarised in \apref{graphterms}.

The topology of an electric circuit can be described as an undirected graph $\G$,
 where a circuit node is represented by a vertex, and if there is a circuit element, or branch, between two nodes, then there is an edge between the corresponding vertices. 
Since two nodes can be connected by more than one element, $\G$ can have multiple edges between two vertices. We assume no element is connected to one node only, so $\G$ does not contain loops (edges from a vertex to itself).
\smallskip

A {\em *graph} $\G$ is an ordered triple $\G =(V, \iE,\alpha)$ with sets $V$ and $\iE$ of vertices and edges, respectively, and a map $\alpha: \iE \to \set{\{u, v\}}{$u,v \in V$, $u\ne v$}$, associating with each edge an unordered pair of distinct vertices. 
If $e\in \iE$ and $\alpha(e) = \{u,v\}$, we say $e$ {\em joins} $u$ and $v$, and call them the {\em ends} of $e$.

Graph $\G'=(V',\iE',\alpha')$ is a {\em subgraph} of $\G=(V, \iE, \alpha)$ if $V'\subseteq V$, $\iE'\subseteq \iE$ and for any $e\in \iE'$, $\alpha'(e)$ is the same unordered pair as $\alpha(e)$.
That is, $\alpha'$ is the restriction of $\alpha$ to $\iE'$, written as $\alpha|_{\iE'}$.

For any $\iE'\subseteq\iE$, we denote by $\G[\iE']$ the {\em *induced subgraph} $(V, \iE', \alpha|_{\iE'})$, that is, with a subset of edges but keeping all the vertices.
Thus $\G[\iE']$ can have isolated vertices.
A {\em trail } is a sequence of alternating vertices and distinct edges $(v_0, e_1, v_1,\dots , v_{l-1}, e_l, v_l)$, where edge $e_i$ has ends $v_{i-1}$ and $v_{i}$, $1\le i\le l$.
A {\em path} is a trail with no repeated vertices. We say the path is {\em between} $v_{0}$ and $v_{l}$.
A graph is {\em connected} if for every pair of distinct vertices there is a path between them. 
A {\em cycle} is a trail with no repeated vertices except the first and last vertices.
Only the cyclic order matters, e.g. $(v_0,e_{1},v_{1},\dots, v_{l-1}, e_l, v_{l}=v_{0},e_{1},v_{1})$ is the same cycle as $(v_0,e_{1},v_{1},\dots, v_{l-1}, e_l, v_{0})$.
We shall also use the term cycle to mean just the set of edges $\{e_1, \dots , e_l\}$.
A graph without cycles is {\em acyclic} and we call a set of edges $\iE'$ acyclic if $\G[\iE']$ is acyclic.
A {\em cut} of $\G$ is a partition of $V$ into two disjoint nonempty subsets, which we call the two {\em sides} of the cut.
The associated {\em cutset} is the set of edges with ends in both sides.
If $\G$ is connected, each cutset is nonempty.
 
A {\em *spanning} subgraph of (a necessarily connected) $\G$ is a connected subgraph that contains all the vertices of $\G$.
We call a set $\iE'\subseteq \iE$ of edges spanning if $\G[\iE']$ is spanning, equivalently if it is connected.
As customary in circuit literature (see e.g. \cite{ChuDK87,Ria08,SesR61}), a {\em *tree} $\iT$ of $\G$ always means a {\em spanning tree}, namely an acyclic spanning subgraph of $\G$.
Its {\em cotree} $\iN$ is the subgraph of $\G$ obtained by removing the edges of $\G$ that are in $\iT$.
We usually identify a tree with its set of edges, and the cotree similarly.
In that sense $\iN$ is the complement of $\iT$, i.e.,
$\iN = \iE\setminus\iT$. 
As in \cite{Ria08} we call tree edges {\em twigs} and cotree edges {\em *links}.

An {\em *orientation} of a graph is an (arbitrary but fixed) assignment of an orientation to each edge. That is, an edge $e$ with ends $u$ and $v$ is oriented from $u$ to $v$ or from $v$ to $u$, denoted here by $u\To v$ or $v\To u$, respectively. 
In $u\To v$, we refer to $u$ as {\em start} vertex and to $v$ as {\em end} vertex; similarly for $v\To u$. 
Assume an orientation of $\G$.
In a cycle $(v_0, e_1, v_1,\dots , v_{l-1}, e_l, v_l)$, edges $e_{i}$ and $e_{j}$, $i\ne
 j$, have the {\em same orientation}, if $v_{i-1}\To v_{i}$ and $v_{j-1}\To v_{j}$ or $v_{i}\To v_{i-1}$ and $v_{j}\To v_{j-1}$, and {\em opposite orientation} otherwise. 
In a cutset, two edges have the {\em same orientation}, if their start vertices are in the same side, and {\em opposite orientation} otherwise.

\subsubsection{Incidence matrix}
For a graph $\G$ defined as above, denote $\nv = |V| $ and $\nb = |\iE|$. 
We label the vertices $i=1,\dots,\nv$ and the edges $j=1,\dots,\nb$. By abuse of terminology, we also identify ``edge'' with ``number in $1,\dots,\nb$'' and ``vertex'' with ``number in $1,\dots,\nv$''.
Given an orientation of $\G$, we represent $\G$ by the (oriented) {\em incidence matrix} $A\in\R^{\nv\x\nb}$ with entries
\begin{align*}
 a_{ij} =\begin{cases}
 \phantom{-} 1 &\text{if vertex $i$ is start of edge $j$},\\
 -1 &\text{if vertex $i$ is end of edge $j$, and} \\
 \phantom{-} 0 &\text{otherwise}.
 \end{cases}
\end{align*}
Thus if the $j$th edge is $i\To i'$, then the corresponding column $a_j$ of $A$ has $+1$ in the start vertex position $a_{ij}$, and $-1$ in the end vertex position $a_{i'j}$, and zeros elsewhere.

The rows of $A$ sum to 0, so $A^\T\1 {=} 0$ where $\1 = (1,\dots, 1)^\T\in\R^\nv$.
That is, $\1$ is in $A^\T$'s null space:
\begin{align*}
 \1\in\ker( A^\T) = \set{x\in\R^\nv}{$A^\T x = 0$}.
\end{align*}

\begin{theorem}[Basic relations of graph topology to linear algebra, see e.g.\ {\cite[p. 199]{Ria08}}]\label{th:A1}~
 \begin{compactenum}[(i)]
  \item The following are equivalent:
  \begin{compactitem}
   \item $\G$ is connected;
   \item $\ker(A^\T) = \Span\{\1\}$, i.e.\ is just the scalar multiples of $\1$;
   \item $A$  
   has rank $\nT$.
  \end{compactitem}
  \item A set of edges is acyclic if and only if the corresponding columns of $A$ are linearly independent, and spanning if and only if they span the column space of $A$. 
 \end{compactenum}
\end{theorem}

\subsubsection{Trees, fundamental cutsets, and fundamental cycles}

Now assume $\G$ is connected. It has at least one tree. 
Identify an edge $j$ with the corresponding column $a_j$ of $A$.
Then by \thref{A1} a tree $\iT$ of $\G$ ``is'' a set of columns of $A$ that forms a basis of the column space.
A basic result of Linear Algebra is that if in a finite-dimensional linear space a linearly independent set $S$ is contained in a spanning set $S'$, one can find a basis $B$ between them: $S\subseteq B\subseteq S'$.
In graph terms by \thref{A1}, if an acyclic set of edges $\iS$ is contained in a spanning set of edges $\iS'$, then there exists a tree $\iT$ between them: $\iS\subseteq \iT\subseteq \iS'$.

Removing any twig from $\iT$ results in a subgraph that is not connected, in fact has exactly two components. (A component might be a single isolated vertex.) 
Their vertex sets are a partition of $\G$'s vertices, thus a cut.
The associated cutset is the {\em fundamental cutset} defined by the twig.

Adding any link to $\iT$ results in a subgraph of $\G$ that contains exactly one cycle, called the {\em fundamental cycle} defined by the link. 
Thus for each tree, there are exactly $\nT$ fundamental cutsets and 
$\nN$ fundamental cycles (see e.g.~\cite[p.~706]{ChuDK87}).
For each $\tee\in\iT$ and $\ell\in\iN$, define 
\begin{align*} 
\begin{aligned}
 \Cutset(\tee) &= \{\text{all edges in the cutset defined by twig $\tee$, excluding $\tee$}\}\quad\text{and}\\
  \Loop(\ell) &= \{\text{all edges in the cycle defined by link $\ell$, excluding $\ell$}\}. 
\end{aligned} 
\end{align*}

Here, $\Cutset(\tee)$ is the links whose fundamental cycles contain $\tee$, and $ \Loop(\ell)$ is the twigs whose fundamental cutsets contains $\ell$
 (see e.g. \cite[p. 204]{gross2005graph}), i.e., 
\begin{align*}
\Cutset(e) &= \{ e^*\in \iT^* \mid e \in \Loop(e^*) \},\\
\Loop(e^*) &= \{ e\in \iT \mid e^* \in \Cutset(e) \}
\end{align*}
and hence 
$e \in \Loop(e^*)\iff e^* \in \Cutset(e)$.

{\bf Notation.}
We use $\edges{s}, \edges{S}, \Edges{S}$ to denote sets of edges, i.e., subsets of $\iE=\{1,\ldots,\nb\}$.
For such a set $\iS$, if $z=[z_i]\in\mathbb{R}^{\nb}$ and $A=[a_{ij}]\in\mathbb{R}^{\nv\times\nb}$, then $z_\iS=[z_i]_{i\in\iS}$ and $A_\iS=[a_{ij}]_{i=1:\nv,j\in\iS}$  denote the sub-vector and sub-matrix corresponding to $\iS$, respectively.
The term $\iS$-cycle shall mean a cycle with edges contained in $\iS$; similarly an $\iS$-cutset.

\begin{theorem}\label{th:A3}\cite[p.~199]{Ria08}
  Let $\iS$ be any set of edges of a connected $\G$ with incidence matrix $A$.
  \begin{compactenum}[(i)]
   \item $\G$ contains no $\iS$-cycles if and only if $A_{\iS}$ has full column rank, i.e., $\ker{A_\iS}=\{0\}$.
   \item $\G$ contains no $\iS$-cutsets if and only if $A_{\iE\setminus \iS}$ has row rank $\nv-1$, i.e., $\ker{A^\T_{\iE\setminus \iS}}=\Span\{\1\}$.
  \end{compactenum}
  \end{theorem}

\subsubsection{The \ccmatrix and Kirchhoff's laws}%\label{ss:Fmatrix}

By \thref{A1}, for any tree $\iT$, the columns of sub-matrix $A_\iT$ (the $\iT$-columns, for short) are a basis of $A$'s column space.
The \ccmatrix $F$ holds the unique coefficients that express the remaining columns of $A$ as linear combinations of the $\iT$-columns.
Enumerating the twigs as $j_1,\ldots,j_\nT$ and the links as $i_1,\ldots,i_\nN$, the matrix $F=(f_{rs})\in\R^{(\nN)\x(\nT)}$ has
\begin{align}\label{eq:ai}
  a_{i_r} = -\sum_{s=1}^\nT f_{rs} a_{j_s} \quad\text{for $r=1,\ldots,(\nN)$},
\end{align}
where $a_k$ is the $k$th column of $A$. The sign is chosen to make the same $F$ as in \cite[\S6.1.2]{Ria08}.
The following facts, though not explicitly stated, can be extracted from results in \cite[\S5.1]{Ria08}.
\begin{theorem}%\label{th:F}
Given a tree $\iT$ of a connected graph $\G$ with cotree $\iN$, arbitrary orientation of $\G$,  and matrix $F$ as above:
\begin{compactenum}[(i)]
\item The nonzeros of $F$ are all $\pm1$, and define all fundamental cycles and cutsets given by $\iT$ and $\iN$. Namely, for any link $i_r\in\iN$ and any twig $j_s\in\iT$ we have
  \begin{align*}
    f_{rs}=\pm1 \iff j_s\in\Loop(i_r) \iff i_r\in\Cutset(j_s).
  \end{align*}
If $f_{rs}=1$ then $i_r,j_s$ have the same orientation around the cycle $\{i_r\}\cup\Loop(i_r)$, and the opposite orientation across the cutset $\{j_s\}\cup\Cutset(j_s)$.
If $f_{rs}=-1$, this is reversed. 
\item Kirchhoff's current and voltage laws can be written as the $\nT$ equations $\i_\iT = F^\T\i_\iN$ and the $\nN$ equations $\v_\iN = -F\v_\iT$, equivalently
  \begin{align}\label{eq:KL}
    0  = \bmx{f_\iT\\ f_\iN} =  \bmx{\i_\iT - F^\T\i_\iN\\ \v_\iN + F\v_\iT} = 
    \bmx{\i_\iT\\ \v_\iN} - \bmx{0 &F^\T\\ -F& 0} \bmx{\v_\iT\\ \i_\iN}.
  \end{align}
\end{compactenum}
\end{theorem}

\begin{example}\label{ex:ex1}
As a running example, we consider the circuit in \fgref{CpH8by5ex0} consisting of a resistor, a conductor, two capacitors and two inductors, as well as a voltage and a current source. %
The  corresponding directed graph with five vertices $V=\{1,\dots,5\}$ and eight edges $\edges{E}=\{1,\dots,8\}$ is given in \fgref{CpH8by5ex_graph0}.
To distinguish the different types of edges we identify $\edges{E}$ with $\edges{E}=\{\edges{V},\edges{C}_1,\edges{C}_2,\edges{G},\edges{R},\edges{L}_1,\edges{L}_2,\edges{I}\}$.
\begin{figure}[ht]
\begin{minipage}[t]{0.48\textwidth}
\begin{center}
\begin{circuitikz}[scale=0.8]
 \draw
 (0,0) -- (0,0.15) to [C, l_=$C_2$,*-o] (0,1.85) -- (0,3) to
 [generic, l_=$R$,*-o] (0,5) -- (0,5) -- (5,5) -- (5,5) to [L, l=$L_2$, *-o] (5,3) -- (5,1.85) to [L, l_=$L_1$,*-o] (5,0.15) -- (5,0) -- (0,0);
 \draw (0,2.5) -- (0.15,2.5) to [V, l=$V$, *-o] (2,2.5) -- (2.5,2.5);
 \draw (5,2.5) to [I, l_=$I$,*-o] (3,2.5) -- (2.5,2.5);
 \draw (2.5,0)--(2.5,0.15) to [C, l_=$C_1$,*-o] (2.5,2);
 \draw (2.5,5) to [generic, l_=$G$,*-o] (2.5,3)--(2.5,2);
\end{circuitikz}
\caption{RLC circuit example \label{fg:CpH8by5ex0}}
\end{center}
\end{minipage}
\begin{minipage}[t]{0.48\textwidth}
 \begin{center}
\nc\linkedge{edge [dashed]}
\nc\treeedge{edge}
\begin{tikzpicture}[->,>=stealth',shorten >=1pt,auto,node distance=2cm,
   thick,main node/.style={circle,fill=blue!10,draw}]
   \node[main node] (1) {1};
   \node[main node] (2) [above of=1] {2};
   \node[main node] (3) [left of=1] {3};
   \node[main node] (4) [below of=1] {4};
   \node[main node] (5) [right of=1] {5};
   \path[every node/.style={font=\sffamily\small}]
      (1) \treeedge node [below] {$\edges{V}$}  node [above] {\st1} (3)
      (1) \linkedge node [right] {$\edges{G}$} node [left] {\st4} (2)
      (2) \treeedge node [right] {$\edges{R}$} node [left] {\st5} (3)
      (3) \linkedge node [right] {$\edges{C}_2$} node [left] {\st3} (4)
      (1) \treeedge node [right] {$\edges{C}_1$} node [left] {\st2} (4)
      (5) \linkedge node [right] {$\edges{L}_2$} node [left] {\st7} (2)
      (4) \treeedge node [right] {$\edges{L}_1$} node [left] {\st6} (5)
      (1) \linkedge node [below] {$\edges{I}$}   node [above] {\st8} (5);
\end{tikzpicture}
\caption{Graph with edge orientations. Tree edges solid, cotree edges dashed.}\label{fg:CpH8by5ex_graph0}
\end{center}
\end{minipage}
\end{figure}
We choose the tree $\iT=\{\edges{V},\edges{C}_1,\edges{R},\edges{L}_1\}$ with cotree $\iN=\{\edges{C}_2,\edges{G},\edges{L}_2,\edges{I}\}$ (depicted in \fgref{CpH8by5ex_graph0} with solid and dashed lines).
The incidence matrix and  the \ccmatrix  are shown in \rf{ex0}.
In $F$, for example,
the entry $(\edges{C}_{2}, \edges{V})$ is 1: these two edges have the same orientation on the cycle  $\{ \edges{C}_{1}, \edges{C}_{2},  \edges{V}\} $,
 and opposite orientation in the cutset $\{ \edges{V}, \edges{C}_{2}, \edges{G}, \edges{L}_{2}\}$. 
Entry $(\edges{G},\edges{V})$ is $-1$: $\edges{V}$ and $\edges{G}$ have opposite orientation on the cycle $\{\edges{V}, \edges{G}, \edges{R}\}$ and the same orientation on the cutset 
$\{ \edges{V}, \edges{C}_{2}, \edges{G}, \edges{L}_2\}$.
\begin{align}\label{eq:ex0}
A= 
  \begin{blockarray}{rrrrrrrr}
    \st1 &\st2 &\st3 &\st4 &\st5 &\st6 &\st7 &\st8\\
     \begin{block}{[rrrrrrrr]}
     1&1&0& 1&0&0&0&1\\
     0&0&0&-1&1&0&-1&0\\
     -1&0&1&0&-1&0&0&0\\
     0&-1&-1&0&0&1&0&0\\
     0&0&0&0&0&-1&1&-1\\
     \end{block}
    \s{\iV}&\s{\iC_1} & \s{\iC_2}&\s{\iG}&\s{\iR}& \s{\iL_1}&\s{\iL_2}&\s{\iI}\\
  \end{blockarray},
\qquad
 F= 
  \begin{blockarray}{rrrrrr}
     &\st1&\st2 &\st5&\st6&  \\
    \begin{block}{r[rrrr]r}
    \s{3}  & 1 &-1 & 0 & 0 & \s{\iC_2} \\  
    \s{4}  &-1 & 0 & 1 & 0 & \s{\iG} \\
    \s{7}  &-1 & 1 & 1 & 1 & \s{\iL_2} \\
    \s{8}  & 0 &-1 & 0 &-1 & \s{\iI}   \\
\end{block}  
     &\s{\iV}&\s{\iC_1} &\s{\iR}&\s{\iL_1}  \\
\end{blockarray}.
\end{align}
\end{example}
  
\subsection{Circuit assumptions}\label{ss:assump}

We consider lumped electrical circuits containing dissipative elements (also called {\em dissipators} in the following) such as resistors and conductors, and energy storing elements such as capacitors and inductors, any of which might be nonlinear, as well as voltage sources and current sources.
Coupling among capacitors, among inductors and among dissipators are possible, but not between an element in one of these groups and an element outside the group.
Various named sets of edges will be used.
First, split the edge set $\iE$ as a disjoint union%
\footnote{Meaning ``union of subsets of a given larger set, that happen to be disjoint'', not the category theory meaning.}
\begin{align}
  \iE = &\iVv\dcup\iCc\dcup\iDd\dcup\iLl\dcup\iIi, \label{eq:edgesplit1}
\end{align}
the subsets related to voltage sources, capacitors, dissipators, inductors and current sources. Let
$n_\iVv,\; n_{\iCc},\;n_{\iDd},\; n_{\iLl},\;   n_\iIi$
be the number of elements in them, respectively.
The set of dissipators can be further partioned as $\iDd=\iRr\dcup\iGg$, if we need to distinguish between the set of resistors and the set of conductors.

The total internal energy of the system is the sum of capacitor energy $\Ham_{\iCc}$ (a function of the $n_{\iCc}$-vector $q$ of charges of all capacitors), and inductor energy $\Ham_{\iLl}$ (a function of the $n_{\iLl}$-vector $\phi$ of fluxes of all inductors)
\begin{align}\label{eq:Ham1}
  \Ham(q,\phi)=\Ham_{\iCc}(q)+\Ham_{\iLl}(\phi).
\end{align}
We assume that $\Ham$ is twice continuously differentiable.
For example, a capacitor with capacitance $C$ leads to  $\Ham_{\iCc}(q)=q^2/2C$; similarly an inductor with inductance $L$ leads to $\Ham_{\iLl}(\phi)=\phi^2/2L$.
\medskip

We make six assumptions; some theory based on the first two is needed to state the third.
\begin{compactenum}
\item[\Asmp1] The circuit has no $\iVv$-cycles, containing voltage source edges only, i.e., $\rank A_\iVv = n_\iVv$.

\item[\Asmp2] The circuit has no $\iIi$-cutsets, containing current source edges only, i.e., $\rank A_{\iE\setminus\iIi} = \nv-1$.
\end{compactenum}

These say the circuit is {\em well-posed}, and are needed because the circuit elements are ideal.
\Asmp1 forbids, e.g., joining voltage sources $V_1(t),V_2(t)$ in parallel: needed because such a circuit is underde\-termined if $V_1(t)=V_2(t)$, or contradictory otherwise.\footnote{Joining two real car batteries in this way is fine---precisely because they are {\em not} ideal sources.}
See \thref{A3} for rank properties.

We now reorder $A$ into element-related matrices
\begin{align}
  A = &[A_\iVv, A_{\iCc},  A_{\iDd},A_{\iLl},  A_\iIi] \label{eq:A}
\end{align}
and define cumulative ranks:
\begin{align}\label{eq:Aranks}
\begin{array}{rlrlr}
r_{\iVv} &:= \rank A_\iVv                          &&= n_\iVv, &\text{ using \Asmp1},\\
r_{\iVv \iCc} &:= \rank [A_\iVv, A_{\iCc}]         &&\leq n_\iVv + n_{\iCc},\\
r_{\iVv \iCc \iDd} &:= \rank [A_\iVv, A_{\iCc}, A_{\iDd}] &&\leq n_\iVv + n_{\iCc} + n_{\iDd},\\
r_{\iVv \iCc \iDd \iLl} &:= \rank [A_\iVv, A_{\iCc}, A_{\iDd}, A_{\iLl}] 
  = \nT &&\leq n_\iVv + n_{\iCc} + n_{\iDd} + n_{\iLl}, &\text{ using \Asmp2}.
\end{array}
\end{align}

\begin{definition}[cf. \cite{Ria08,Sch18}]\label{df:normalT}
A tree $\iT$ of $\G$ is {\em normal} if it contains all $n_\iVv$ edges in $\iVv$, $r_{\iVv \iCc} - r_\iVv$ edges from $\iCc$, $r_{\iVv \iCc \iDd} - r_{\iVv \iCc}$ edges from $\iDd$, $r_{\iVv \iCc \iDd \iLl} - r_{\iVv \iCc \iDd}$ edges from $\iLl$ and no edges from $\iIi$.
\end{definition}

\begin{theorem}\label{th:normalT}
Let $\G$ be the graph of a well-posed circuit, i.e.\ one fulfilling assumptions \Asmp1 and \Asmp2.
Then there exists a normal tree $\iT$ of $\G$. 
\end{theorem}

A graph theoretical proof of \thref{normalT} is given in \cite{Bro63}, without the explicit counts.
For our purpose a proof based on linear algebra is more suitable. 

\begin{proof}
Consider any incidence matrix $A$ of $G$, ordered as in \rf{A}.
By elementary row operations (e.g.\ Gauss Jordan elimination with partial pivoting), $A$ can be transformed to {\em Reduced Row Echelon Form} (RREF) $C$ where for some nonsingular matrix $S\in\R^{\nv\x\nv}$
{\small\[
  C = SA = \left[\begin{array}{p{4mm}p{4mm}p{4mm}|p{4mm}p{4mm}p{4mm}|p{4mm}p{4mm}p{4mm}|p{4mm}|p{4mm}p{4mm}p{4mm}}
1&\mc2{|p{8mm}|}{\quad$\ast$}& 0&\mc2{|p{8mm}|}{ }& 0&\mc2{|p{8mm}|}{ }&   &  0&\mc2{|p{8mm}}{ } \\\cline{1-3}
\mc3{p{12mm}|}{}&1&\mc2{|p{8mm}|}{\quad$\ast$ }&0&\mc2{|p{8mm}|}{ }& &&\mc2{|p{8mm}}{ } \\\cline{4-6}
\mc3{p{12mm}|}{}&\mc3{p{12mm}|}{ }&1&\mc2{|p{8mm}|}{\quad$\ast$ }& & &\mc2{|p{8mm}}{ } \\\cline{7-9}
&0& &\mc3{p{12mm}|}{ }&\mc3{p{15mm}|}{ }&$\ddots$ & \vdots&\mc2{|p{8mm}}{ } \\
\mc3{p{12mm}|}{}&&0&&\mc3{p{12mm}|}{ }& & 0&\mc2{|p{8mm}}{ } \\
\mc3{p{12mm}|}{}&\mc3{p{12mm}|}{ }& &0&& & 1&\mc2{|p{8mm}}{\quad$\ast$}\\\hline
&0&&&0&&&0&&0&&0&
    \end{array}\right].
\]}\noindent
Here $\ast$ denotes some (zero or nonzero) block, and 0 denotes a block of zeros.
Namely, there exists a sequence of natural numbers $j_1,\dots,j_{\nT}$ (the ``pivot'' columns), where $1\leq j_1<\dots<j_{\nT}\leq \nb$, such that  $c_{i,j_i}=1$ for $1\leq i\leq \nT$, and all other entries in column $j_i$ are zero.
The last row of $C$ contains only zeros, as $\rank{A}=\nT$.
Note that no column swaps are done.
By \Asmp1, $A_\iVv$ has full column rank, so all of the first $n_\iVv$ columns of $A$ are pivot columns.
By \Asmp2 the matrix $[A_\iVv\;A_\iCc\;A_\iDd\;A_\iLl]$ has rank $\nv-1$, so $j_{\nT}\leq n_\iVv+n_\iCc+n_\iDd+n_\iLl$, whence no columns of $A_\iIi$ are pivot columns.
Owing to the partial pivoting, $r_{\iVv\iCc} - r_\iVv$ columns from $A_\iCc$ are selected as pivot columns, followed by $r_{\iVv\iCc\iDd} - r_{\iVv\iCc}$ columns from $A_\iDd$: the greatest possible number in each case. 
This leaves the smallest possible number of pivots, namely $r_{\iVv\iCc\iDd\iLl} - r_{\iVv\iCc\iDd}$ to be selected from  $A_\iLl$.
These counts depend only on the ranks defined in \rf{Aranks} and are independent of the ordering of columns within each block (ordering of edges in each subset).
Because $C=SA$ with $S$ nonsingular, the columns $c_j$ of $C$ and $a_j$ of $A$ have the same linear dependences, i.e.\ for any scalars $\lambda_1,\ldots,\lambda_b$ we have $\sum_j\lambda_j c_j =0$ iff $\sum_j\lambda_j a_j =0$.
The pivot columns of $C$ are a basis for its column space, hence the same is true for $A$.
So by \thref{A1}, the edges of $\G$ corresponding to the pivot columns form a tree $\iT$ of $G$ having the properties stated in the Theorem.
\end{proof}

The proof shows that in a normal tree the number of twigs taken from $\iCc$ is the most possible, and the number of twigs taken from $\iLl$ is the fewest possible, in any tree of $\G$.
\smallskip

We can now give the remaining assumptions.
\begin{compactenum}
\item[\Asmp3] The tree $\iT$ is normal.

\item[\Asmp4] The Hessian $\tdbd{^2 \Ham_{\iCc}}{q^2}(q)$ of $\Ham$  in \rf{Ham1} with respect to $q$ is pointwise \PD.

\item[\Asmp5] The Hessian $\tdbd{^2 \Ham_{\iLl}}{\phi^2}(\phi)$ of $\Ham$ in \rf{Ham1} with respect to $\phi$ is pointwise \PD.

\item[\Asmp6] The $n_\iDd$-vectors $\i_\iDd$ and $\v_\iDd$ of currents and voltages on the dissipative edges are related by a continuously differentiable system of $n_\iDd$ equations $r(\i_\iDd,\v_\iDd) = 0$ such that $(\tdbd{ r}{\i_\iDd}, \tdbd{ r}{\v_\iDd})$ is a pointwise positive definite matrix pair (see \dfref{PD}).
\end{compactenum}
\medskip

\Asmp3 is key to obtaining an SA-amenable DAE, while \Asmp4--\Asmp6 state {\em passivity}, i.e., non-source elements can't create energy and only dissipators lose energy.
See \ssrf{PD} for positive definiteness properties in \Asmp4--\Asmp6.
In particular, \Asmp6 relates to {\em strictly local passivity} as defined in \cite[p.\ 267]{Ria08}, and allows full coupling of dissipators which is used e.g. in the modelling of gyrators or ideal transformers (cf. \cite[p.\ 209]{Ria08}).

The matrix $F$ can be directly obtained from the RREF.
With the assumptions and notation of \thref{normalT}, and identifying edges with indices in $\iE = \{1,\ldots,\nb\}$, the tree is the pivot indices $\iT=\{j_1,\ldots,j_\nT\}$ and the cotree is the non-pivot indices $\iN=\iE\setminus\iT = \{i_1,\ldots,i_\nN\}$.
For definiteness, assume both are enumerated in ascending order.
\begin{corollary}\label{co:F}
  $F = -\widetilde{C}^\T$ where $\widetilde{C}$ is the first $\nT$ rows of the matrix $C_\iN$ that comprises the non-pivot columns of $C$.
\end{corollary}
\begin{proof}
It is easily verified that with this definition of $F$, the non-pivot columns of $C$ are given in terms of its (unit vector) pivot columns by
 $c_{i_r} = -\sum_{s=1}^\nT f_{rs} c_{j_s}$ for $r=1,\ldots,\nN$.
Since the columns of $C$ and of $A$ have the same linear dependences, the same holds when $c$'s are replaced by $a$'s, which is equation \rf{ai} that defines $F$.
\end{proof}

\begin{remark}\rm (See \cite[p. 207]{Ria08}.)
The number of {\em link} capacitors $n_\iCc-r_{\iVv \iCc} + r_\iVv$ in a normal tree (i.e., $n_\iC$ in \ssref{model1})  equals the number of linearly independent $\iCc\iVv$-cycles (including in particular $\iCc$-cycles). Analogously, the number $r_{\iVv \iCc \iDd \iLl} - r_{\iVv \iCc \iDd}$  of {\em twig} inductors in a normal tree (i.e., $n_\il$ in \ssref{model1}) equals the number of linearly independent $\iLl\iIi$-cutsets (including $\iLl$-cutsets).
\end{remark}

\section{The \cpH method}\label{sc:CpH}

For the detailed definition and notation with the ingredients (pH1) to (pH4) of a \pH system see \apref{pHdyn}.

The compact \pH (CpH) method consist of the following steps.
First we choose the state variables that describe the energy storage (pH1), and define the dissipative structure (pH2) and the external ports (pH3) of the \pH system.
Then we specify the Dirac structure (pH4) linking the parts together and convert the whole system to a DAE, which brings in extra DAE variables besides the state variables.
In \scref{CpHmodels} we present three circuit models that differ in how the external ports and the dissipative structure are specified.
The energy storing part and the Dirac structure are the same for all three  models.

\subsection{State variables and DAE variables}

The Hamiltonian, giving the total energy stored by the capacitors and inductors, is to be an algebraic function $\Ham(z)$ of the chosen {\em state variables} $z_i$ forming the vector $z$, see \ssrf{pHsystem}.
For an electrical circuit, following \pH practice, e.g. \cite[Appendix B, Table B.1]{vdSchJ14}, we take the $z_i$ as the {\em charges} $q$ on capacitors and the {\em flux linkages} $\phi$ on inductors.
Thus we split $z$ as
\begin{align}\label{eq:pHcirc1}
  z &= (z_\iCc, z_\iLl) = (q,\phi)
\end{align}
and correspondingly the {\em state space} $\ZS$ in (pH1) (see \ssrf{pHsystem}) becomes the product $\ZS = \ZS_\iCc \x \ZS_\iLl$.

This choice of $z$ gives economy in representation because by definition $\dot{q}$ is capacitor currents and $\dot{\phi}$ is inductor voltages.
It can be shown (cf. \cite{GerHRvdS21}) this makes $\tdbd{\Ham}{q}$ the capacitor voltages  and $\tdbd{\Ham}{\phi}$ the inductor currents, even for nonlinear, mutually dependent storage elements.
A capacitor with linear behaviour has $\v_\iCc=\tdbd{\Ham}{q}=q/C$; similarly a linear inductor has $\i_{\iLl}=\tdbd{\Ham}{\phi}=\phi/L$.

A pHDAE is a system of the first order form $f(t,x,\dot x) = 0$ in a vector $x$ of {\em DAE variables}.
In the CpH method, $x$ is all the $z_i$ variables plus some {\em currents or voltages} on non-storage elements---dissipative elements and external ports---depending on the chosen model.

\subsection{Dirac structure}
We write an inner product pairing as $\<a,b>$ when between abstract linear spaces, and $\(a,b)$ when between copies of some $\R^d$. For the latter we do not distinguish between column and row vectors in this section.

With the state space $\ZS$ defined via \rf{pHcirc1},
the Dirac structure linking energy-storage, dissipative structure and external ports is given, cf.\ \rf{pHdyn3}, by
\footnote{By abuse of notation: strictly $T_z\ZS_\iCc$ should read $T_{z_\iCc}\ZS_\iCc$, etc.}
\begin{align*}
  \Dir \;&\subset\; T_z\ZS_\iCc \x T_z^*\ZS_\iCc \,\x\, T_z\ZS_\iLl \x T_z^*\ZS_\iLl \,\x\, \Flo_\Res\x\Eff_\Res \,\x\, \Flo_\iPort\x\Eff_\iPort,
%  \label{eq:pHcirc2}
\end{align*}
where the right side consists of tuples
$(\flo_\iCc, \eff_\iCc,\; \flo_\iLl,\eff_\iLl,\; \flo_\Res,\eff_\Res,\; \flo_\iPort, \eff_\iPort)$.
Reordering, see \rf{pHdyn4}, to
\begin{align}
  &\bigl((\flo_\iCc, \flo_\iLl, \flo_\Res, \flo_\iPort),\; (\eff_\iCc, \eff_\iLl, \eff_\Res, \eff_\iPort)\bigr)
  \label{eq:pHcirc3}
\end{align}
this can be identified with the total flow-effort space $\Flo\x\Eff$.

$\Dir$ is to be defined by Kirchoff's laws for the circuit, but there is a small incompatibility.
Namely, flows $\flo_\iCc,\flo_\Res,\flo_\iPort$ are currents and $\flo_\iLl = -\dot\phi$ is a voltage; efforts $\eff_\iCc,\eff_\Res,\eff_\iPort$ are voltages and $\eff_\iLl = \tdbd{\Ham}{\phi}$ is a current (because of the duality between inductors and capacitors).
To resolve this, a ``swap'' operation is needed.
Let $\ol\Flo$ and $\ol\Eff$ denote the spaces of current vectors $\i$ and voltage vectors $\v$ on the whole circuit, so that Kirchhoff's laws define a linear subspace $\ol\Dir$ of the current-voltage space $\ol\Flo\x\ol\Eff$.
The latter, split up compatibly with \rf{pHcirc3}, consists of tuples
\[
  \bigl((\i_\iCc, \i_\iLl, \i_\Res, \i_\iPort),\; (\v_\iCc, \v_\iLl, \v_\Res, \v_\iPort)\bigr).
\]
Then the operator $K$ (for Kirchhoff) that swaps $\flo_\iLl$ and $\eff_\iLl$, namely
\begin{align*}
  K: (\flo,\eff) &= \bigl((\flo_\iCc, \flo_\iLl, \flo_\Res, \flo_\iPort),\; (\eff_\iCc, \eff_\iLl, \eff_\Res, \eff_\iPort)\bigr) 
  \mapsto \bigl((\flo_\iCc, \eff_\iLl, \flo_\Res, \flo_\iPort),\; (\eff_\iCc, \flo_\iLl, \eff_\Res, \eff_\iPort)\bigr) \\
  &=: \bigl((\i_\iCc, \i_\iLl, \i_\Res, \i_\iPort),\; (\v_\iCc, \v_\iLl, \v_\Res, \v_\iPort)\bigr) = (\i,\v)
\end{align*}
is a linear isomorphism $\Flo\x\Eff\to\ol\Flo\x\ol\Eff$.
Since $\(\flo_\iLl,\eff_\iLl) = \(\eff_\iLl,\flo_\iLl)$, the operator $K$ preserves the inner product:
\begin{align}\label{eq:pHcirc4}
  K(\flo,\eff)  = (\i,\v) \implies \<\flo,\eff> = \(\i,\v).
\end{align}
We define $\ol\Dir$ by \rf{KL}, where $F$ is the \ccmatrix of some tree $\iT$ of the circuit graph 
\begin{align}
  \ol\Dir &= \set{\i\in\ol\Flo}{$\i_\iT = F^\T\i_\iN$} \x \set{\v\in\ol\Eff}{$\v_\iN = -F\v_\iT$}
  \subset\ol\Flo\x\ol\Eff. \label{eq:pHdyn10}
\inter{Then $\Dir$ in \rf{pHdyn3} is the inverse image of $\ol\Dir$ under $K$:}
  \Dir &= K\`\, \ol\Dir = \set{(\flo,\eff)\in\Flo\x\Eff}{$K(\flo,\eff)\in\ol\Dir$}.\nonumber %\label{eq:pHdyn11}
\end{align}
Under our assumptions, $\ol\Dir$ and $\Dir$ are constant, independent of $z$ and $t$, which need not hold for the general definition of pH dynamics in \rf{pHdyn5a}.
\bigskip
\begin{lemma}~%\label{lm:Dirac}~
  \begin{compactenum}[(i)]
  \item $\ol\Dir$ is a Dirac structure ({\em Tellegen's Theorem}).
  \item $\Dir$ is a Dirac structure.
  \end{compactenum}
\end{lemma}
\begin{proof}
For (i), take any $(\i,\v)\in\ol\Dir$ and split $\i,\v$ in tree and cotree parts, then \rf{pHdyn10} gives
\[ \(\i,\v) = \(\i_\iT,\v_\iT) + \(\i_\iN,\v_\iN) = \((F^\T\i_\iN ),\v_\iT) - \(\i_\iN,(F\v_\iT)) = 0. \]
Clearly, $\ol\Dir$ is a $\nb$-dimensional subspace of the $2\nb$-dimensional space $\Flo\x\Eff$, so we have a Dirac structure.
Then (ii) is immediate from \rf{pHcirc4} and $K$ being linear and preserving inner product.
\end{proof}

\subsection{The general CpH circuit model}

\begin{definition}
  A {\em general CpH circuit model} is a pHDAE in which
  \begin{compactenum}[\rm (CpH1)]
  \item the state variables $z$ (describing energy storage) are charge $q$ for capacitors and flux linkage $\phi$ for inductors; and
  \item Kirchhoff's laws are applied in the form \rf{KL} where $F$ is the \ccmatrix of some tree $\iT$ of the circuit graph.
  \end{compactenum}
\end{definition}

Algorithm \ref{al:cphdae} gives an informal algorithm to evaluate the residual vector of the DAE resulting for a general CpH circuit model. 
\begin{algo}{\bf CpH DAE}\label{al:cphdae}
\smallskip
\!\!{\bf Input:} $t,x,\xp$, tree $\iT$.
\begin{compactenum}[Step 1.]
\item \label{alg:step1} Extract $q,\phi$ from $x$, $\qp,\phip$ from $\xp$, and
  use the energy store's constitutive relations to evaluate
  \begin{align}\label{eq:zconstiteqns}
  \begin{aligned}
     \i_\iCc &= \qp,
    &\i_\iLl &= \tdbd{\Ham}{\phi},
    &\v_\iCc &= \tdbd{\Ham}{q},
    &\v_\iLl &= \phip.
  \end{aligned}
  \end{align}

\item  \label{alg:step2} Extract $\v_\iDd, \i_\iDd$ that belong to $\Res$ and $\v_\iPort, \i_\iPort$ that belong to $\iPort$ from $x$ or evaluate them from constitutive relations to get the complete $\v$ and $\i$ vectors as functions of the inputs $(t,x,\xp)$.

\item  \label{alg:step3} 
Split the vectors $\v$ and $\i$ into tree and cotree parts; insert into \rf{KL} to form the residual
\begin{align}\label{eq:CpH}
\bmx{f_\iT\\f_\iN} 
  = \bmx{\i_\iT - F^\T \i_\iN \\ \v_\iN + F \v_\iT}.
\end{align}

\item  \label{alg:step4} Append to \rf{CpH} the residuals of any implicit equations used in specifying $\Res$ and $\iPort$.
\end{compactenum}
{\bf Output:} $f=f(t,x,\xp)$, the vector of residuals resulting from Step 3 and 4.
\end{algo}
This algorithm requires minor adjustment if control input-output form is used, see \ssrf{ctrlIOform}.
The assumption at the start of \ssrf{assump}, of no cross-dependence of dissipators on capacitors or inductors, implies that $\xp$ is not used in Step 2 or Step 4.

Up to this point no specific properties for the tree $\iT$ of the circuit graph have been required.
However, normality of $\iT$ is crucial for proving SA-amenability of the resulting DAE, see \scref{CpHmodels}, thus, from now on we assume  that \Asmp3 holds.
Then we can sub-split the edge set in \rf{edgesplit1} into tree $\iT$ and cotree $\iN$ edges to get
\begin{align}\label{eq:8kinds}
  \iT &= \iv\dcup\ic\dcup\id\dcup \il, \quad
  \iN  = \iC\dcup \iD\dcup \iL\dcup \iI,
\end{align}
where sets with lowercase letters $\iv,\ic,\id,\il$ index voltage sources, capacitors, dissipators and inductors on $\iT$, while sets with uppercase letters $\iC,\iD,\iL,\iI$ index capacitors, dissipators, inductors and current sources on $\iN$.
Correspondingly, we split the vectors $\v_\iT$, $\v_\iN$, $\i_\iT$ and $\i_\iN$ into
\begin{align}\label{eq:vecpartn1}
  \v_\iT = \bmx{\v_\iv\\ \v_\ic\\ \v_\id\\ \v_\il}, \quad
  \v_\iN = \bmx{\v_\iC\\ \v_\iD\\ \v_\iL\\ \v_\iI}, \quad
  \i_\iT = \bmx{\i_\iv\\ \i_\ic\\ \i_\id\\ \i_\il}, \quad
  \i_\iN = \bmx{\i_\iC\\ \i_\iD\\ \i_\iL\\ \i_\iI}.
\end{align}

The matrix $F$ (with a suitable re-ordering of the edges) then takes the block form
\begin{align}\label{eq:blockF}
  F =\
  \begin{array}{c@{}c}
    &\text{tree} \\
  \rotatebox[origin=c]{90}{cotree} &
  \begin{blockarray}{rccccl}
      & \iv &    \ic   &   \id          & \il       \\
  \begin{block}{r[cccc]l}
  \iC & \$\iC\iv & \$\iC\ic &   0 &      0         \\
  \iD & \$\iD\iv& \$\iD\ic       & \$\iD\id &0  \\
  \iL & \$\iL\iv& \$\iL\ic   & \$\iL\id & \$\iL\il \\
  \iI & \$\iI\iv& \$\iI\ic  & \$\iI\id & \$\iI\il  \\
  \end{block}
  \end{blockarray}.
  \end{array}
\end{align}
\vspace{-4ex}

Important for later analysis is that the $\iC\id, \iC\il, \iD\il$ blocks must be zero because any nonzero in them contradicts the assumption $\Asmp3$ that $\iT$ is normal, cf.~\cite[Theorem 3]{Bro63}.
E.g.\ a nonzero in block $\iC\il$ means some twig inductor edge $e\in\il$ is in the fundamental loop of some link capacitor edge $e^*\in\iC$; 
but then $\iT\setminus\{e\}\cup\{e^*\}$ is another tree with fewer inductors (and more capacitors), contradicting the assumption that $\iT$ has the fewest possible inductors.

We illustrate the steps of Algorithm \ref{al:cphdae} with our example.

\begin{example}\label{ex:ex2}
We continue with Example \ref{ex:ex1}. 
Of several possible trees satisfying \Asmp3, we  chose $\iT=\{\edges{V},\edges{C}_1,\edges{R},\edges{L}_1\}$ with cotree $\iN=\{\edges{C}_2,\edges{G},\edges{L}_2,\edges{I}\}$. Each of the sets in \rf{8kinds} has just one member:
 \begin{align*}
\iv=\{\edges{V}\}, \;  \ic=\{\edges{C}_1\}, \; \id=\{\edges{R}\}, \; \il=\{\edges{L}_1\}, \quad
  \iC=\{\edges{C}_2\}, \; \iD=\{\edges{G}\}, \; \iL=\{\edges{L}_2\}, \; 
  \iI=\{\edges{I}\}.
 \end{align*}
Then \rf{vecpartn1} becomes 
 \begin{align*}
   \v_\iT = [\v_\iV; \v_{\iC_1}; \v_{\iR}; \v_{\iL_1}],\quad
   \v_\iN = [\v_{\iC_2}; \v_{\iG}; \v_{\iL_2}; \v_{\iI}],\quad
   \text{and $\i_\iT,\i_\iN$ similarly}.
 \end{align*}
 
Step~\ref{alg:step1}. 
We form  $q = [q_{\iC_1}; q_{\iC_2}]$ and  $\phi = [\phi_{\iL_1}; \phi_{\iL_2}]$.
The Hamiltonian is
\begin{align*}
  \Ham &=\frac{q_{\iC_1}^2}{2C_1}+ \frac{q_{\iC_2}^2}{2C_2}+\frac{\phi_{\iL_1}^2}{2L_1}+\frac{\phi_{\iL_2}^2}{2L_2},
  \end{align*}
 and \rf{zconstiteqns} becomes
 \begin{align}
  \begin{aligned}\label{eq:ex2_1}
    \bmx{\i_{\iC_1}\\ \i_{\iC_2}} &= \bmx{\qp_{\iC_1}\\ \qp_{\iC_2}},
    &\bmx{\i_{\iL_2}\\ \i_{\iL_1}} &= \bmx{ {\phi_{\iL_2}}/L_2\\{\phi_{\iL_1}}/L_1},
    &\bmx{\v_{\iC_1}\\ \v_{\iC_2}} &= \bmx{{q_{\iC_1}}/{C_1}\\ {q_{\iC_2}}/C_2},
    &\bmx{\v_{\iL_2}\\ \v_{\iL_1}} &= \bmx{\phip_{\iL_1}\\ \phip_{\iL_2}}.
  \end{aligned}
  \end{align}
  
Step \ref{alg:step2}. 
The remaining components of $x$ are  $  
 \v_\iR, \v_\iG, \i_\iR, \i_\iG, \v_\iI, \i_\iV, \v_\iV, \i_\iI $, and  we can evaluate $\v_\iV = V(t)$, $\i_\iI = I(t)$. 
  
Step~\ref{alg:step3}.
With the \ccmatrix $F$ in \rf{ex0}, corresponding to $\iT$, Kirchhoff's laws \rf{KL} take the form
\begin{align}\label{eq:ex2_2}
0 & = \bmx{
\i_{\iV}\\
\i_{\iC_1}\\ 
\i_{\iR}\\  
\i_{\iL_1} } 
  -\bmx{
1&-1&-1&0\\
-1&0&1&-1\\
0&1&1&0\\
0&0&1&-1
} 
   \bmx{\i_{\iC_2}\\ \i_{\iG}\\ \i_{\iL_2}\\ \i_\iI },
&
0 & = \bmx{\v_{\iC_2}\\ \v_{\iG}\\ \v_{\iL_2}\\  \v_\iI }
  +\bmx{
1 &-1 & 0  & 0 \\  
-1 & 0 & 1  & 0 \\
-1 & 1 & 1  & 1 \\
0 &-1 & 0  &-1
}
   \bmx{\v_\iV\\ \v_{\iC_1}\\ \v_{\iR}\\  \v_{\iL_1}}.
\end{align}
Using \rf{ex2_1}, \rf{ex2_2} is
\begin{align}\label{eq:ex2_3}
0 & = \bmx{
\i_\iV\\
\qp_{\iC_1}\\ 
\i_{\iR}\\  
\phi_{\iL_1}/L_1 } 
  -\bmx{
 1&-1&-1&0\\
-1&0&1&-1\\
0&1&1&0\\
0&0&1&-1      
       } 
   \bmx{\qp_{\iC_2}\\ \i_{\iG}\\ \phi_{\iL_2}/L_2\\ I(t) },
&
0 & = \bmx{q_{\iC_2}/C_2\\ \v_{\iG}\\ \phip_{\iL_2}\\  \v_\iI }
  +\bmx{
1 &-1 & 0  & 0 \\  
-1 & 0 & 1  & 0 \\
-1 & 1 & 1  & 1 \\
0 &-1 & 0  &-1
      }
   \bmx{V(t)\\ q_{\iC_1}/C_1\\ \v_{\iR}\\  \phip_{\iL_1}}.
\end{align}

Step~\ref{alg:step4}.
We have two types of dissipative elements: a {\em voltage-controlled} conductor $G$ and a {\em current-controlled} resistor $R$ (cf. \ssref{resistexplicit}), which are defined by
\begin{align}\label{eq:ex2_4}
  \i_{\iG}&= G\,\v_{\iG}, \quad  \v_{\iR}=R\,\i_{\iR},
\end{align}
with $G>0$ and $R>0$.
Thus, \rf{ex2_3} and \rf{ex2_4} result in a system of 10 equations in 10 variables
\[
\underbrace{q_{\iC_1}, q_{\iC_2}, \phi_{\iL_1}, \phi_{\iL_2}}_\iStor, \ 
   \underbrace{\v_{\iG}, \i_{\iG},   
   \v_{\iR}, \i_{\iR}}_\Res, \ 
   \underbrace{ \i_\iV, \v_\iI}_\iPort.
\]
\end{example}

In the next sections we present three variations of the general CpH model  that differ in how the external ports and the dissipative structure are specified, and put detail on Steps \ref{alg:step2} and \ref{alg:step4} of \alref{cphdae} for these.

\subsection{Variations of the general CpH model}

\subsubsection{Explicit form of the dissipative structure}\label{ss:resistexplicit}

\nc\rholoc{\rho}
\nc\rhoexpl{\rho}
\nc\xhatR{\^x^{\scalebox{.6}{$\Res$}}}
\nc\Pd{P_\id}

The variables $x$ of the CpH DAE in \alref{cphdae} comprise the vectors $q,\phi$ describing energy storage, summarized in the energy state variables $z=(q,\phi)$,
and additional variables, call them  $\^x$, holding currents or voltages of the non-storage elements. The variables $\^x$ are needed to model the dissipative structure $\Res$ and external ports $\iPort$; their specific forms depend on the chosen model.
If the constitutive relations for the dissipative structure $\Res$, comprising the edge subset $\iDd$ in \rf{edgesplit1}, are implicitly given by the $n_\iDd$ equations $r(\i_\iDd,\v_\iDd)=0$ in $2n_\iDd$ variables $\i_\iDd$, $\v_\iDd$, then all the latter must be in $\^x$, forming 
 a subvector $\xhatR$.
That is, each dissipative element contributes two DAE variables.
However, for each edge $j\in\iDd$ one can choose arbitrarily either $\i_j$ or $\v_j$, and assuming \Asmp6, $r=0$ can be solved locally for the remaining $n_\iDd$ variables.
Namely, let $\iDd=\iRr\dcup\iGg$ where one chooses $\i_j$ for $j\in\iRr$ and $\v_j$ for $j\in\iGg$, then it follows from \lmref{PDMP1}(iv) that the matrix pair
\[
  \left(\bmx{\dbd{r}{\i_\iRr},\dbd{r}{\v_\iGg}}, \bmx{\dbd{r}{\v_\iRr},\dbd{r}{\i_\iGg}}\right)
\]
is pointwise positive definite.
Then from  \lmref{PDMP1}(ii), $\bmx{\dbd{r}{\i_\iRr},\dbd{r}{\v_\iGg}}$ is nonsingular at any point $(\i_\iDd^*,\v_\iDd^*)$ that solves $r(\i_\iDd^*,\v_\iDd^*)=0$.
By the Implicit Function Theorem there exists a continuously differentiable mapping $\rholoc$ such that $r(\i_\iDd,\v_\iDd) = 0$ can be represented locally  as
\begin{align}\label{eq:mixedForm}
 (\v_\iRr,\i_\iGg) &= \rholoc(\i_\iRr,\v_\iGg).
 \end{align} 
The Jacobian of $\rholoc$ is given by
\begin{align}\label{eq:jacR}
\rholoc':=\bmx{\dbd{\rholoc}{\i_\iRr},\dbd{\rholoc}{\v_\iGg}}
  =-\bmx{\dbd{r}{\i_\iRr},\dbd{r}{\v_\iGg}}\` \bmx{\dbd{r}{\v_\iRr},\dbd{r}{\i_\iGg}}.
\end{align}
By \lmref{PDMP1}(ii) this matrix is always positive definite.
Thus, one can represent the dissipative structure by taking $\xhatR$ as $(\i_\iRr,\v_\iGg)$ instead of $(\i_\iDd,\v_\iDd)$---now each dissipative element contributes just one DAE variable.
See \cite[\S6.2.2]{Ria08} for a detailed derivation, also \cite{BalB69,Chu80,GanC78,OhtW69}.
It is useful to represent this in vector form using the permutation matrix $P_\iGg$ such that $P_\iGg[\i_\iDd;\v_\iDd]$ swaps the $\i_j$ and $\v_j$ when $j\in\iGg$.
Then, given a value of $\xhatR$, one reconstructs all of $\i_\iDd$ and $\v_\iDd$ by
\begin{align}\label{eq:mixed2}
 \bmx{\i_\iDd\\ \v_\iDd} = P_\iGg\bmx{\xhatR\\ \rho(\xhatR)}.
\end{align}

If $\iRr=\iDd$ and $\iGg=\emptyset$,
we get the {\em current-controlled representation} $\v_\iDd=\rhoexpl (\i_\iDd)$ with $\xhatR = \i_\iDd$.
If $\iRr=\emptyset$ and $\iGg=\iDd$, we  get the {\em voltage-controlled representation} $\i_\iDd=\rhoexpl (\v_\iDd)$
with $\xhatR = \v_\iDd$. 
In the case of independent linear dissipators this simplifies to $\v_\iDd=R \i_\iDd$ and $\i_\iDd= G \v_\iDd$,
with  positive diagonal {\em resistance}  $R\in\R^{n_\iDd\x n_\iDd}$ and  {\em conductance}   $G\in\R^{n_\iDd\x n_\iDd}$ matrices.

\subsubsection{Control input-output form}\label{ss:ctrlIOform}
Suppose the set $\iPort$ of external port elements consists only of independent sources described by the $n_\iIi$- and $n_\iVv$-vectors $I(t)$ and $V(t)$ of
current and voltage sources, respectively. The source functions $I(t)$ and $V(t)$ are assumed to be continuously differentiable.
We write $\iPort=\iVv\dcup\iIi$ with voltage sources $\iVv$ and current sources $\iIi$.
Then the voltages $\v_\iVv$ of voltage sources and the currents $\i_\iIi$ of current sources are considered as inputs $u$, while the currents $\i_\iVv$ of voltage sources and the voltages $\v_\iIi$ of current sources are considered as outputs $y$, i.e., we have
\begin{align}\label{eq:ctrlIOform1}
  u &=\bmx{\v_\iVv\\\i_\iIi}=\bmx{V(t)\\I(t)}, \qquad y =\bmx{\i_\iVv\\\v_\iIi}.
\end{align}
To make this explicit in the equations, choose the tree $\iT$ and its cotree $\iN$ so that $\iVv\subseteq\iT$ and $\iIi\subseteq\iN$.
This is possible by \Asmp1, \Asmp2; we need not assume here that \Asmp3 holds.
Then in \rf{KL}, $\i_\iVv$ and $\v_\iIi$ are part of $\i_\iT$ and $\v_\iN$, respectively, and occur nowhere else.

Denote by $\iXx = \iCc\dcup\iLl\dcup\iDd$ the set of the remaining edges, and by $\ix= \ic\dcup\il\dcup\id$ and $\iX= \iC\dcup\iL\dcup\iD$ the intersections of $\iXx$ with tree and cotree. Then, $F$ in \rf{blockF} can be represented as
\begin{align*}
  F =\
  \begin{array}{c@{}c}
    &\text{tree} \\
  \rotatebox[origin=c]{90}{cotree} &
  \begin{blockarray}{rccl}
      & \iv &    \ix        \\
  \begin{block}{r[cc]l}
  \iX & \$\iX\iv& \$\iX\ix \\
  \iI & \$\iI\iv& \$\iI\ix \\
  \end{block}
  \end{blockarray}
  \end{array}
\end{align*}
\vspace{-4ex}

and \rf{CpH} takes the form 
\begin{align}
  0 & =\bmx{f_\ix\\ f_\iX} =
    \bmx{\i_\ix - (\$\iX\ix^\T \i_\iX + \$\iI\ix^\T I(t)) \\
         \v_\iX + (\$\iX\ix \v_\ix + \$\iX\iv V(t)) },
   \label{eq:ctrlIOform2}
\inter{together with the output equations}
  y &= \bmx{\i_\iv \\ \v_\iI} =
    \bmx{(\$\iX\iv^\T \i_\iX + \$\iI\iv^\T I(t)) \\ -(\$\iI\ix\v_\ix +  \$\iI\iv V(t))}.
    \label{eq:ctrlIOform3}
\end{align}

Thus, the DAE takes a {\em control input-output form}
\begin{subequations}\label{eq:ctrlIOform4}
\begin{align}
  0 &=f(t,x,\xp,u),\label{eq:ctrlIOform4a}\\
  y &= g(t,x,\xp,u),\label{eq:ctrlIOform4b}
\end{align}
\end{subequations} 
where $x=(q,\phi,\i_\iDd,\v_\iDd)$, $u$, $y$ are given by \rf{ctrlIOform1}, and functions $f:\Tivl\x\R^{N}\x \R^{N}\x\R^{n_\iPort}\to\R^{N}$, $g:\Tivl\x\R^{N}\x\R^{N}\x\R^{n_\iPort}\to\R^{n_\iPort}$ with $N=n_\iCc{+}n_\iLl+2n_\iDd$ or $N=n_\iCc{+}n_\iLl+n_\iDd$ depending on the chosen model (see \scrf{CpHmodels}), and $n_\iPort{=}n_\iVv{+}n_\iIi$ on $\Tivl\subset\R$.

\begin{example}
In \exref{ex2}, the state equation \rf{ctrlIOform4a} is
\begin{small}
\begin{align}\label{eq:ex4_1}
0 & = \bmx{
\qp_{\iC_1}\\ 
\i_{\iR}\\  
\phi_{\iL_1}/L_1 } 
  -\bmx{
       -1 &0&1 & -1  \\
       0 &1& 1 & 0  \\
       0 &0& 1 & -1  } 
   \bmx{\qp_{\iC_2}\\ \i_{\iG}\\ \phi_{\iL_2}/L_2\\ I(t) },
&
0 & = \bmx{q_{\iC_2}/C_2\\ \v_{\iG}\\ \phip_{\iL_2}\\
}
  +\bmx{
      1& -1 & 0  & 0 \\
      -1& 0 & 1  & 0 \\
      -1& 1 & 1  & 1 
        }
   \bmx{V(t)\\ q_{\iC_1}/C_1\\ \v_{\iR}\\  \phip_{\iL_1}}\\
   \label{eq:ex4_2}
    0 &= \i_{\iG}- G\,\v_{\iG}, \quad  0 = \v_{\iR}-R\,\i_{\iR},
\end{align}
\end{small}
and the output equation \rf{ctrlIOform4b} is 
\begin{align*}
\bmx{\i_V \\ \v_I} = 
\bmx{ \qp_{\iC_2}-\i_{\iG}-\phi_{\iL_2}/L_2\\
q_{\iC_1}/C_1 + \phip_{\iL_1}
} .
\end{align*}
This is an  example of Model 1 in \ssref{model1}. 
The equations in \rf{ex4_2} give an explicit representation of dissipating elements as in \rf{mixedForm}.
If we replace in \rf{ex4_1} $\i_{\iG}$ by $G\,\v_{\iG}$ 
and $\v_{\iR}$ by $R\,\i_{\iR}$, 
then  \rf{ex4_1} becomes a system of 6 equations in 6 variables; this is an example of Model 2 in \ssref{model2}.

The tree $\iT$ chosen in \exref{ex2} is normal, thus satisfying $\iVv\subseteq\iT$ and $\iIi\subseteq\iN$.
Note that choosing a non-normal tree, e.g., $\iTt=\{\edges{V},\edges{R},\edges{L}_1,\edges{L}_2\}$ and cotree $\iNt=\{\edges{C}_1,\edges{C}_2,\edges{G},\edges{I}\}$  satisfying $\iVv\subseteq\iTt$ and $\iIi\subseteq\iNt$ results in a system in input-output form that, however, is not necessarily SA-amenable.
Indeed, with $\iTt$, $\iNt$ as above the resulting state equation \rf{ctrlIOform4a} is given by
\begin{small}
\begin{align*}
0 & = \bmx{
\i_{\iR}\\  
\phi_{\iL_1}/L_1 \\
\phi_{\iL_2}/L_2 }
  -\bmx{
       1 &1& 1 & 1  \\
       1 &1& 0 & 0  \\
       1 &1& 0 & 1  } 
   \bmx{\qp_{\iC_1}\\ \qp_{\iC_2}\\ G\v_{\iG}\\ I(t) },
&
0 & = \bmx{q_{\iC_1}/C_1\\q_{\iC_2}/C_2\\ \v_{\iG} }
   +\bmx{
       -1& 1 & 1  & 1 \\
        0& 1 & 1  & 1 \\
       -1& 1 & 0  & 0 
         }
    \bmx{V(t)\\  R\i_{\iR}\\  \phip_{\iL_1}\\  \phip_{\iL_2}},
\end{align*}
\end{small}
which is SA-unamenable.
\end{example}

\section{Specific CpH Circuit Models}\label{sc:CpHmodels}
In the following, we present three specific CpH circuit models that differ in how the dissipative structure (explicit or implicit form) and external ports (independent or controlled sources) are described. 
We will show that all three models are SA-amenable. 

\subsection{CpH Circuit Model 1}\label{ss:model1}

For the CpH Circuit Model 1 we make the following assumptions.
\begin{compactenum}[\bf{M1-}1.]
\item The assumptions \Asmp1--\Asmp6 hold.
\item The external port $\iPort$ comprises only independent current and voltage sources. 
\end{compactenum}
\smallskip
To obtain a well-defined DAE system we append the implicit relation
\begin{align*}
r(\i_\iDd,\v_\iDd) = 0
\end{align*}
as the constitutive relations for the dissipative elements to the system \rf{CpH}.
Due to {\bf M1-2}, we can cast the DAE to input-output form as in
\rf{ctrlIOform4a,ctrlIOform4b}.
For a given input $u$, i.e., given source functions, the state equation  \rf{ctrlIOform4a} is a well-defined DAE system that can be solved for $(t,x,\dot x)$ with a numerical integration scheme (cf. \scref{Numerics}). 
Once a solution for $(t,x,\dot x)$ is obtained, the output $y$ can be computed using \rf{ctrlIOform4b} if required.
Therefore, in the following we restrict our analysis to the CpH DAE
\begin{align}\label{eq:CpH_Model1}
 0=f(t,x,\dot x,u)
\end{align}
that is given by the state equation \rf{ctrlIOform4a} of the control input-output form with a given input $u$.
In \rf{CpH_Model1} the vector $x$ consists of the energy state variables $z=(q,\phi)$ plus the $\^x$ components $\v_\id$, $\v_\iD$, $\i_\id$, $\i_\iD$.
This results in one DAE variable for each capacitor and inductor, and 
two variables for each dissipative element in the circuit, for a total DAE size of
$N = n_\iCc+n_\iLl+ 2 n_\iDd$.

\begin{theorem}\label{th:SAamenable1} 
If the assumptions {\bf M1-1} and {\bf M1-2} hold, then the CpH DAE \rf{CpH_Model1} is SA-amenable.
\end{theorem}

The proof follows standard lines of the \Smethod summarised in \apref{sigmethod}.
We construct the \sigmx of \rf{CpH_Model1}.
We choose suitable offsets, and show in \lmref{nonsingJ}---deploying \Asmp3, \Asmp4, \Asmp5 and the normality of the tree---that the \sysJ defined by these offsets is nonsingular.
The Theorem then follows.
\smallskip

Using \rf{blockF} in \rf{CpH} we get
{\rnc\arraystretch{1.1}
\begin{align}\label{eq:CpH_2}
0  = \left \{
\begin{array}{lc}
 \bmx{ f_\iv\\f_\ic\\ f_\id\\ f_\il} = 
\bmx{\i_\iv\\ \i_\ic\\ \i_\id\\ \i_\il}
    - \bmx{
    \$Cv\+ &\$Dv\+&\$Lv\+ &\$Iv\+ \\
    \$Cc\+ &\$Dc\+&\$Lc\+ &\$Ic\+ \\ 
    0&\$Dd\+&\$Ld\+ &\$Id\+ \\ 
    0&0&\$Ll\+ &\$Il\+ }
       \bmx{\i_\iC\\ \i_\iD\\ \i_\iL\\ \i_\iI}  \\[6ex]
   \bmx{f_\iC\\ f_\iD\\f_\iL\\  f_\iI} = 
      \bmx{\v_\iC\\ \v_\iD\\ \v_\iL\\ \v_\iI}
     + \bmx{
     \$Cv &\$Cc &0    &0   \\
     \$Dv &\$Dc &\$Dd &0   \\
     \$Lv &\$Lc &\$Ld &\$Ll\\  
     \$Iv &\$Ic &\$Id &\$Il  
     }
     \bmx{\v_\iv\\ \v_\ic\\ \v_\id\\ \v_\il} \\[6ex]
     r(\i_\id, \v_\iD, \v_\id, \i_\iD)
     \end{array}
     \right.
 \end{align}
}\noindent
We remove the $f_\iv$ and $f_\iI$ rows in \rf{CpH_2}, as they give 
\rf{ctrlIOform4b}, 
 set $\v_\iv = V(t)$, $ \i_\iI = I(t)$, and using 
 \rf{zconstiteqns} write  the 
 resulting equations constituting \rf{CpH_Model1} as  

{\small
\rnc\arraystretch{1.5}
\begin{align}\label{eq:model1block}
\begin{array}{l@{\;=\;}|c|c|c|c}
  \mc1c{} &\mc1c{\s{q=(q_\iC, q_\ic)}} & \mc1c{\s{\phi=(\phi_\il, \phi_\iL)}} & \mc1c{\s{\^x=(\i_\id, \v_\iD, \v_\id, \i_\iD)}} & \\
  \cline{2-4}
  f_\iC &\dbd{\Ham}{q_\iC}(q_\ic,q_\iC) + \$Cc \dbd{\Ham}{q_\ic}(q_\ic,q_\iC) &&&+ \$Cv V(t) \\
  f_\ic &- \$Cc\+\qp_\iC + \qp_\ic &- \$Lc\+\dbd{\Ham}{\phi_\iL}(\phi_\il,\phi_\iL) &- \$Dc\+\i_\iD & - \$Ic\+I(t) \\ \cline{2-4}
  f_\il &&\dbd{\Ham}{\phi_\il}(\phi_\il,\phi_\iL) - \$Ll\+\dbd{\Ham}{\phi_\iL}(\phi_\il,\phi_\iL) &&- \$Il\+I(t) \\
  f_\iL & \$Lc\dbd{\Ham}{q_\ic}(q_\ic,q_\iC) &+ \$Ll \phip_\il + \phip_\iL &+ \$Ld \v_\id &+ \$Lv V(t) \\
        \cline{2-4}
  f_\id &&- \$Ld\+ \dbd{\Ham}{\phi_\iL}(\phi_\il,\phi_\iL) &+ \i_\id - \$Dd\+\i_\iD &- \$Id\+ I(t) \\
  f_\iD & \$Dc \dbd{\Ham}{q_\ic}(q_\ic,q_\iC) &&+\v_\iD + \$Dd \v_\id &+ \$Dv V(t) \\ 
  r     &&&{r(\i_\id, \v_\iD, \v_\id, \i_\iD)} & \\ \cline{2-4}
\end{array}
\end{align}
}

The corresponding schematic $\Sigma$ matrix is 
{
\rnc\![1]{\colorbox{lightgray!110}{$#1$}}
\rnc\g[1]{\colorbox{lightgray!50}{$#1$}}
\begin{align}\label{eq:model1Sigma} 
  \Sigma \ &=\
\begin{blockarray}{rcccccccclc}
         &  q_\iC     &  q_\ic       &  \phi_\il  &  \phi_\iL    &  \i_\id &  \v_\iD    &  \v_\id  &  \i_\iD  & \s{c_i} &       \\
  \begin{block}{r[cc|cc|cccc]ll}
   f_\iC & \!{\le0}   & \!{\le0}     & -          & -            & -       & -          & -        & -        & \s1     & n_\iC \\
   f_\ic & \!{\le 1}   & \!{\^\I{+}1} &  \g{\le0}    &  \g{\le0}      & -       & -          & -        & \g{\le0}   & \s0     & n_\ic \\ \cline{2-9}        
   f_\il & -          & -            & \!{\le0}   & \!{\le0}     & -       & -          & -        & -        & \s1     & n_\il \\
   f_\iL &  \g{\le0}    &  \g{\le0}      & \!{\le1}   & \!{\^\I{+}1} & -       & -          & \g{\le0}   & -        & \s0     & n_\iL \\ \cline{2-9}        
   f_\id & -          & -            &  \le0      &  \le0        &\!{\^\I} & -          & -        & \!{\le0} & \s1    & n_\id \\
   f_\iD &  \le0      &  \le0        & -          & -            &-        & \!{\^\I}   & \!{\le0} & -        & \s1     & n_\iD \\    
   r     & -          & -            & -          & -            &\!{\le0} & \!{\le0}   & \!{\le0} & \!{\le0} & \s1     & n_\iDd{=}n_\iD{+}n_\id\\
  \end{block}  
  \s{d_j}&  \s1       &  \s1         &  \s1       &  \s1         &  \s1    & \s1  &\s1  & \s1      &          &         &       \\
           &  n_\iC     &  n_\ic       &  n_\il     &  n_\iL       &  n_\id  &  n_\iD     &  n_\id   &  n_\iD   &         &
  \end{blockarray}
\end{align}
}\noindent
All entries $\sij$ are $\ninf$, $0$ or $1$ because at most first derivatives occur in the DAE.
Here $\le 0$ denotes a block with $0$ or $\ninf$ entries; $\le 1$ denotes a block with $1$ or $\ninf$ entries; $-$ is a block of $\ninf$  entries; $\^\I$ is a square block with $0$'s on the main diagonal and $\ninf$ elsewhere; $\^\I+1$ (elementwise addition) is $1$ on the main diagonal and $\ninf$ elsewhere.
We assign provisional offset values $c_i, d_j$. 
Their meaning is blockwise, e.g.\ $c_i=0$ everywhere in the $f_\ic$ and $f_\iL$ rows and $c_i=1$ elsewhere.

\begin{lemma}\label{lm:nonsingJ}
  The \sysJ $\_J$ defined by the provisional offsets is nonsingular.
\end{lemma}

\begin{proof}
$\_J$ is defined by \rf{sysJ} in \apref{sigmethod}.
The blocks with light grey background in \rf{model1Sigma} are those where $d_j-c_i>\sij$, and therefore 
 $\_J_{ij}=0$ in these blocks. Hence  $\_J$ has block-triangular form 
{\rnc\![1]{\colorbox{lightgray}{$#1$}}
\begin{align}\label{eq:Jblocks}
  \_J \ &=\
\begin{blockarray}{rcccl}
          & q_\iC,q_\ic & \phi_\il,\phi_\iL & \i_\id,\v_\iD,\v_\id,\i_\iD  &              \\
  \begin{block}{r[c|c|c]l}
    f_\iC,f_\ic         & \!{\_J_{\cC}}    & 0                             & 0            &\s{n_\iCc=n_\iC{+}n_\ic} \\  \cline{2-4}
    f_\il,f_\iL         & 0                & \!{\_J_{\cL}}                 & 0            &\s{n_\iLl=n_\il{+}n_\iL} \\\cline{2-4}
    f_\id,f_\iD,r       & \ast            & \ast                         &\!{\_J_{\cD}} &\s{2n_\iDd=n_\id{+}n_\iD{+}n_\iDd} \\
  \end{block}
                        & \s{n_\iCc}       & \s{n_\iLl}                    & \s{2n_\iDd}
  \end{blockarray}
\end{align}
}where $\ast$ marks a possibly nonzero block. 
It suffices to show each of $\_J_{\cC},\_J_{\cL},\_J_{\cD}$ is nonsingular.

\para{The $\_J_{\cC}$ block}
Since $d_j-c_i$ equals $0$ in the $f_\iC$ rows and $1$ in the $f_\ic$ rows,  by \rf{sysJ}
\begin{align*}
  \_J_{\cC} &= \bmx{\ds\dbd{f_\iC}{q_\iC}& \ds\dbd{f_\iC}{q_\ic}
  \\[2ex]
  \ds\dbd{f_\ic}{\dot q_\iC}& \ds\dbd{f_\ic}{\dot q_\ic}} 
   = \bmx{\ds\dbd{^2 \Ham}{q_\iC^2} + \$Cc \ddbd{\Ham}{q_\iC}{q_\ic} & \ds\ddbd{\Ham}{q_\ic}{q_\iC} + \$Cc \dbd{^2 \Ham}{q_\ic^2} \\[3ex]
  -\$Cc\+ & \I} \\
 \nonumber
\end{align*}
Here, $\I$ denotes an identity matrix of suitable size.
By \Asmp4,  the matrix 
\begin{align*}
  \cCl  & = \bmx{\ds\dbd{^2 \Ham}{q_\iC^2}  & \ds\ddbd{\Ham}{q_\iC}{q_\ic}  \\[3ex]
 \ds\ddbd{\Ham}{q_\ic}{q_\iC}& \ds\dbd{^2 \Ham}{q_\ic^2} } 
\end{align*}
is PD. Setting  $M =  \cCl $ and $N=-\$Cc\+$ in \lmref{PD2}, $\_J_{\cC} $ is of the form \rf{PD2} and hence   nonsingular.

\para{The $\_J_{\cal L}$ block}
Similarly,  
\begin{align*}
  \_J_{\cal L} &= 
  \bmx{\ds\dbd{f_\il}{\phi_\il}& \ds\dbd{f_\il}{\phi_\iL}\\[2ex] \ds\dbd{f_\iL}{\dot \phi_\il}& \ds\dbd{f_\iL}{\dot \phi_\iL}}
  =\bmx{
  \ds\dbd{^2 \Ham}{\phi_\il^2} -
  	\$Ll\+  \ddbd{\Ham}{\phi_\iL}{\phi_\il} & 
 \ds  \ddbd{\Ham}{\phi_\il}{\phi_\iL} -
  	\$Ll\+  \dbd{^2 \Ham}{\phi_\iL^2}
  \\[3ex]
  \$Ll & \I}.
 \end{align*}
By \Asmp5 the matrix 
\begin{align*}
  \cLl & = \bmx{\ds\dbd{^2 \Ham}{\phi_\il^2}  & \ds\ddbd{\Ham}{\phi_\il}{\phi_\iL}  \\[3ex]
 \ds\ddbd{\Ham}{\phi_\iL}{\phi_\il}& \ds\dbd{^2 \Ham}{\phi_\iL^2} } 
\end{align*}
is PD. 
Setting   $M=\cLl$ and $N=\$Ll$ in \lmref{PD2}, $\_J_{\cal L}$ is of the form \rf{PD2} and hence  nonsingular.

\para{The $\_J_\cD$ block}
We have
\[
  \dbd{(f_\id, f_\iD)}{(\i_\id, \v_\iD, \v_\id, \i_\iD)} = 
  \bmx{\I & 0 & 0    &-\$Dd\+\\
       0 & \I & \$Dd &  0      } =: \bmx{\I&S} \in\R^{n_\iDd\x2n_\iDd},
\]
where $S$ is skew-symmetric.
By \Asmp6, $(H,H') {=} \bigl(\tdbd{r}{(\v_\iDd)}, \tdbd{r}{(\i_\iDd)}\bigr) {=} \bigl(\tdbd{r}{(\v_\id, \v_\iD)}, \tdbd{r}{(\i_\id, \i_\iD)}\bigr)$ is a positive definite matrix pair (PDMP).
By \lmref{PDMP1}(iv), $(\cD,\cD') = \bigl(\tdbd{r}{(\i_\id, \v_\iD)},\linebreak \tdbd{r}{(\v_\id, \i_\iD)}\bigr)$, which is obtained by swapping some columns of $H$ with corresponding columns of $H'$, is also a PDMP.
Thus Jacobian $\_J_\cD$ has the form
\begin{align*}
 \_J_\cD &= 
{\rnc\arraystretch{2}
 \bmx{
\ds\dbd{f_\id}{\i_\id}&\ds\dbd{f_\id}{\v_\iD}& \ds\dbd{f_\id}{\i_\iD}& \ds\dbd{f_\id}{\v_\id}\\
\ds\dbd{f_\iD}{\i_\id}&\ds\dbd{f_\iD}{\v_\iD}&\ds\dbd{f_\iD}{\i_\iD}&\ds\dbd{f_\iD}{\v_\id}\\[1ex]\cline{1-4}
\ds\dbd{r}{\i_\id}&\ds\dbd{r}{\v_\iD}& \ds\dbd{r}{\i_\iD}& \ds\dbd{r}{\v_\id}
}
\mx{n_\id \\n_\iD \\n_\iDd}
 = \bmx{\I &S \\\cline{1-2} \cD &\cD'}\mx{n_\iDd \\n_\iDd}
}
 \end{align*}
satisfying the conditions of \lmref{HK} and hence is nonsingular.
\end{proof}

{\bf Proof of \thref{SAamenable1}.}\\
\lmref{nonsingJ} has shown that the Jacobian $\_J$ defined by the provisional offsets $c_i$, $d_j$ given in \rf{model1Sigma} is nonsingular.
Thus, there exists a  transversal $\calT$  such that $\_J_{ij}\ne0$ for all $(i,j)\in\calT$, which means the DAE \rf{CpH_Model1} is structurally well-posed (has a finite transversal).
Namely $d_j-c_i = \sigma_{ij}$ at all such $(i,j)$, see \rf{sysJ}, and by  \cite[Theorem 3.4(iii)]{Pry01}, $\calT$ is an HVT.
The value of $\calT$ is 
\begin{align*}
 \val\calT = \sum_{(i,j)\in \calT}\sij = \sum_j d_j - \sum_i c_i = n_\ic+n_\iL, 
\end{align*}
which is also the number of degrees of freedom.
Since the offsets $c_i,d_j$ satisfy \rf{offsets} in \apref{sigmethod} they are valid offsets.
Hence, the DAE \rf{CpH_Model1} is SA-amenable.
\qed

\medskip   
\begin{theorem}\label{th:index_pHC1}
Assume that {\bf M1-1} and {\bf M1-2} hold.
If $n_\iC=n_\il=n_\id=n_\iD=0$ the system \rf{CpH_Model1} is an ODE.
Otherwise, \rf{CpH_Model1} is a DAE of index 1.
\end{theorem} 

\begin{proof}
If $n_\iC=n_\il=n_\id=n_\iD=0$ it is immediately clear that the system \rf{CpH_Model1} simplifies to an ODE,
as \rf{model1block} does not contain algebraic equations.
In the other case---at least one of $n_\iC, n_\il, n_\id, n_\iD$ does not vanish---since the DAE is SA-amenable by \thref{SAamenable1} we can apply the {\em standard solution scheme (SSS)}, see \rf{SSS} in \scrf{sigmethod} for solving \rf{CpH_Model1}.
By \thref{Jrank} the Jacobian $\_J_k$ is of full row rank for $k=-1,0$ such that the generally underdetermined system in stage $k=-1$ can be solved for $x$ if $\DOF= n_\ic+n_\iL$ initial conditions are provided. 
In stage $k=0$ we get a nonsingular system for all of $\xp$.
Thus, one differentiation of equations is needed to convert the original DAE to an implicit ODE, and the DAE has index 1.
\end{proof}

\begin{corollary}
The following statements are equivalent:
\begin{compactenum}[(i)]
\item system \rf{CpH_Model1} is an ODE;
\item the circuit contains neither dissipators $\iDd$, ${\iLl\iIi}$-cutsets or ${\iCc\iVv}$-cycles;
\item $A_{\iDd}=\emptyset$ and $\rank[A_{\iVv}\;A_{\iCc}]=\nT$ as well as $\ker[A_{\iVv}\;A_{\iCc}]=\{0\}$.
\end{compactenum}
\end{corollary} 
\begin{proof}
Under \Asmp3 the circuit contains ${\iCc\iVv}$-cycles if and only if $n_\iC>0$ (since otherwise components from $\iC$ could be swapped with other component from the tree contradicting  the optimality of \Asmp3). The same argument holds for ${\iLl\iIi}$-cutsets.
The corresponding rank conditions on $A$ follow from \thref{A3}. 
\end{proof}

\begin{remark}\rm
If one is interested in controlling the circuit, i.e., the input $u$ is treated as an unknown itself, then the whole system \rf{ctrlIOform4} has to be considered.
The same holds if the interconnection of different \pH systems via inputs and outputs is taken into account.
The result in \thref{SAamenable1} also holds for this more general setting, although with  different offsets and index results.
Since we are interested in the simulation (and not control) of the circuit equations, we do not present these more general results here.
\end{remark}

\begin{remark}\rm
The offsets provided in \rf{model1Sigma} are shown to be valid offsets, however they are non-canonical, i.e., not the smallest ones possible (see \scrf{sigmethod}).
Using these offsets from \rf{sindex} we get a structural index $\nu_S$ which corresponds to the results of \thref{index_pHC1}: $\nu_S=0$ if $n_\iC=n_\il=n_\id=n_\iD=0$, and  $\nu_S=1$ otherwise.
If we use the canonical offsets, i.e., $c_i=0$ and $d_j=0$ in the rows/columns belonging to the dissipative structure, we get a structural index which overestimates the index by one.
\end{remark}

\subsection{CpH Circuit Model 2}\label{ss:model2}
For the CpH Circuit Model 2 we make the following assumptions.
\begin{compactenum}[\bf{M2-}1.]
\item The assumptions \Asmp1--\Asmp6 hold.
\item The external port $\iPort$ comprises only independent current and voltage sources. 
\item The dissipative structure is described by a mixed form as in \ssrf{resistexplicit}, \rf{mixedForm}.
\end{compactenum}
\smallskip

As in \ssrf{ctrlIOform}, we cast the DAE into an input-output form and restrict to the CpH DAE
\begin{align}\label{eq:CpH_Model2}
 0=f(t,x,\xp,u),
\end{align}
where $u$ is given input.
Port $\iPort$ contributes no DAE variables, so \alref{cphdae} has $x=(z,\^x)$ with energy state variables $z=(q,\phi)$, dissipative variables $\^x$---one variable for each dissipative edge, for a total DAE size of
\begin{align}\label{eq:CpH_Model2b}
N = n_\iCc+n_\iLl+ n_\iDd = \nb - n_\iPort,
\end{align}
i.e., the total number of branches that are not external ports.

\begin{theorem}
If the assumptions {\bf M2-1}--{\bf M2-3} hold, then the CpH DAE \rf{CpH_Model2} is SA-amenable.
\end{theorem}

\begin{proof}
The proof is similar to that of \thref{SAamenable1}.
The version of \rf{model1block} for Model 2 is
{\small
 \rnc\arraystretch{1.5}
 \begin{align}\label{eq:model2block}
 \begin{array}{l@{\;=\;}|c|c|c|c}
  \mc1c{} &\mc1c{\s{q=(q_\iC, q_\ic)}} & \mc1c{\s{\phi=(\phi_\il, \phi_\iL)}} & \mc1c{\s{\^x=(\i_\id, \v_\iD)}} & \\
 \cline{2-4}
   f_\iC 
         &\dbd{\Ham}{q_\iC}(q_\ic,q_\iC) {+} \$Cc \dbd{\Ham}{q_\ic}(q_\ic,q_\iC) &&&+ \$Cv V(t) \\
   f_\ic 
         &- \$Cc\+\qp_\iC + \qp_\ic &- \$Lc\+\dbd{\Ham}{\phi_\iL}(\phi_\il,\phi_\iL) &- \$Dc\+\i_\iD & - \$Ic\+I(t) \\ \cline{2-4}
   f_\il 
         &&\dbd{\Ham}{\phi_\il}(\phi_\il,\phi_\iL) {-} \$Ll\+\dbd{\Ham}{\phi_\iL}(\phi_\il,\phi_\iL) &&- \$Il\+I(t) \\
   f_\iL 
         & \$Lc\dbd{\Ham}{q_\ic}(q_\ic,q_\iC) &+ \$Ll \phip_\il + \phip_\iL &+ \$Ld \v_\id &+ \$Lv V(t) \\
         \cline{2-4}
   f_\id 
         &&- \$Lr\+ \dbd{\Ham}{\phi_\iL}(\phi_\il,\phi_\iL) &+ \i_\id {-} \$Dd\+\i_\iD &- \$Id\+ I(t) \\
   f_\iD 
         & \$Dc \dbd{\Ham}{q_\ic}(q_\ic,q_\iC) &&+\v_\iD {+} \$Dd \v_\id &+ \$Dv V(t) \\ \cline{2-4}
 \end{array}
 \end{align}
 }\noindent

Here, we use that by \rf{mixed2} $\i_\id, \i_\iD, \v_\id, \v_\iD$ are explicit functions of the $\^x=\^x^\Res$ variables.
This gives the schematic form \rf{model2Sigma} of the \sigmx:
{\rnc\![1]{\colorbox{lightgray!110}{$#1$}}
\rnc\g[1]{\colorbox{lightgray!50}{$#1$}}
\begin{align}\label{eq:model2Sigma}
  \Sigma \ &=\
  \begin{blockarray}{rccccclc}
          &    q_\iC  &    q_\ic   & \phi_\il  &   \phi_\iL &  \^x    & \s{c_i} &       \\ \begin{block}{r[cc|cc|c]lc}
    f_\iC &\!{\le0}   &\!{\le0}    &  -        &  -         &    -    & \s1     & n_\iC \\
    f_\ic &\!{\le1}   &\!{\^\I{+}1}&\g{\le0}     &\g{\le0}      &  \g{\le0} & \s0     & n_\ic \\ \cline{2-6}        
    f_\il &  -        &  -         &\!{\le0}   &\!{\le0}    &    -    & \s1     & n_\il \\
    f_\iL &\g\le0     &\g{\le0}      &\!{\le1}   &\!{\^\I{+}1}&  \g{\le0} & \s0     & n_\iL \\ \cline{2-6}
    f_\id&  -        &  -         &\le0       &\le0        &\!{\le0} & {\s1}   & n_\id \\
    f_\iD &\le0       &\le0        &  -        &  -         &\!{\le0} & {\s1}   & n_\iD \\  \end{block}  
   \s{d_j}&  \s1      &  \s1       &  \s1      &  \s1       & {\s1}   &         &       \\
          &   n_\iC   &  n_\ic     &   n_\il   &   n_\iL    &  n_\iDd &         &
  \end{blockarray}
\end{align}
}\vspace{-4ex}

As in \rf{model1Sigma}, provisional offsets are added, and the blocks with light-grey background are those that cannot contribute nonzeros to the associated \sysJ.
The $\_J_{\cC}$, $\_J_{\cal L}$ blocks are the same as before.
We need to show that the Jacobian 
{ \rnc\arraystretch{1.5}
\[\_J_{\cD}=\bmx{\tdbd{f_\id}{\^x}\\\tdbd{f_\iD}{\^x}}\]}
is nonsingular.
Let $P_1$ be the permutation matrix that swaps $\i_\iD$ and $\v_\iD$ in $[\i_\iDd;\v_\iDd]=[\i_\id;\i_\iD;\v_\id;\v_\iD]$, and let $P=P_1P_\iGg$ where $P_\iGg$ is as in \rf{mixed2}. Denote $S=\smallbmx{0& -\$Dd\+ \\ \$Dd & 0}$.

The part of the $f_\id$, $f_\iD$ equations that depends on $\^x$ is
\begin{align*}
\bmx{\i_\id\\ \v_\iD} + \bmx{0& -\$Dd\+ \\ \$Dd & 0}\bmx{\v_\id\\ \i_\iD}
=[\I,\;S]P_1\, \bmx{\i_\iDd\\ \v_\iDd}
= [\I,\;S]\, P \bmx{\^x\\ \rho(\^x)}. 
\end{align*}
Denote $\rho'=\tdbd{\rho}{\^x}$. Then
\begin{align*}
 \_J_{\cD} 
 &= \dbd{}{\^x} \bmx{f_\id\\f_\iD}
 = [\I,\;S]\; P\; \bmx{\I\\ \rho'},
\end{align*}
where $\rho'$ is PD by assumption.
By \lmref{PDMP2} $(I,\rho'^\T)$ is a NDMP.
Let $[C_1, C_2] = [\I,\rho'^\T]P$. By \lmref{PDMP1}, $(C_1, C_2)$ is a NDMP, $C_1$ and $C_2$ are nonsingular, and
$C^{-1}_2 C_1$ is PD.
Then
\[\_J_{\cD} = [\I, S]\bmx{C^\T_1\\C_2^\T}= C_1^\T + SC_2^\T = (C_1 - C_2S)^\T = [C_2(C^{-1}_2 C_1 - S)]^\T.\]
Since $S$ is skew-symmetric the nonsingularity of $\_J_{\cD}$ follows from
\lmref{PDone}.
\end{proof} 

In a similar way as before we get the following results.
\begin{theorem}
Assume that {\bf M2-1} to {\bf M2-3} hold. If $n_\iC=n_\il=n_\id=n_\iD=0$, the
system \rf{CpH_Model2} is an ODE.
Otherwise, \rf{CpH_Model2} is a DAE of index $1$.
\end{theorem} 

\begin{example}
The running example can be described in the form of Model 2.
Below is shown the signature matrix, alongside the actual equations given by \rf{ex4_1} with \rf{ex4_2} inserted.
 {\rnc\![1]{\colorbox{lightgray}{$#1$}}
 \begin{align}\label{eq:Sigma2}
   \Sigma &=
   \raisebox{-1.2ex}
   {$\begin{blockarray}{rccccccl}
            & q_{\iC_2}    & q_{\iC_1}     & \phi_{\iL_1} & \phi_{\iL_2}  & \i_{\iR} & \v_{\iG}  &\s{c_i}     \\
   \begin{block}{c[cc|cc|cc]l}
f_{\iC_2} &\!{ 0  }   &\!{ 0 }     &  -        &  -         &    -    &    -       &\s1 \\
f_{\iC_1} &\!{ 1  }   &\!{ 1 }     &  -        &  0         &    -    &    -       &\s0 \\ \cline{2-7}        
f_{\iL_1} &  -        &  -         &\!{ 0  }   &\!{ 0  }    &    -    &    -        &\s1 \\
f_{\iL_2} &  -        &  0         &\!{ 1  }   &\!{ 1  }    &    0    &    -      &\s0 \\ \cline{2-7}        
f_{\iR} &  -        &  -         &   -       &   -        &\!{ 0  } &\!{ 0  }     &\s1         \\
f_{\iG} &  -        &  -         &   -       &   0        &\!{ 0  } &\!{ 0  }     &\s1 \\
   \end{block}  
\s{d_j}&  \s1      &   \s1      &  \s1      &  \s1       &  \s1    &  \s1              &     %\\
   \end{blockarray}
   $},\quad
   {\rnc\arraystretch{1.13}
   \left\{\begin{array}{r@{\;=\;}l@{}c}
 0& q_{\iC_2}/C_2-q_{\iC_1}/C_1             &    +V(t) \\
 0&    \dot{q}_{\iC_1}   +\dot{q}_{\iC_2}-\phi_{\iL_2}/L_2   & + I(t)\\
 0&    \phi_{\iL_1}/L_1  -\phi_{\iL_2}/L_2 & +I(t)        \\
 0&   \dot\phi_{\iL_2}  +q_{\iC_1}/C_1+R\i_{\iR} +\dot\phi_{\iL_1}    &-V(t) \\
 0&   \v_{\iG}        +R\i_{\iR}       & -V(t)  \\
 0&   \i_{\iR}- G\v_{\iG}-\phi_{\iL_2}/L_2       &    
   \end{array}\right.
   }
 \end{align}
 }
From \fgref{CpH8by5ex_graph0} one can verify each equation is KCL applied to a cutset, or KVL applied to a cycle, of tree $\iT$.
\end{example}

\subsection{CpH Circuit Model 3}%\label{ss:model3}
For the CpH Circuit Model 3 we make the following assumptions.
\begin{compactenum}[\bf {M3-}1.]
\item The assumptions \Asmp1--\Asmp6 hold.
\item The behavior of the external ports $\iPort$ is described by 
an implicit relation $p(\i_\iPort,\v_\iPort)=0$.
\item The dissipative structure is described by a user-chosen mixed form \rf{mixedForm} (cf. \ssrf{resistexplicit}).
\end{compactenum}
\medskip
 
We restrict to $\iPort=\iVv\dcup\iIi$ and
further split these sets into controlled and independent sources as
\[\iVv=\iVv^\text{c}\dcup\iVv^\text{i},\quad \iIi=\iIi^\text{c}\dcup\iIi^\text{i}.\]
Independent sources are given by explicit source functions $\i_{\iI^\text{i}}=I(t)$ and $\v_{\iv^\text{i}}=V(t)$ and can again be seen as input $u=\bmx{\v_{\iv^\text{i}}\\\i_{\iI^\text{i}}}$, with corresponding output $y=\bmx{\i_{\iv^\text{i}}\\\v_{\iI^\text{i}}}$.
Moreover, we pose the following assumptions for 
controlled voltage source (CVS) and controlled current source (CCS):
\begin{compactenum}[\bf {VS-}1.]
\item No CVS is part of a $\iCc\iVv$-cycle.
\item Controlling voltages of CVS are voltages of capacitances or independent voltage sources.
\item Controlling currents of CVS are currents of inductances or independent current sources.
\end{compactenum}
\begin{compactenum}[\bf {CS-}1.]
\item No CCS is part of an $\iLl\iIi$-cutset.
\item Controlling voltages of CCS are voltages of capacitances, independent voltage sources or resistances.
\item Controlling currents of CCS are currents of inductances, independent current sources or resistances. 
\end{compactenum}

They are adapted from the conditions for controlled sources given in \cite{Tis01}.
Under these assumptions and 
using the constitutive relations $\v_\iCc = \tdbd{\Ham}{q}$, $\i_\iLl = \tdbd{\Ham}{\phi}$, the source functions of independent sources, as well as
\rf{mixedForm} according to {\bf {M3-}3}, the controlled sources are described by
\begin{align*}
 \v_{\iv^\text{c}}&=\gamma_{\iv}\left(q_\iC,q_\ic,\phi_\iL,\phi_\il,t\right),\\
\i_{\iI^\text{c}}&=\gamma_{\iI}\left(q_\iC,q_\ic,\phi_\iL,\phi_\il,\v_\iGg,\i_\iRr,t\right).
\end{align*}
Casting the DAE to input-output form as in \ssrf{ctrlIOform} and restricting 
to the state equation given by \rf{ctrlIOform4a} with a given input 
$u=\bmx{\v_{\iv^\text{i}}\\\i_{\iI^\text{i}}}$ leads to the CpH DAE
\begin{align}\label{eq:CpH_Model3}
0=f(t,x,\dot x,u)
\end{align}
of the form
\begin{align*}
0= &\bmx{ f_{\iv^\text{c}}\\
f_\ic\\ f_\id\\ f_\il} = \bmx{\i_{\iv^\text{c}}\\
\dot q_\ic\\ \i_\id\\ \tdbd{\Ham}{\phi_\il}}
    - \bmx{
    0 &F_{\iD \iv^\text{c}}\+&F_{\iL \iv^\text{c}}\+ &F_{\iI^\text{c} \iv^\text{c}}\+ \\
    \$Cc\+ &\$Dc\+&\$Lc\+ &F_{\iI^\text{c} \ic}\+ \\ 
    0&\$Dd\+&\$Ld\+ &F_{\iI^\text{c} \id}\+ \\ 
    0&0&\$Ll\+ &0 }
\bmx{\dot q_\iC\\ \i_\iD\\ \tdbd{\Ham}{\phi_\iL}\\ \i_{\iI^\text{c}}}
  - \bmx{F_{\iI^\text{i} \iv^\text{c}}\+ \\F_{\iI^\text{i} \ic}\+\\ F_{\iI^\text{i} \id}\+ \\F_{\iI^\text{i} \il}\+ }I(t),\\[.5ex]
0= &\bmx{f_\iC\\ f_\iD\\f_\iL\\  f_{\iI^\text{c}}
} = \bmx{\tdbd{\Ham}{q_\iC}\\ \v_\iD\\ \dot\phi_\iL\\ \v_{\iI^\text{c}}
}
      + \bmx{
       0  &\$Cc &0    &0   \\
      F_{\iD \iv^\text{c}}  &\$Dc &\$Dd &0   \\
      F_{\iL \iv^\text{c}}  &\$Lc &\$Ld &\$Ll\\  
      F_{\iI^\text{c} \iv^\text{c}}  &F_{\iI^\text{c} \ic} &F_{\iI^\text{c} \id} &0 
      }
      \bmx{\v_{\iv^\text{c}}\\ \tdbd{\Ham}{q_\ic}\\ \v_\id\\ \dot\phi_\il}+
      \bmx{F_{\iC \iv^\text{i}}\\F_{\iD \iv^\text{i}}\\F_{\iL \iv^\text{i}}\\ F_{\iI^\text{c} \iv^\text{i}}}V(t),\\
  0=& \bmx{
  p_{\iI^\text{c}}\\
  p_{\iv^\text{c}}
  } = \bmx{
  \i_{\iI^\text{c}}\\
  \v_{\iv^\text{c}}
  }-\bmx{
  \gamma_{\iI}\left(q_\iC,q_\ic,\phi_\iL,\phi_\il,\v_\iGg,\i_\iRr,t\right)\\
  \gamma_{\iv}\left(q_\iC,q_\ic,\phi_\iL,\phi_\il,t\right)
  }, 
\end{align*}
where, similar as before we use that $\i_\id, \i_\iD, \v_\id, \v_\iD$ are explicit functions of $\v_\iGg$ and $\i_\iRr$.
The assumptions {\bf VS-1} and {\bf CS-1} imply that $F_{\iC\iv^\text{i}}=0$ and $F_{\iI^\text{i}\il}=0$.
In \rf{CpH_Model3} the vector $x$ consists of the energy state variables $z=(q,\phi)$ plus the $\^x$ components $\v_\iGg$, $\i_\iRr$, 
$\v_{\iv^\text{c}}$, $\i_{\iI^\text{c}}$, $\i_{\iv^\text{c}}$, and $\v_{\iI^\text{c}}$ resulting in one DAE variable for each capacitor, inductor and  dissipator and 
two DAE variables for each controlled source in the circuit, for a total DAE size of
$ N = n_\iCc+n_\iLl+ n_\iDd + 2(n_{\iv^\text{c}}+n_{\iI^\text{c}})$.

\begin{theorem}
If the assumptions {\bf M3-1}--{\bf M3-3}, {\bf VS-1}--{\bf VS-3}, and {\bf CS-1}--{\bf CS-3}  hold,
and in addition $F_{\iI^\text{c}\id}=0$ (i.e., there are no cycles containing twig resistors and controlled current sources),
then the CpH DAE \rf{CpH_Model3} is SA-amenable.
\end{theorem}

\begin{proof}
For Model 3 the \sigmx takes a schematic form as in \fgref{model3sigma}, where $\Sigma$ has been ordered such that the corresponding Jacobian is in block triangular form.
\begin{figure}
\rnc\![1]{\colorbox{lightgray!110}{$#1$}}
\rnc\g[1]{\colorbox{lightgray!50}{$#1$}}
\begin{footnotesize}
\begin{align*}
\Sigma \ &=\
  \begin{blockarray}{rccccccccclc}
& q_\iC & q_\ic & \phi_\il& \phi_\iL&\v_{\iv^\text{c}}&(\i_\iRr,\v_\iGg)&\i_{\iI^\text{c}}& \i_{\iv^\text{c}} &\v_{\iI^\text{c}}
&\s{c_i}& \\
\begin{block}{r[cc|cc|c|c|c|cc]lc}
 f_\iC &\!{\le0} &\!{\le0} & - & - & - &  - & - & - & -  &\s1 & n_\iC \\
 f_\ic &\!{\le1} &\!{\^\I{+}1} &\g{\le0} & \g{\le0} &-&\g{\le0} &\g{\le0}& -  &-    &\s0 & n_\ic \\ \cline{2-11}        
 f_\il & - & - &\!{\le0} &\!{\le0} & - & - & - &- & - &\s1 & n_\il \\
 f_\iL &\g{\le0}&\g{\le0}&\le1 &\!{\^\I{+}1} &\g{\le0} & \g{\le0} & -  & -  &- &\s0 & n_\iL \\ \cline{2-11}
  p_{\iv^\text{c}}&\le0&\le0&\le0&\le0&\!{\^\I}&-&-&-&-&\s1  &n_{\iv^\text{c}} \\\cline{2-11}
 f_\id & - &- &\le0 &\le0 &- &\!{\le0}  & -  &-   &-     &\s1      & n_\id \\
 f_\iD &\le0 &\le0 & - & - &\le0&\!{\le0}  & - & - & -   &\s1      & n_\iD \\ \cline{2-11}
 p_{\iI^\text{c}}&\le0&\le0&\le0&\le0&-&\le0&\!{\^\I}&-&-&\s1  &n_{\iI^\text{c}}\\\cline{2-11}
 f_{\iv^\text{c}}& - &- & \le0 &\le0 &-&\le0  &\le0&\!{\^\I}   &-      & \s1     & n_{\iv^\text{c}} \\
 f_{\iI^\text{c}} &\le0 &\le0 & - &-&\le0      &\le0     &-& -    &\!{\^\I}   & \s1     & n_{\iI^\text{c}} \\
\end{block}  
\s{d_j} & \s1 &\s1 & \s1 &\s1 &\s1  & \s1 &  \s1   &  \s1&\s1&     \\
&n_\iC& n_\ic & n_\il& n_\iL&n_{\iv^\text{c}} & n_\iDd&n_{\iI^\text{c}}&n_{\iv^\text{c}} &n_{\iI^\text{c}}   &
  \end{blockarray}
\end{align*}\vspace{-6ex}
\end{footnotesize}
\caption{Schematic \sigmx for Model 3.\label{fg:model3sigma}}
\end{figure}
Again, a transversal can be found in the diagonal shaded blocks. 
Blocks with light-grey background cannot contribute nonzero entries to the \sysJ $\_J$.
Thus, $\_J$ takes the form 
{\rnc\![1]{\colorbox{lightgray}{$#1$}}
\begin{align*}
  \_J \ &=\
\begin{blockarray}{rccccccl}
          & q_\iC,q_\ic & \phi_\il,\phi_\iL & \v_{\iv^\text{c}}&\i_\iRr,\v_\iGg  &\i_{\iI^\text{c}}&\i_{\iv^\text{c}},\v_{\iI^\text{c}}              \\
  \begin{block}{r[c|c|c|c|c|c]l}
    f_\iC,f_\ic         & \!{\_J_{\cC}}    & 0 &0&0&0                            & 0            &\s{n_\iCc=n_\iC{+}n_\ic} \\  \cline{2-7}
    f_\il,f_\iL         & 0                & \!{\_J_{\cL}}                 & 0  &0&0&0          &\s{n_\iLl=n_\il{+}n_\iL} \\\cline{2-7}
 p_{\iv^\text{c}}&  \ast&\ast&\!{\I}&0&0&0&\s{n_{\iv^\text{c}}} \\\cline{2-7}   
    f_\id,f_\iD   & \ast  & \ast& \ast &\!{\_J_{\cD}}&0&0 &\s{n_\iDd=n_\iRr{+}n_\iGg} \\\cline{2-7}
p_{\iI^\text{c}}&  \ast&\ast&0&\ast&\!{\I}&0&\s{n_{\iI^\text{c}}}\\\cline{2-7}
f_{\iv^\text{c}},f_{\iI^\text{c}}&  \ast & \ast &\ast &\ast&\ast &\!{\I}       &\s{n_{\iv^\text{c}}{+}n_{\iI^\text{c}}}        \\  
  \end{block}
& \s{n_\iCc}       & \s{n_\iLl}  &\s{n_{\iv^\text{c}}}& \s{n_\iDd} &\s{n_{\iI^\text{c}}}&\s{n_{\iv^\text{c}}}+\s{n_{\iI^\text{c}}}                    
  \end{blockarray}
\end{align*}
}
The $\ast$ denote some (zero or nonzero) blocks that play no role for the further analysis. 
Under the given assumptions the diagonal blocks, and thus $\_J$, are nonsingular as discussed before.
\end{proof}

\begin{theorem}
Assume that {\bf M3-1}--{\bf M3-3}, {\bf VS-1}--{\bf VS-3}, and {\bf CS-1}--{\bf CS-3} hold, and $F_{\iI^\text{c}\id}=0$.
If $n_\iC=n_\il=n_\id=n_\iD=n_{\iv^\text{c}}=n_{\iI^\text{c}}=0$ the system \rf{CpH_Model3} is an ODE,
otherwise a DAE of index 1.
\end{theorem}

\begin{proof}
 Along the same lines as the proof of \thref{index_pHC1}.
\end{proof}

\subsection{Reduction to ODE}\label{ss:RedODE}
An advantage of knowing that the DAE is SA-amenable is that the \Smethod's {\em standard solution scheme} (see \rf{SSS}) tells us which equations need to be differentiated in order to find a numerical solution.
This can be used in several ways, e.g.\ to construct a scheme for reducing the DAE to an explicit ODE.
It is explained in \cite{McKP15} why this might be desirable, e.g.\ for storage reasons, when the number of DOF is small compared to the total number of variables.
We show this process for a Model 2 formulation; with small changes it also works for the other model formulations.

{\rnc\~[1]{\overline{#1}}
\nc\fn{\mathop{\textrm{\sl fun}}}
\begin{algo}{\bf Reduction of Model 2 form to explicit ODE}\label{al:explicitode}
Consider the Model 2 equations \rf{model2block}.
To simplify notation we do not distinguish between column and row vectors.
\begin{compactitem}[]
\item {\em I. Corresponding to stage $k{=}-1$ of the \Smethod:}
\begin{compactenum}[{I.}1.]
\item As $\_J_\cC$ is nonsingular, its $f_\iC$ rows $\tdbd{f_\iC}{q}$ (with $q=(q_\iC,q_\ic)$) have full row rank. Choose columns to form a nonsingular $n_\iC\x n_\iC$ matrix.
Split $q$ as $(\~q,\^q)$ with $\^q$ belonging to the chosen columns.

\item By the IFT we can solve $f_\iC=0$ for $\^q$ in terms of $\~q$ and the right-hand-side item $V(t)$.

Let $\fn$ denote some function that does not need a specific name, then:
\[ \^q = \fn(\~q,V(t)) \]

\item Similarly choose columns of $\tdbd{f_\il}{\phi}$ (with $\phi=(\phi_\il,\phi_\iL)$) to form a nonsingular $n_\il\x n_\il$ matrix.
Split $\phi$ as $(\~\phi,\^\phi)$ with $\^\phi$ belonging to the chosen columns.
\item By the IFT we can solve $f_\il=0$ for $\^\phi$ in terms of $\~\phi$ and the right-hand-side item $I(t)$:
  \[ \^\phi = \fn(\~\phi,I(t)). \]
\item Recall $\i_\id, \i_\iD, \v_\id, \v_\iD$ are explicit functions of $\^x$.
As $\_J_\cD = \tdbd{(f_\id,f_\iD)}{\^x}$ is nonsingular, by the IFT we can solve $(f_\id,f_\iD)=0$ for $\^x$ in terms of $V(t)$, $I(t)$, as well as of $q_\iC,q_\ic,\phi_\il,\phi_\iL$ which reordered are the same as $\~q,\^q,\~\phi,\^\phi$.
That is,
  \[ \^x = \fn(\~q,\^q,\~\phi,\^\phi,V(t),I(t)). \]
\item Now $V(t),I(t)$ comprise the control input $u(t)$ so write the results of steps I.2, I.4, I.5 as
\begin{align}\label{eq:explicitode1}
  \^q &= \fn(\~q,u(t)), &\^\phi &= \fn(\~\phi,u(t)),
  &\^x &= \fn(\~q,\^q,\~\phi,\^\phi,u(t)).
\end{align}
\item Eliminate $\^q,\^\phi$ by substituting the first two equations of \rf{explicitode1} into the third to get
\begin{align}\label{eq:explicitode2}
  \^x = \fn(\~q,\~\phi,u(t)).
\end{align}
\end{compactenum}

\item {\em II. Corresponding to stage $k{=}0$ of the \Smethod, so we differentiate $f_\iC,f_\il$:}
\begin{compactenum}[{II.}1.]
\item Solve $(\.1f_\iC,f_\ic)=0$ as a nonsingular linear system with matrix $\_J_\cC$ for $\qp = (\.1{\~q},\.1{\^q})$ as a function of $(\~q,\^q,\~\phi,\^\phi,\^x,\.1V(t),I(t))$.
Discard $\.1{\^q}$ which will not be used.
\item Solve $(\.1f_\il,f_\iL)=0$ as a nonsingular linear system with matrix $\_J_\cL$ for $\phip = (\.1{\~\phi},\.1{\^\phi})$ as a function of $(\~q,\^q,\~\phi,\^\phi,\^x,V(t),\.1I(t))$.
Discard $\.1{\^\phi}$ which will not be used.
\item Regard $V(t),I(t)$ as comprising $u(t)$ as above.
Eliminate $\^q,\^\phi$ and $\^x$ using \rf{explicitode1,explicitode2} to get $\.1{\~q}$ and $\.1{\~\phi}$ as functions of $(\~q,\~\phi,u(t),\.1u(t))$ yielding a first-order ODE of size $n_\ic{+}n_\iL=\DOF$
\begin{align*}%\label{eq:model2ode}
  \left.\begin{aligned}
  \.1{\~q}    &= \fn(\~q,\~\phi,u(t),\.1u(t)), \\
  \.1{\~\phi} &= \fn(\~q,\~\phi,u(t),\.1u(t)),
  \end{aligned}\right\}.
\end{align*}
\end{compactenum}
\end{compactitem}
\end{algo}

Thus, by three (in general nonlinear) solves corresponding to an $n_\iC\x n_\iC$ part of the $\cC$ block, an $n_\il\x n_\il$ part of the $\cL$ block, and all of the $\cD$ block,
followed by two linear solves corresponding to the whole $\cC$ and $\cL$ blocks, we convert the DAE to an ODE in $(\~q,\~\phi)$, with $u(t), \.1u(t)$ as control input.
The remaining DAE variables, which comprise $\^q,\^\phi$ and $\^x$, can be retrieved from a solution of the ODE by the purely algebraic equations \rf{explicitode1}.

\begin{remark}\rm
\alref{explicitode} is the \Smethod's version of the dummy derivatives process \cite{MatS93}, specialised to the present case. The discarded variables $\.1{\^q},\.1{\^\phi}$ correspond to the ``dummy derivatives'' in the sense of \cite{MatS93}.
See \cite{McKP15} for more details.
In the linear time-invariant case the columns to make nonsingular matrices can be chosen once for all.
But in general a matrix might become singular as solution proceeds, and one needs a new choice of columns (dummy pivoting).
A cheap way to find the needed $n_\iC$ columns of the $f_\iC$ rows, and $n_\il$ columns of the $f_\il$ rows, is the Scholz--Steinbrecher method \cite{SchS16b}: using some HVT of $\Sigma$, choose the columns where it intersects the $f_\iC$ rows and $f_\il$ rows respectively.
In general, this gives nonsingular matrices with high probability, but is not guaranteed.
For CpH DAEs, it always provides a right choice if each inductor and each capacitor is independent. Namely, it chooses $q_\iC$ and $\phi_\il$ for ${\^q}$, ${\^\phi}$ and the resulting nonsingular matrices are the inverses of the corresponding conductance or inductance matrix, respectively.
\smallskip

In code, packing up the repeated nonlinear and linear solves within an ``\li{odefun}'' routine to give to a standard ODE solver makes for opaque programming.
It is more convenient to use a ``reverse communication'' interface so that the solves can be done inline in the calling program.
The NAG Library in 2015 provided a reverse communication Runge--Kutta solver for just such applications, replaced since 2018 by a variable-order variable-step Adams code \li{d02qgf} \cite{NAG18a}.
\end{remark}

\begin{example}
The running example's $\Sigma$-matrix and equations are in \rf{Sigma2}.
The steps of \alref{explicitode} are:
\begin{compactitem}
\item[I.1,I.3] Each circuit element is independent, so the Scholz--Steinbrecher method is guaranteed to succeed.
The main diagonal is an HVT. With this HVT, Scholz--Steinbrecher selects $\^q=q_{\iC_2}$, $\^\phi=\phi_{\iL_1}$.
The ODE variables will then be $(\~q,\~\phi) = (q_{\iC_1},\phi_{\iL_2})$.
\item[I.2] Solve $q_{\iC_2} = \fn(q_{\iC_1},V(t)) 
  $.
\item[I.4] Solve $\phi_{\iL_1} = \fn(\phi_{\iL_2},I(t)) 
  $.
\item[I.5] Solve $(\i_{\iR},\v_{\iG}) = \^x = \fn(\phi_{\iL_2},V(t))$.
In fact $\v_{\iG}$ is not used.
\item[I.6,I.7] $\^x$ is already in the form $\^x = \fn(\~q,\~\phi,u(t))$.
\item[II.1] Solve $(\.1f_\iC,f_\ic)=0$ for $(\qp_{\iC_2},\qp_{\iC_1}) = \fn(\phi_{\iL_2},\.1V(t),I(t))$. Discard $\qp_{\iC_2}$.
\item[II.2,II.3] Solve $(\.1f_\il,f_\iL)=0$ for $(\phip_{\iL_1},\phip_{\iL_2}) = \fn(q_{\iC_1},\i_{\iR},\.1I(t),V(t))$. Discard $\phip_{\iL_1}$.
Eliminate $\i_{\iR}$ by I.5.
\end{compactitem}
This results in an ODE 
$(\qp_{\iC_1},\phip_{\iL_2}) = \fn(q_{\iC_1},\phi_{\iL_2}, u(t),\.1u(t))$
as follows.
\begin{align*}
  \qp_{\iC_1} &= \left(\left(\frac{\phi_{\iL_2}}{L_2}-I(t)\right)/C_2 + \.1V(t)\right)
  \left/\left(\frac1{C_1}+\frac1{C_2}\right)\right., \\
  \phip_{\iL_2} &= \left(-\frac{q_{\iC_1}}{C_1L_1}
  - \frac{R}{L_1}\left(\frac{GV(t) + \phi_{\iL_2}/L_2}{1+GR}\right) + \frac{V(t)}{L_1} + \.1I(t) \right)
  \left/\left(\frac1{L_1}+\frac1{L_2}\right)\right..
\end{align*}
In this symbolic form it is somewhat elaborate. Program code to perform the algorithm steps numerically would be simpler.
\end{example}
}

\subsection{Spectral analysis for linear time-invariant systems}\label{ss:LTI}

When the non-source circuit elements are constant and linear, the system becomes a {\em linear time-invariant (LTI)} system
\begin{align}\label{eq:LTI1}
  0 = Ax - B\xp - u(t),
\end{align}
with constant matrices $A$ and $B$.
One may view \rf{LTI1} as linearization of a nonlinear DAE system around an equilibrium point.
Eigenvalues and eigenvectors of electrical networks can be used to to examine stability of DC operating points; to determine the exponential stability of equilibrium points; or to study bifurcations \cite{Ria06}, \cite[\S6.3.3]{Ria08}.
Several qualitative properties of equilibria can be characterized in terms of the spectrum of the corresponding matrix pencil $A-\lambda B$, see e.g. \cite{Reich95,Ria02,Ria04,RiaT07}.
The condition $\Re\lambda < 0$ guarantees asymptotic stability of the equilibrium in the sense of Lyapunov.
For several reasons, it can be of interest to ensure there exists a unique  equilibrium, which is the case if and only if $A$ is invertible (or all eigenvalues are nonzero). 
Non-trivial, purely imaginary eigenvalues are important in linear circuit applications since they characterize the existence of periodic solutions describing oscillations.
Also linear stability analysis of a power system via \evls can be used to predict its degree of stability \cite{Kun94}.
Sometimes a system starts out as passive but discretisation or model reduction may destroy passivity---whether this has happened can be observed by checking the \evls.

As an example of finding basic LTI system data we show here how to compute eigenfunctions for a Model 2 formulation,
it works in the same way for Models 1 or 3.
One might adapt \alref{explicitode}---which in the LTI case yields a linear constant-coefficient ODE $\xp=Cx+g(t)$---and compute \evls and eigenvectors of the matrix $C$.
However, it is simpler to derive \rf{LTI1} directly and solve the corresponding {\em generalised \evl problem} (GEVP).

In the Model 2 formulation, the DAE variables are $x = (q,\phi,\^x)$, with
$\^x=(\i_\iRr,\v_\iGg)$ according to the splitting of $\iDd=\iRr\dcup\iGg$ into current-controlled (resistor) and voltage-controlled (conductor) edges.
The Hamiltonian Hessians in \Asmp4 and \Asmp5 are now constant and denoted by $\cCl\in\R^{n_\iCc\x n_\iCc}$ and $\cLl\in\R^{n_\iLl\x n_\iLl}$, respectively. Similarly, we denoted the constant Jacobian  of $\rho(\i_\iRr,\v_\iGg)$ in \rf{jacR} by 
$R\in\R^{n_\iDd\x n_\iDd}$.
The relation \rf{mixed2} becomes
\begin{align*}
  \bmx{\i_\iDd\\ \v_\iDd} &= P_\iGg \bmx{\^x\\ R\^x} = P_\iGg\bmx{\I\\R}\^x 
  = \bmx{\cDl_\i \\\cDl_\v}\^x, \quad\text{say},
\end{align*}
with $\cDl_\i, \cDl_\v \in\R^{n_\iDd\x n_\iDd}$, where $P_\iGg$ is the flip-matrix of the split $\iDd=\iRr\dcup\iGg$.
Let the edge set $\iXx = \iCc\dcup\iLl\dcup\iDd$ be split into $\ix\dcup\iX$ as in \ssrf{ctrlIOform}.
The constitutive relations, being linear, give the current and voltage vectors on $\iXx$ as
\begin{align}
  \i_\iXx &=\bmx{\i_\iCc\\ \i_\iLl\\ \i_\iDd} =\bmx{\qp\\ \cLl\phi\\ \cDl_\i\^x }
  = \bmx{0&&\\ &\cLl&\\ &&\cDl_\i}x - \bmx{-\I&&\\ &0&\\ &&0}\xp &= A_\i x - B_\i\xp,
  \text{ say}, \label{eq:LTI5i} \\
  \v_\iXx &=\bmx{\v_\iCc\\ \v_\iLl\\ \v_\iDd} =\bmx{ \cCl q\\ \phip\\ \cDl_\v\^x}
  = \bmx{\cCl&&\\ &0&\\ &&\cDl_\v}x - \bmx{0&&\\ &-\I&\\ &&0}\xp &= A_\v x - B_\v\xp,
  \text{ say}. \label{eq:LTI5v}
\end{align}
Here $A_\i, B_\i, A_\v, B_\v$ are $N\x N$ matrices, where $N=n_\iCc+n_\iLl+n_\iDd$ as in \rf{CpH_Model2b} is the DAE size.

Define the coordinate projection matrices $\Pi_\ix\in\R^{n_\ix\x N}$, $\Pi_\iX\in\R^{n_\iX\x N}$ that extract the $\ix$ and $\iX$ coordinates of a column $N$-vector.
They consist of the $\ix$-rows and $\iX$-rows respectively, of the $N\x N$ identity matrix. Then $\i_\ix = \Pi_\ix\i_\iXx$, $\i_\iX = \Pi_\iX\i_\iXx$  and so on.
Thus, \rf{ctrlIOform2} may be written as
{\rnc\arraystretch{1.1}
\begin{align}
  0 &=\bmx{f_\ix\\f_\iX} 
    =\bmx{\Pi_\ix\, \i_\iXx - (\$\iX\ix\+ \Pi_\iX\, \i_\iXx + \$\iI\ix\+ I(t)) \\
          \Pi_\iX\, \v_\iXx + (\$\iX\ix   \Pi_\ix\, \v_\iXx + \$\iX\iv   V(t)) }
    =\bmx{F_\i\, \i_\iXx \\ F_\v\, \v_\iXx}
    +\bmx{-\$\iI\ix\+ I(t) \\ \quad\$\iX\iv V(t)}, \notag
\inter{where $F_\i = \Pi_\ix - \$\iX\ix\+ \Pi_\iX$, $F_\v = \Pi_\iX + \$\iX\ix \Pi_\ix$. With \rf{LTI5i} and \rf{LTI5v} this gives}
 0 &=\bmx{F_\i A_\i \\ F_\v A_\v}x - \bmx{F_\i B_\i \\ F_\v B_\v}\xp
      +\bmx{-\$\iI\ix\+ I(t) \\ \quad\$\iX\iv V(t)} \notag\\
  &= Ax - B\xp - u(t), \label{eq:LTI6}
\end{align}
}
which has the form \rf{LTI1}, with
\begin{align}\label{eq:LTI7}
   A &= \bmx{F_\i\,A_\i \\F_\v A_\v}, 
  &B &= \bmx{F_\i\,B_\i \\F_\v B_\v},
  &u(t) &= \bmx{-\$\iI\ix\+ I(t) \\ \quad\$\iX\iv V(t)}.
\end{align}
Matrix rows, i.e.\ equations, have undergone reordering during this process, but this does not change the DAE's solutions so we may regard \rf{LTI6,LTI7} as a way, for the LTI case, to write \rf{model2block}.
\bigskip

Eigenfunctions are solutions of the corresponding homogeneous system $Ax - B\xp = 0$ that have the form $x(t) = ve^{\lam t}$ for some constant vector $v\ne0$ and complex scalar $\lam$.
Then $\lam,v$ solve the GEVP $Av = \lam Bv$, so that $\lam$ is a zero of the characteristic polynomial $p(\lam) = \det(A-\lam B)$.
As a result of the tree being normal according to \Asmp3, derivatives occur only in $f_\ic$ and $f_\iL$ of \rf{model2block}.
Only the corresponding rows of $B$ can have nonzeros, a total of $n_\ic+n_\iL=\DOF$ rows.
Hence by the formula for determinant as a sum of products over matrix transversals, the degree $\deg p$ of $p(\lam)$ cannot exceed \DOF.
Non-singularity of the \sysJ (a result of \Asmp4--\Asmp6) implies that $(A,B)$ is a regular matrix pencil that has exactly $n_\ic+n_\iL$ eigenvalues or,
equivalently, that $\det(A-\lam B)$ is a polynomial in $\lambda$ with degree $n_\ic+n_\iL$\footnote{This results from the structure of $B$ and the Kronecker canonical form, see e.g. \cite{Ria02}.}.
These are the finite \evls of the GEVP; its remaining $(N{-}\DOF)$ \evls are $\infty$ and are discarded (see \cite[p.90ff]{Dua10}).
\smallskip

For simplicity assume the problem has \DOF linearly independent \evcs $v_i$ (true e.g.\ if the finite \evls are distinct).
Then  each solution of $Ax - B\xp = 0$ has a unique expansion 
\begin{align*}%\label{eq:LTI9}
  x(t) = \sum_{i=1}^{\DOF} c_i e^{\lam_i t} v_i.
\end{align*}

\begin{example}
The running example equations in \rf{Sigma2}  can be written in the LTI form \rf{LTI1} as
{\footnotesize\begin{align*}
0 &=
\bmx{
 1/C_2 &-1/C_1 &   0  &   0   & 0  & 0  \\
   0   &   0   &   0  &-1/L_2 & 0  & 0  \\
   0   &   0   &1/L_1 &-1/L_2 & 0  & 0  \\ 
   0   &1/C_1  &   0  &   0   &R & 0  \\ 
   0   &   0   &   0  &-1/L_2 & 1  &-G\\ 
   0   &   0   &   0  &   0   &R & 1
 }
 \bmx{q_{C_2}\\ q_{C_1}\\ \phi_{L_1}\\ \phi_{L_2}\\ \i_{R}\\ \v_{G}}
 -
 \bmx{
  0 & 0 & 0 & 0 & 0 & 0 \\
 -1 &-1 & 0 & 0 & 0 & 0 \\
  0 & 0 & 0 & 0 & 0 & 0 \\
  0 & 0 & -1&-1 & 0 & 0 \\
  0 & 0 & 0 & 0 & 0 & 0 \\
  0 & 0 & 0 & 0 & 0 & 0 \\
 }
 \bmx{\qp_{C_2}\\ \qp_{C_1}\\ \phip_{L_1}\\ \phip_{L_2}\\ \.1\i_{R}\\ \.1\v_{G}}
 - \bmx{-V(t)\\ -I(t)\\ -I(t)\\ V(t)\\ V(t)\\ 0} \\
 &= Ax - B\xp - u(t).
\end{align*}}\noindent

The characteristic polynomial is $p = \det(A-\lam B) = c(\alpha\lam^2 + \beta\,\lam + \gamma)$ where
\begin{align*}
   c &=-R/(C_1C_2L_1L_2),
  &\alpha &=(G+1/R)(C_1+C_2)(L_1+L_2),
  &\beta &=C_1+C_2,
  &\gamma &=G+1/R.
\end{align*}
The finite \evls are
\begin{align*}
  \lam_1,\lam_2 &=-\delta \pm \sqrt{\delta^2 -\eta}
  &\text{where }  \delta &= \frac1{2(G+1/R)(L_1+L_2)} = \frac\beta{2\alpha}, 
  &\eta &= \frac1{(C_1+C_2)(L_1+L_2) }= \frac\gamma\alpha.
\end{align*}
\end{example}

\section{Relation to other circuit formulations}\label{sc:otherforms}

Commonly used methods for modelling electrical circuits are the {\em Modified Nodal Analysis (MNA)} \cite{VlaS94}, the {\em Modified Loop Analysis (MLA)} \cite{Wei84}, or the analysis via {\em Branch-oriented Models (BOM)} \cite{Ria08} (or tree-based formulations) also known as {\em hybrid analysis} in the circuit literature \cite{Ama62,Bra62, Har63,Kro39}.
It has been shown in \cite{IwaT10, IwaTT10, IwaTT12, TakI10} that the  index of the DAE resulting from the hybrid analysis is at most one for strictly locally passive components, in contrast to MNA and other nodal techniques, for which certain configurations (as $\iCc\iVv$-cycles or $\iLl\iIi$-cutsets) yield index two systems . 

While MNA (and also MLA) naturally gives a DAE of \pH structure, see e.g.\ \cite{GerHRvdS21,GueBJR20} for related results, simple examples show that the resulting DAE need not be SA-amenable, see \cite{Sch18}.
In \cite{Sch18} a {\em branch-oriented model} in \pH form based on a normal tree is  presented that is closely related to the results presented here, but uses different variables.
A different definition for {\em \pH descriptor systems} is given in \cite{BeaMXZ18,MehM19}. 
With a slight increase in the dimension of the system
-- by adding the variables $\i_\iC$, $\v_\il$ and the corresponding relations $\v_\il =\dot\phi_\il$, $\i_\iC=\dot q_\iC$ -- our formulation can be transformed into this form.
 Due to \Asmp3 $n_\il$ and $n_\iC$ are usually small.
The MNA formulation is based on Kirchhoff's laws written as
\begin{align}\label{eq:KL_MNA}
\v = \Ap\+\eta\quad\text{ and}\quad \Ap\i = 0,
\end{align}
where  $\Ap\in\R^{(\nT)\x\nb}$ is a {\em reduced incidence matrix}, i.e., $A$ with one row removed, and
$\eta$ denotes the $(\nv-1)$-dimensional vector of node potentials.
By \cite[Proposition 5.4]{biggs1993algebraic} $\Ap_\iT$ is nonsingular, and
\begin{align} \label{eq:Fdef}
  F &= -(\Ap_\iT\` \Ap_\iN)\+, \quad\text{independent of which row was removed}.
\end{align}
A transformation of \rf{KL_MNA} with $\widetilde A_\iT^{-1}$ gives 
$\v = (\widetilde A_\iT^{-1}\Ap)\+ \mu$ and $\widetilde A_\iT^{-1}\Ap\i = 0$, 
 where $\mu=\widetilde A_\iT\+\eta $.
Ordering $\Ap$ into tree and cotree parts as $[\Ap_\iT, \Ap_\iN]$ , similarly the vectors $\v$ and $\i$, then 
with \rf{Fdef} we get \rf{KL} and $\mu$ becomes the tree voltages $\v_\iT$.
Similar relations also holds for the MLA equations. The BOM formulation is similar to $\rf{KL}$ without using the properties of a normal tree.

The MNA equations  contain the variable $\i_\iV$, while in the CpH approach $\i_\iV$ can be either separated out (as e.g.\ in Model 2) or considered together with the complementary variable $\v_\iI$ as output of the \pH system. 
Also note that the CpH method does not require a ground node as do nodal-based methods.

\begin{example}
The MNA formulation of the running example, with vertex 2 chosen as ground, is
\begin{small}
\begin{align*}
C_1(\dot \eta_1 - \dot \eta_4)+ G  \eta_1+\i_{\iL_1}+\i_\iV&=0\\
C_2(\dot \eta_3 -\dot \eta_4)+ (1/R)\eta_3-\i_\iV&=0\\
C_1(\dot \eta_4-\dot \eta_1)+C_2(\dot \eta_4-\dot \eta_3)+\i_{\iL_2}&=0\\
-\i_{\iL_1}-\i_{\iL_2}+I(t)&=0\\
L_1\d\i_{\iL_1}/\d t - \eta_1+\eta_5&=0\\
L_2\d\i_{\iL_2}/\d t - \eta_4+\eta_5&=0\\
\eta_1-\eta_3-V(t)&=0.
\end{align*}
\end{small}
It is SA-unamenable. However, choosing another vertex as ground leads to an SA-amenable system. 
\end{example}

\section{Numerical Implementation and Examples}\label{sc:Numerics}

Since the CpH DAEs in this paper have index $\le1$, they can be solved by a DASSL-type code. Here we use \matlab's index-1 solver \li{ode15i}. 
\ssref{matlabimplem} outlines the main ideas of our \matlab implementation of the CpH formalism. The computation of a normal tree and it's \ccmatrix is discussed in \ssref{computeTandF}.
\ssref{getCICs} discusses computing consistent initial conditions (CICs)
 and
\ssref{Examples} presents numerical results.

{\lstset{keywords={}}
\subsection{\matlab implementation}\label{ss:matlabimplem}
We use three \matlab classes:
\begin{compactitem}[$\bullet$]
\item \li{pHedge} implements independent, possibly nonlinear, circuit elements.
  An object is a vector of edge elements each having three fields:
  \li{type}, indicating the kind of edge, one of \iV, \iC, \iR, \iG, \iL, \iI;
  \li{name}, a text name of the edge;
  \li{value}, a number or a function (handle).
  If \li{value} is a number $\xi$, it defines a constant current or voltage for a \iV\ or \iI\ edge, and linear behaviour for a \iC, \iR, \iG\ or \iL\ edge, e.g.\ $\v=\xi\i$ for an \iR\ edge and $\i=\xi\v$ for a \iG\ edge.
  If \li{value} is a function $f$, it defines $V(t)$ or $I(t)$ for a \iV\ or \iI\ edge, and the relation $\v=f(q)$, $\v=f(\i)$,  $\i=f(\v)$, $\i=f(\phi)$ for a  \iC, \iR, \iG\ or \iL\ edge respectively.
\item A \li{pHcircuit} object is mainly a \li{pHedge} object together with an incidence matrix \li{A} stored in sparse form.
  An important method is \li{pHjoin} which ``solders together'' specified nodes of one or more \li{pHcircuit} objects to create a new one combining all the edges of the input objects.
  One can specify only a subset of a circuit's nodes as ``visible'', i.e., usable by \li{pHjoin}.
\item \li{CpH} represents a circuit by a \li{pHedge} object and its \ccmatrix $F$  plus tree $\iT$ and cotree $\iN$.
  The constructor makes a \li{CpH} object \li{P} from a \li{pHcircuit} object \li{P0} and optionally a proposed tree \iT.
   If \iT\ is not given, this constructor computes one. If a normal tree $\iT$ and cotree $\iN$ have been found, the \ccmatrix $F$ is computed as in \ssref{computeTandF}.
  \li{P} stores sufficient data to reverse the construction, i.e., recover \li{P0} exactly from \li{P}.
\end{compactitem}

There are two methods to construct from a CpH object \li{P} a function $f(t,x,\xp)$ that implements CpH Model 2 and is accepted by \li{ode15i} \cite{Sha97MO}: one that ``interprets'' \li{P} as this $f$, and one that for more efficiency generates a \matlab \li{daefcn} function for this $f$.
We actually solve the DAE by a slight adaptation of \li{ode15i} which at each $t_i$ returns numerical approximations of $\xp_i$ as well as of $x_i$.
This supports computing the control-output $y$, which is a function of $t$, $x$ and $\xp$.

\subsection{Computation of a normal tree and its \ccmatrix}\label{ss:computeTandF}
  
From \thref{normalT}, a normal tree and its \ccmatrix $F$ are a byproduct of computing a suitable RREF of the incidence matrix $A$.
\matlab provides the function \li{[R,T] = rref(A)} to do this using Gauss-Jordan elimination with partial pivoting. 
The resulting \li{R} is the echelon form; \li{T}, the integer vector of pivot column indices, is the tree; and $F$ can be formed from the non-pivot columns, see \coref{F}.
A faster version of \matlab's RREF, \li{frref}, that uses sparse QR algorithm to compute RREF for sparse matrices has been provided by \cite{Ataei21}.
Since floating point is used in the computations, the obtained $F$ can differ from integer values at roundoff level, so we round it to the nearest integer. 

Both \li{rref} and \li{frref} are very inefficient for finding a tree.
A much faster approach is to use Kruskal's algorithm (KA) \cite{Kru56} for finding a minimum spanning tree; see below and Figure~\ref{fig:rref_kruskal}. In this algorithm,  we do not sort the edges but provide them in the order $\iVv, \iCc, \iDd, \iLl, \iIi$.
The ordering in each of these sets does not matter.
Starting with $\iT = \emptyset$, KA adds an edge $e$ to $\iT$ if $\iT\cup\{e\}$ is acyclic.
This results in a normal tree as long as \Asmp1 and \Asmp2 hold.
 If a normal tree $\iT$ and cotree $\iN$ are found, the \ccmatrix $F$ is then computed from \rf{Fdef} by
\li{F = -inv(A(1:n-1,T))*A(1:n-1,coT)},
using \matlab notation.
Surprisingly, using the \li{inv} function is more efficient than finding $F$ through \matlab's backslash operator (see Figure~\ref{fig:kaf}).
\begin{figure}[ht]
\centering
\subfigure[CPU time versus $b$ for computing a normal tree  using \li{rref}, \li{frref}, and KA \label{fig:rref_frref_ka}]{
\includegraphics[width=.47\textwidth]{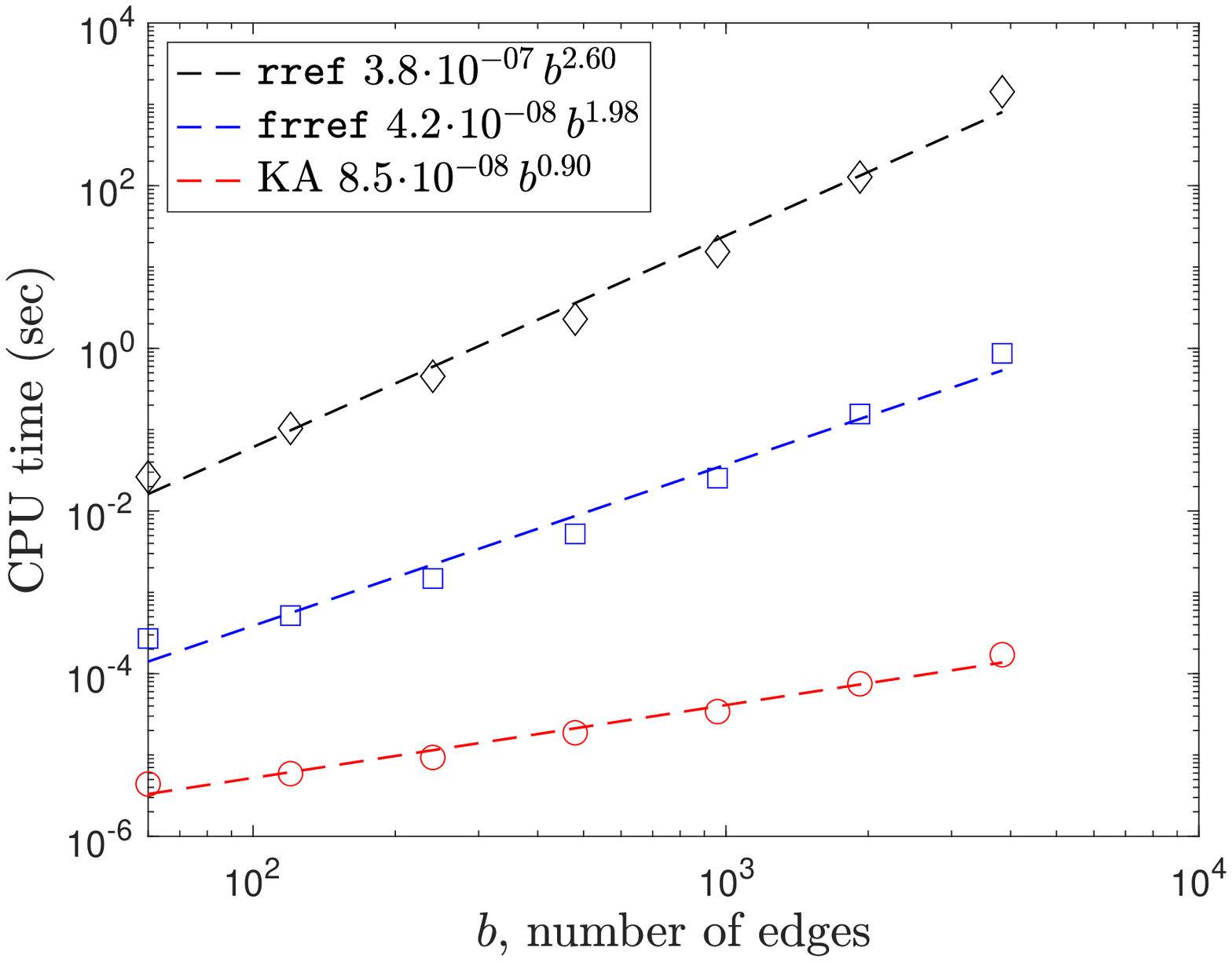}}\hfill
\subfigure[CPU time versus $b$ for {KA} and computing $F$ using \texttt{\textbackslash} and
\li{inv}  
\label{fig:kaf}]{\includegraphics[width=.47\textwidth]{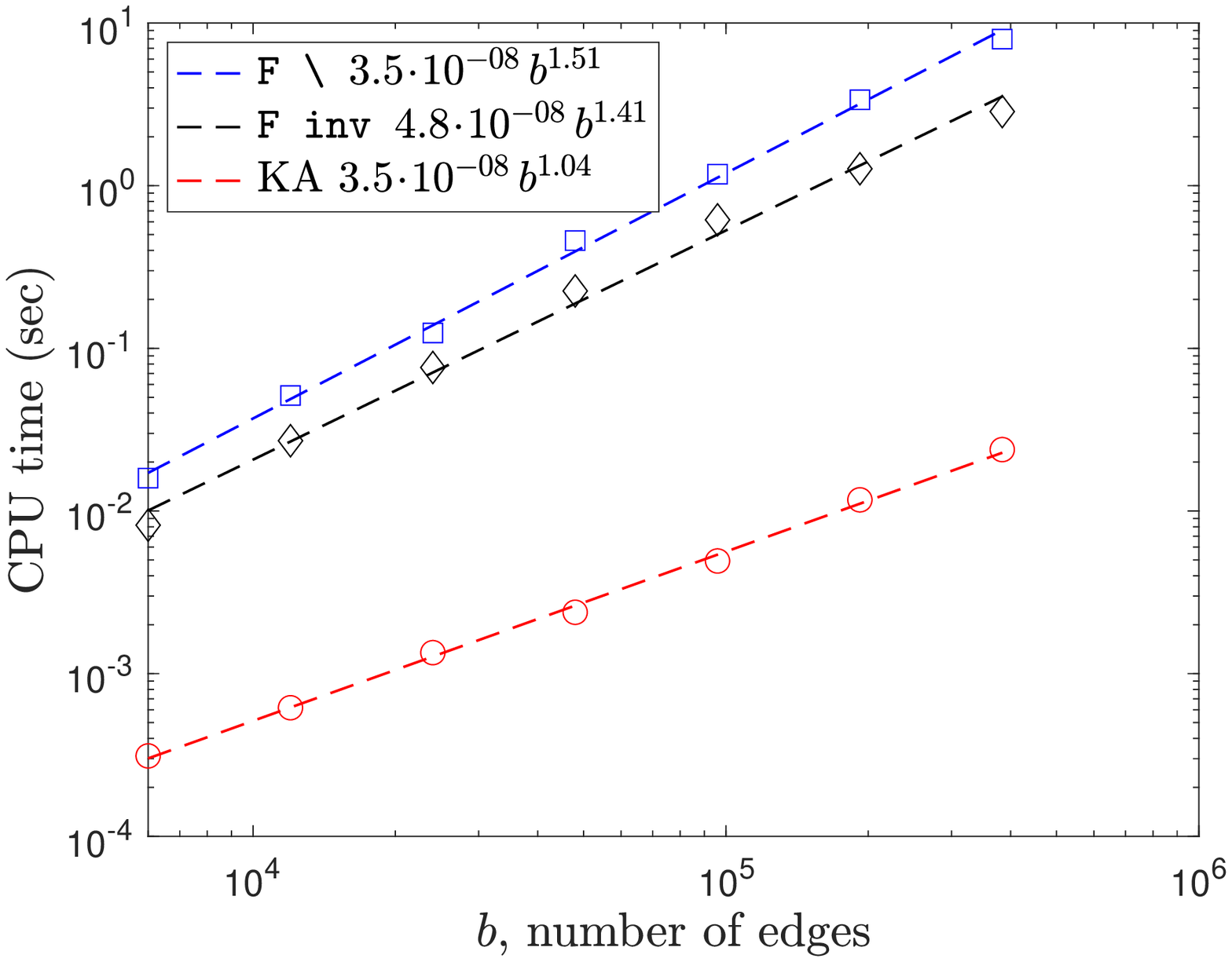}}
\caption{Performance of computing a normal tree and $F$.
CPU times are indicated by the markers; lines are least-squares fits.
\label{fig:rref_kruskal}}
\end{figure}
The timing results  are produced with \matlab version R2021a on Mac OSX with 2.2 GHz Quad-Core Intel Core i7
 and 16 GB 1600 MHz DDR3.
 
\sisetup{group-minimum-digits = 3,group-separator = {,}}

 We use sparse incidence matrices obtained from random connected graphs with $n$ vertices and $b = 3n$ edges. 
In Figure~\ref{fig:rref_frref_ka}, we plot the CPU time 
of  \li{rref}, \li{frref}, and KA for 
$n = 20, 40, \ldots, 1280$.   For each set of data points
(number of edges, CPU time), we also plot  $c_1b^c$, where the constants $c_1, c$ are determined from a least-squares fit of these data. 
In Figure~\ref{fig:kaf}, we plot the CPU time for KA
and computing $F$ using \li{inv} and \texttt{\textbackslash}, 
where $n = 2000, 4000, \ldots, 128000$, and  also plot the corresponding   $c_1b^c$ fits. 
The time for the overall computation of $F$ is dominated by  \li{inv}.
This matrix can be directly derived from the incidence matrix of the tree by an algorithm that is purely based on logical operations \cite[Appendix III]{Bra62}. We are planning to compare the efficiency of this approach to our current method.

}

\subsection{Consistent initial conditions}\label{ss:getCICs}

\nc\RN[2]{{\color{green}#1} {\color{red} #2}}

\nc\xtilde{\widetilde x}
Computing consistent initial conditions (CICs),
as in Step 6 of the \Sigma-method in \scrf{sigmethod}, is challenging for general DAEs. Usually,  this  involves index reduction and  finding hidden constraints \cite{Est02,EstL18}.
A solver such as \li{ode15i} requires consistent values for derivatives as well.~\matlab provides the  function \li{decic} that seeks CICs for an index-1 DAE  prior to solving by \li{ode15i}.~The \li{decic} function works well for a {\em semi-explicit} DAE but a CpH DAE is usually not so, and \li{decic}   may give  wrong $\xp$ values.
An example is $x_1+x_2=0,\;x_1-\xp_2=0$, which comes from a circuit comprising just a resistor $R_1=1$ and capacitor $C_2=1$.
In our CpH approach, having correct approximations for derivatives is crucial, as we insert them in  \rf{IOformb} to obtain the output $y$.

The DAE being SA-amenable index-1, as here, eases the CICs task.
From results of \scrf{CpHmodels}, for both models 1 and 2, one can fix $n_\ic$ charges $q_j$ and $n_\iL$ fluxes $\phi_j$.
In case of independent elements, they can be exactly the tree capacitor charges and cotree inductor fluxes.  In general, a correct choice comes from rank properties of the $f_\iC$ and $f_\il$ rows of $\_J$, see \rf{Jblocks}. 
This ensures a unique solution---locally in the nonlinear case, globally in the linear.

We use the \Smethod's SSS, cf.\ \rf{SSS}. At stage $k=-1$, for a given initial guess $\xtilde_0$ at $t_0$  solve $0=f_\iC,f_\il,f_\id,f_\iD,r$ in \rf{model1Sigma} for the whole vector $x$: a full-rank, generally underdetermined system.
Let $x_0$ be the result.
Now let $x_1(h)$ denote the result of one step of implicit Euler from $(t_0,x_0)$ to $t_1=t_0+h$.
For a DAE with analytic coefficients, this has a series expansion $x_1(h) = x_0 + \xp(t_0)h + c_1h^2 + c_2h^3 + \cdots$ where $x(t)$ is the exact solution through $(t_0,x_0)$.
Hence the central difference $\frac{x_1(h)-x_1(-h)}{2h}$ is an $O(h^2)$ accurate estimate of $\xp(t_0)$.
Error estimates, and higher-order approximations, can be found by Richardson extrapolation.
We aim to detail this in a future paper.

\subsection{Numerical Examples}\label{ss:Examples}

First we consider the circuit from \exref{ex1} with parameter values
 $C_1=C_2=5\cdot 10^{-6}$,  $G=R=1$, $L_1=L_2=0.1$ and source functions $V(t)=10t\sin(200\pi t)$, $I(t)=10\sin(10t)$. 
 We integrated 
 the Model 2 formulation over the time interval $[0,0.2]$ using the  initial guess $\widetilde x_0=[1,\dots,1]\+$.
The numerical solution for $x=[q_{\iC_2},q_{\iC_1},\phi_{\iL_1},\phi_{\iL_2},\i_\iR,\v_\iG]$ is depicted in \fgref{outputRUNEX}(a). The output $y=[\i_\iV,\v_\iI]$ is given in \fgref{outputRUNEX}(b).
The analytical solution can be found for this problem, and the numerical one agrees with it (see red dashed line in \fgref{outputRUNEX}(b)).

\begin{figure}[htbp]
\centering
\subfigure[Numerical solution]
{
\nc\mpw{.32\TW}
\begin{minipage}{.7\TW}
\includegraphics[width=\mpw]{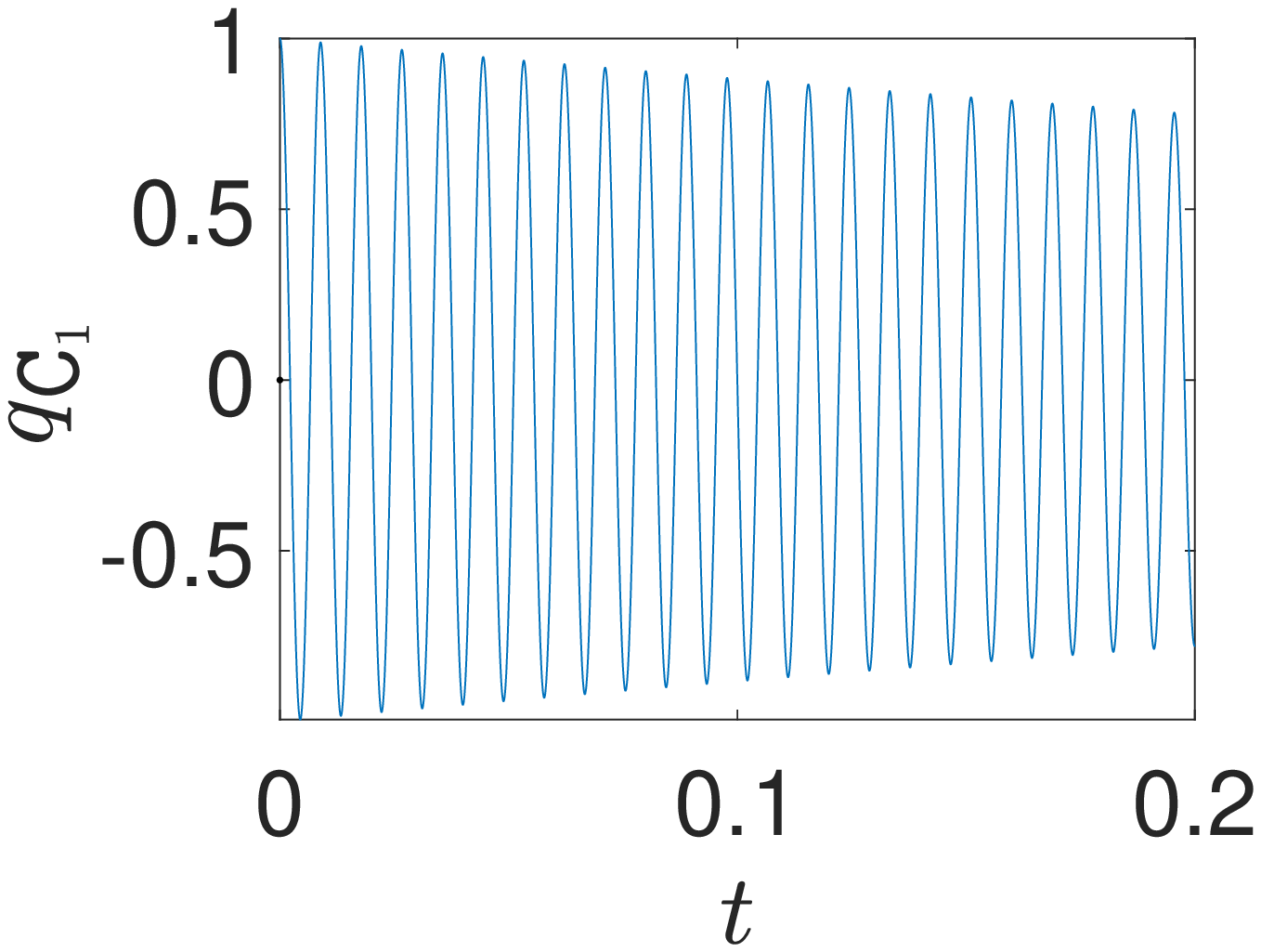}
\includegraphics[width=\mpw]{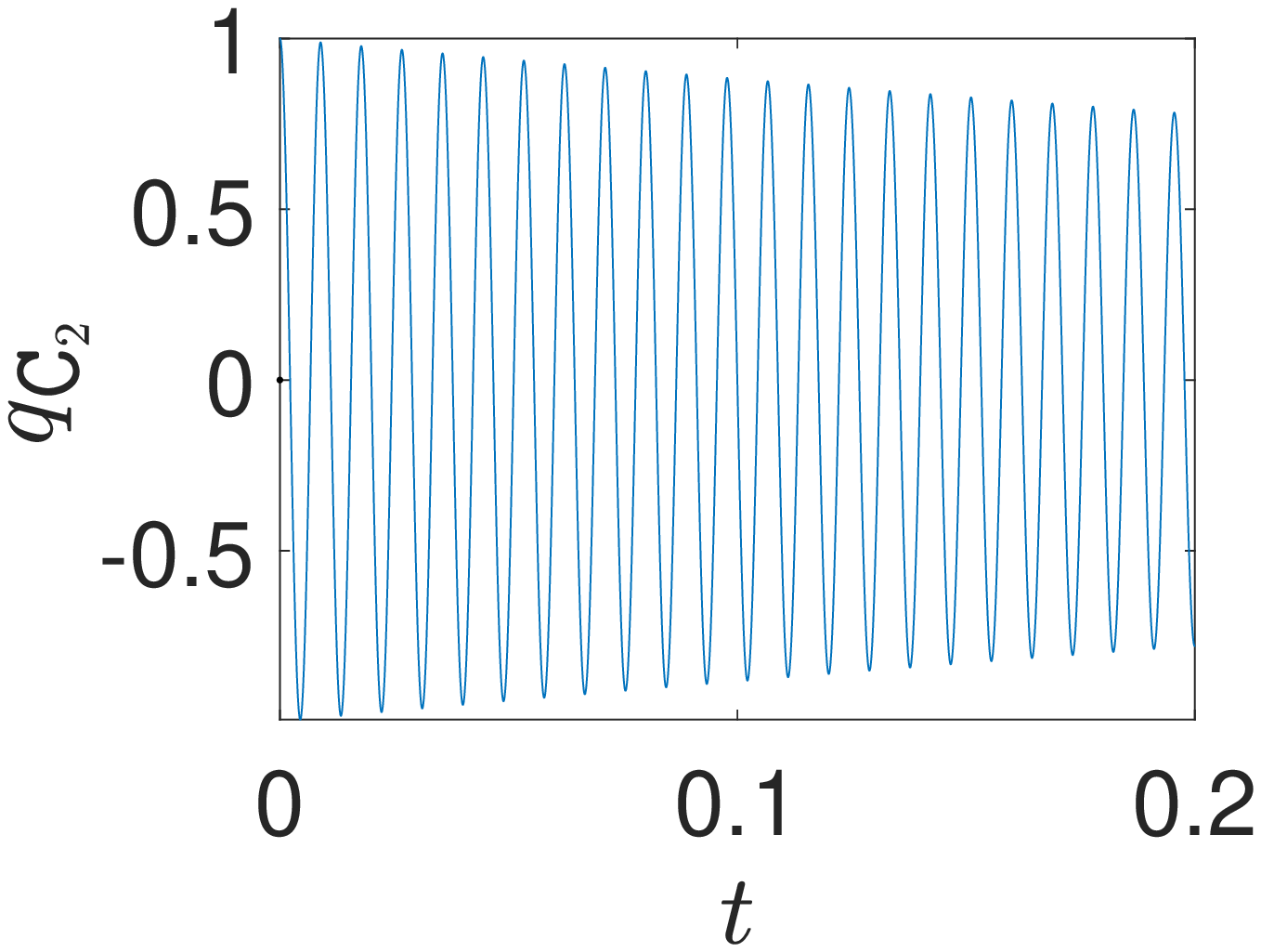}
\includegraphics[width=\mpw]{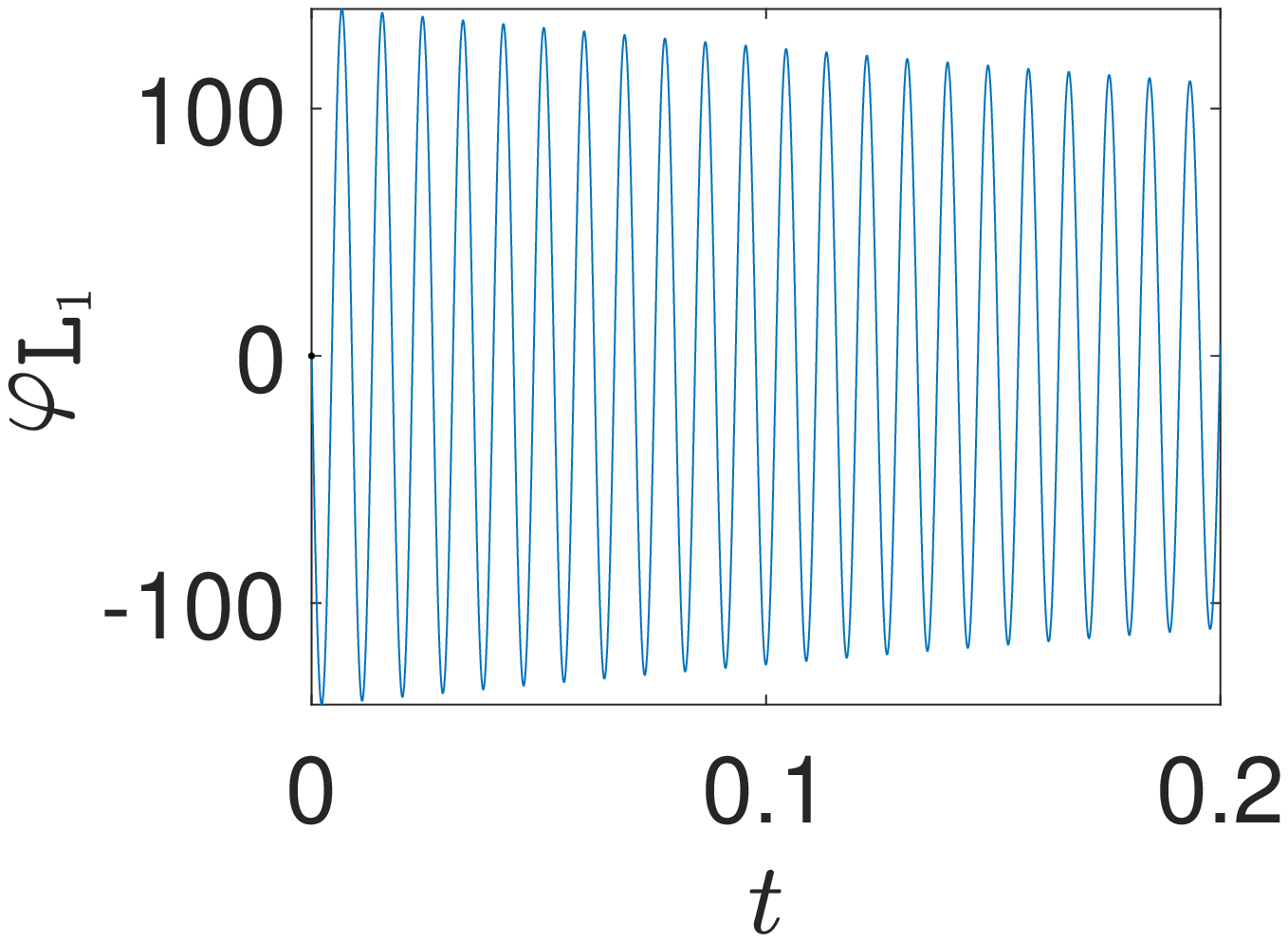}

\includegraphics[width=\mpw]{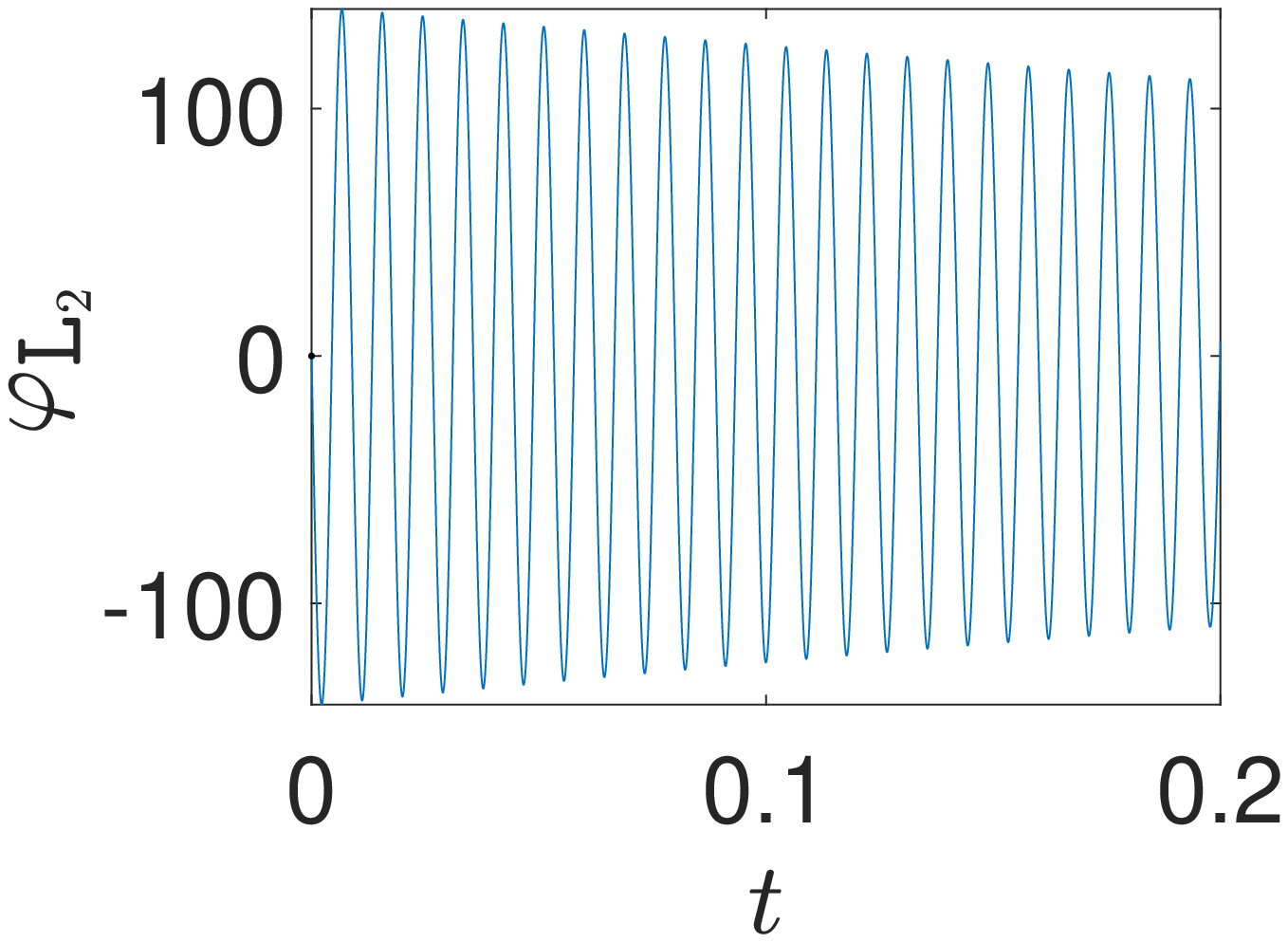}
\includegraphics[width=\mpw]{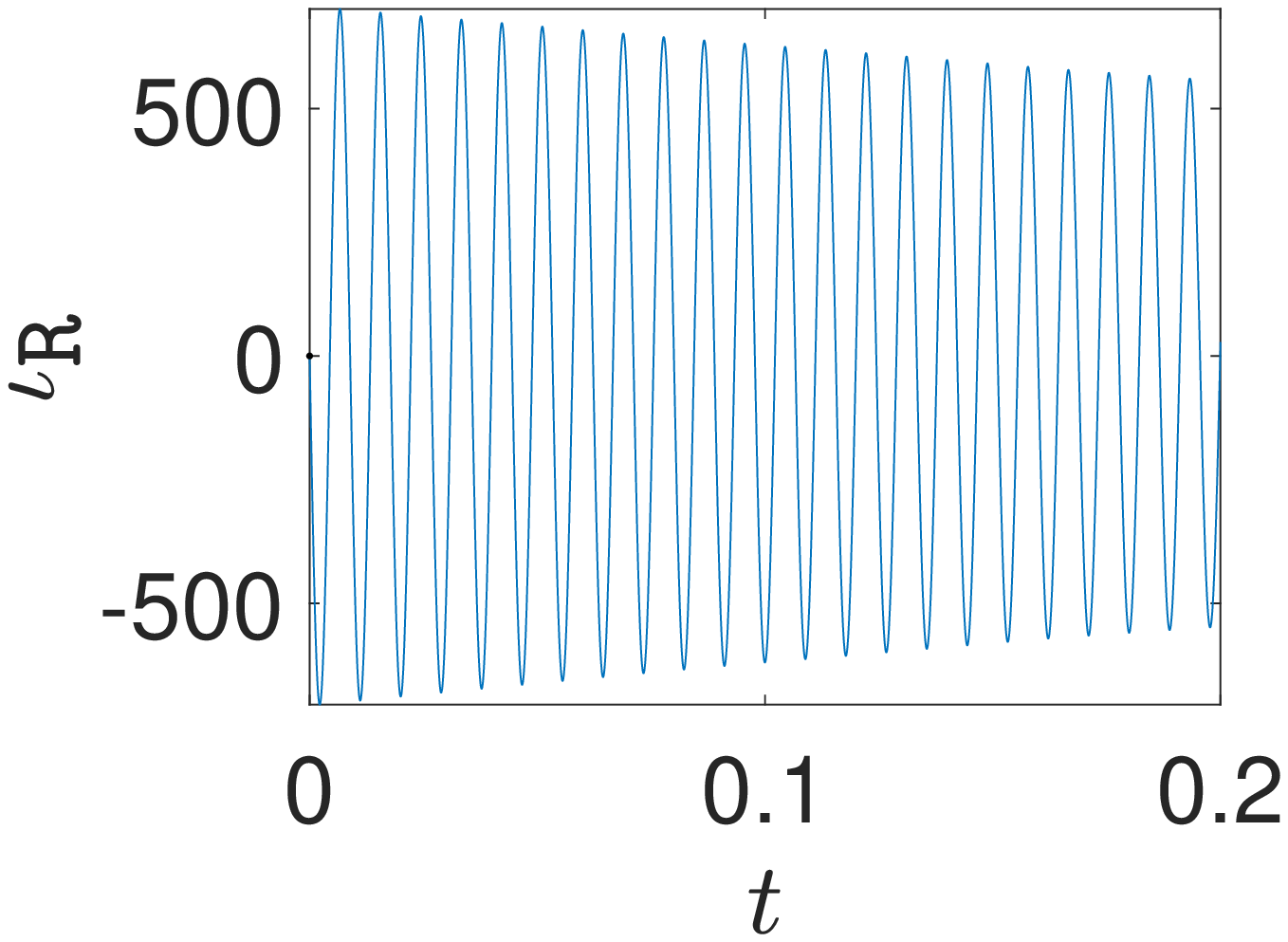}
\includegraphics[width=\mpw]{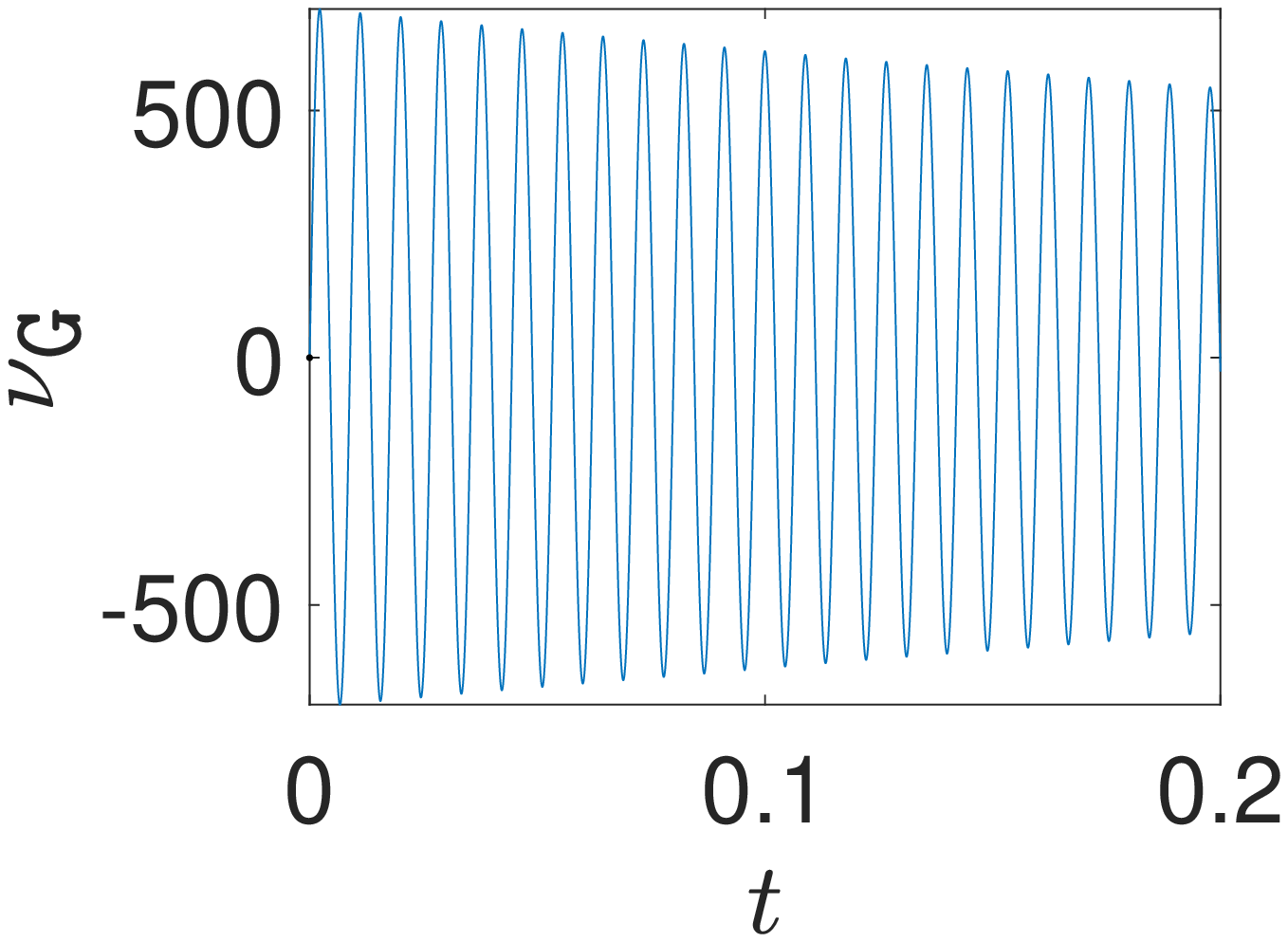}

\end{minipage}
}
\subfigure[Control output]{
\nc\mpw{\textwidth}
\begin{minipage}{.21\TW}
\includegraphics[width=\mpw]{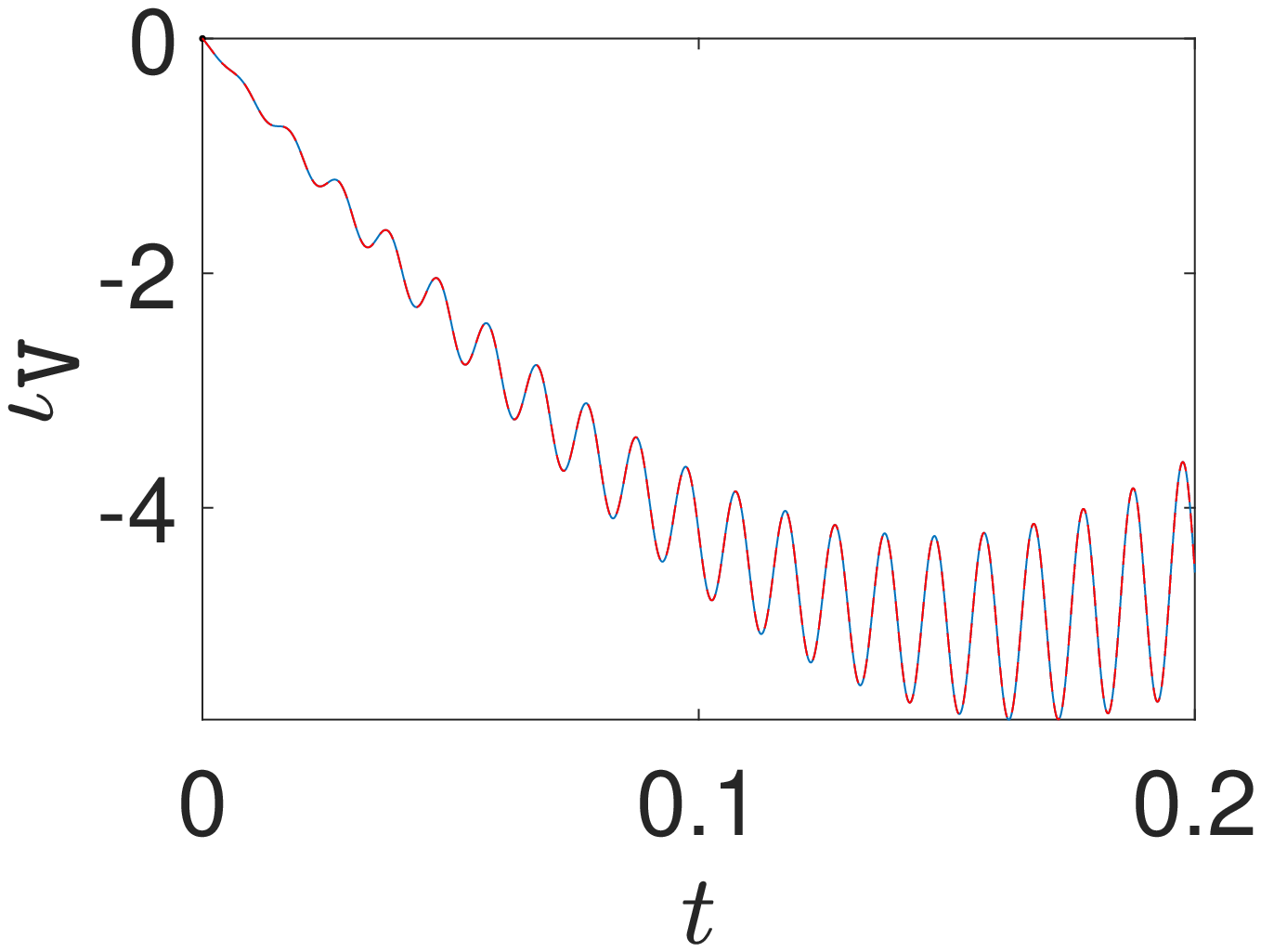}

\includegraphics[width=\mpw]{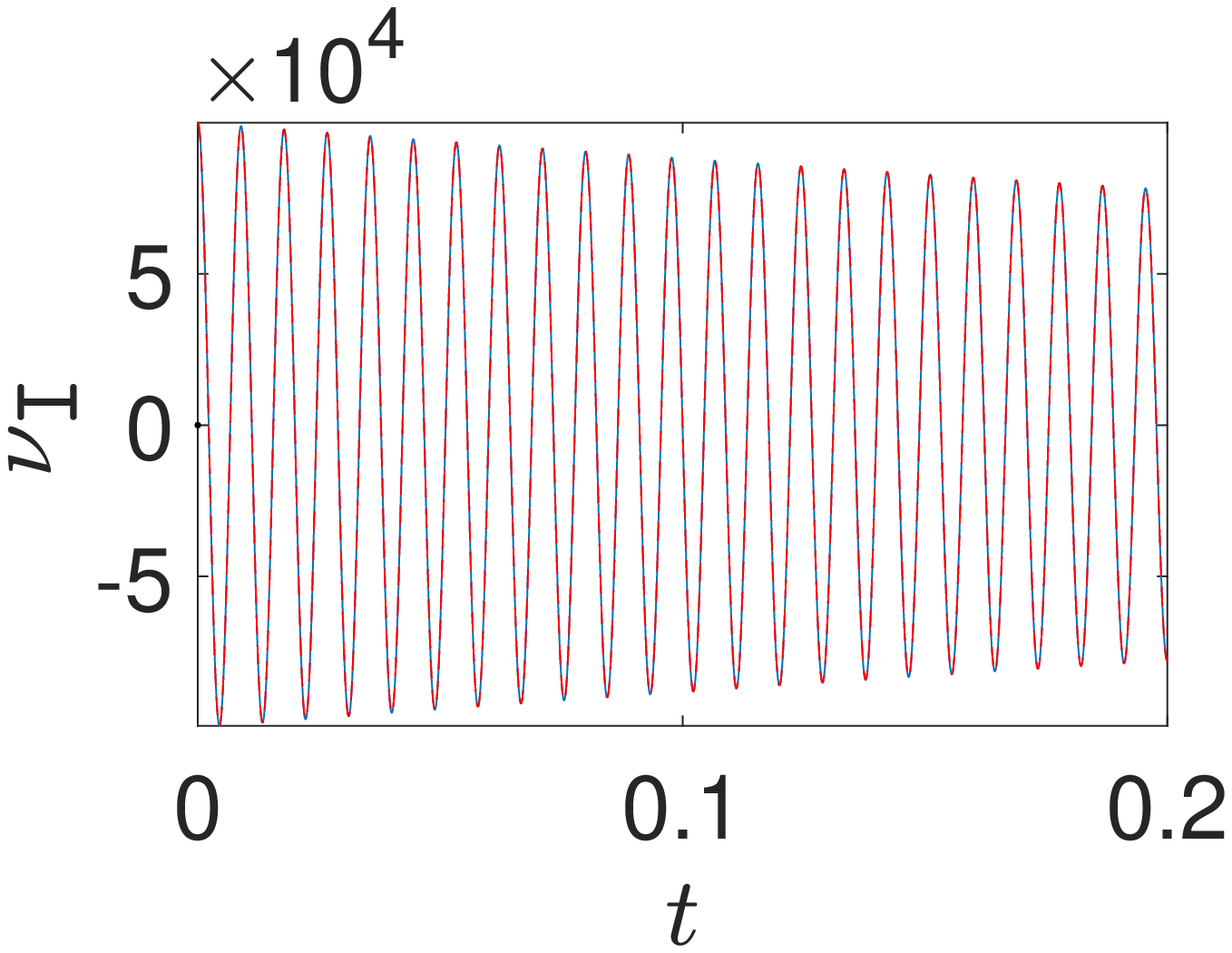}

\end{minipage}
}
\caption{Output of \exref{ex1}
\label{fg:outputRUNEX}}
\end{figure}

The second example is the Diode clipper circuit in \fgref{DiodeClipper}(a) from \cite{FalH16} (it can be found in audio-distortion devices). 
\begin{figure}[ht]
\subfigure[Diagram]{
{
\begin{circuitikz}[scale=0.65]
\draw (0,0) to [V, l_={\footnotesize $V$},*-o] (0,4);
\draw (0,0) -- (4,0) -- (6,0) -- (8,0);
\draw (8,0) to [I, l_={\footnotesize $I$},*-o] (8,4);
\draw (4,0) to [diode, l={\footnotesize $D_1$},*-o] (4,4);
\draw (6,4) to [diode, l={\footnotesize $D_2$},*-o] (6,0);
\draw (0,4) to [generic, , l_={\footnotesize $R$},*-o] (4,4);
\draw (4,4) -- (8,4);
\end{circuitikz}
}
}
\subfigure[Control output]{
\nc\mpw{.25\textwidth}

\includegraphics[width=\mpw]{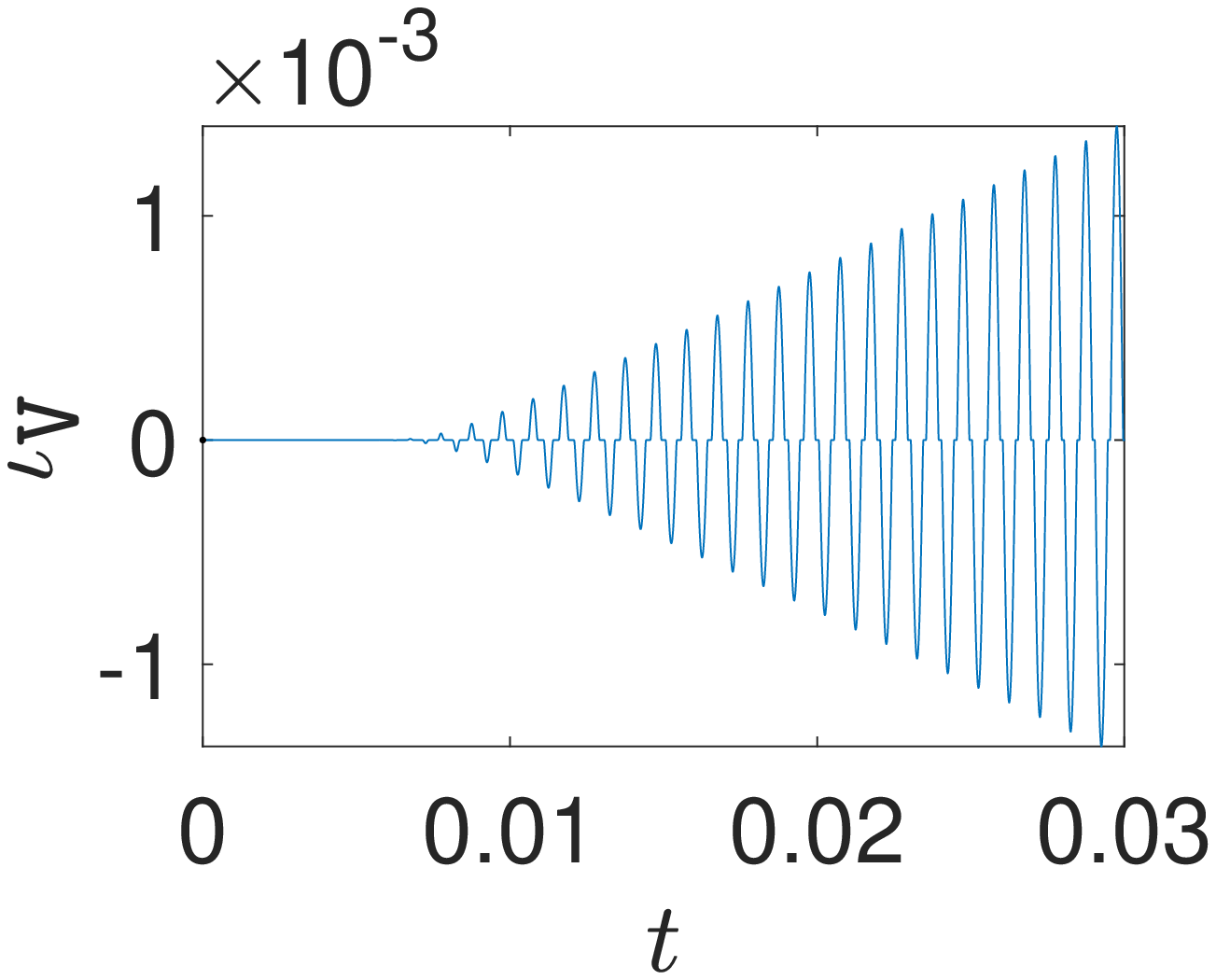}
\includegraphics[width=\mpw]{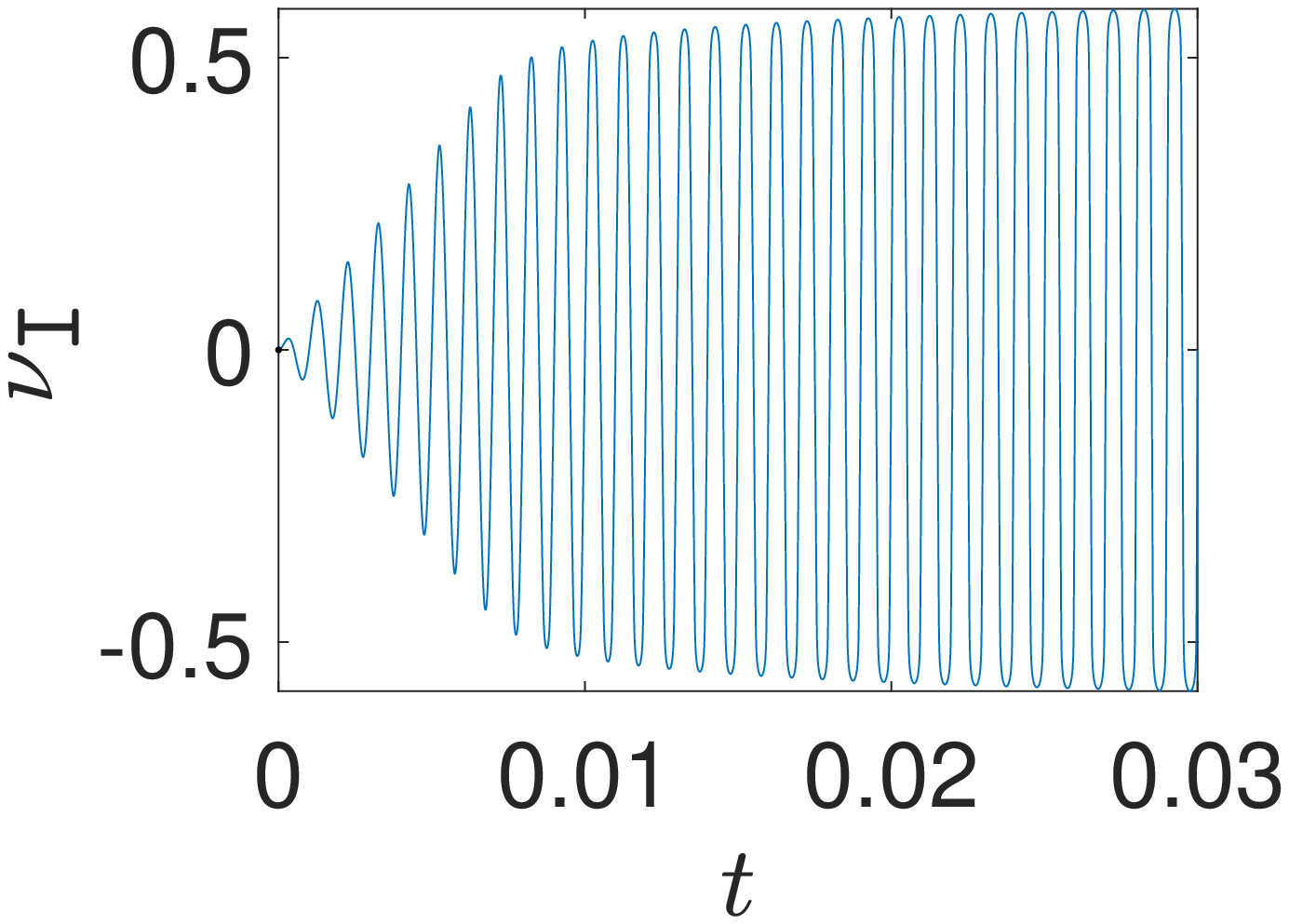}

 }
 \caption{Diode clipper\label{fg:DiodeClipper}}
\end{figure}
It contains one linear resistor ($R=1000$) and two diodes that are treated as nonlinear conductors (voltage-controlled dissipators) modeled by Shockley's equation, i.e., the diode current is given as a function of its branch voltage $\i(\v)=I_s\cdot(e^{\v/\v_T}-1)$, with typical parameters $I_s= 10^{-13}$ (the saturation current), and $\v_T= 0.025$ (the thermal voltage).\footnote{Standard notations for this diode model: the $T$ in $\v_T$ is about temperature, not a tree.}
The circuit contains a voltage source with $V(t)=(\widehat{V}\,t/t_{max})\,\sin(\omega t)$, with $\omega = 2\pi f$, a linearly increasing sinusoidal excitation with frequency $f=\SI{1}{\kilo\hertz}$ and maximum amplitude $\widehat{V}=2$ at $t_{\max}=0.03$, and an artificial current source $I(t)=0$ that acts as voltmeter to measure the output voltage.
The Model 2 formulation is integrated over the interval $[0, t_{\max}]$ using an initial guess $\widetilde x_0=[1,\dots,1]\+$. 
The control-output $y=[\i_\iV,\v_\iI]$ is plotted in \fgref{DiodeClipper}(b).

\matlab code for the two examples is given in \apref{matlabcode}.

\smallskip 
To show larger circuits being handled, we created an example where a graph is generated and populated randomly with circuit elements.
If it satisfies \Asmp1, \Asmp2 (checked during construction) its CpH Model 2 form is constructed and solved with random initial values.
\begin{figure}[ht]
\centering
\subfigure[]{
\begin{minipage}{.50\TW}
\includegraphics[width=\TW]{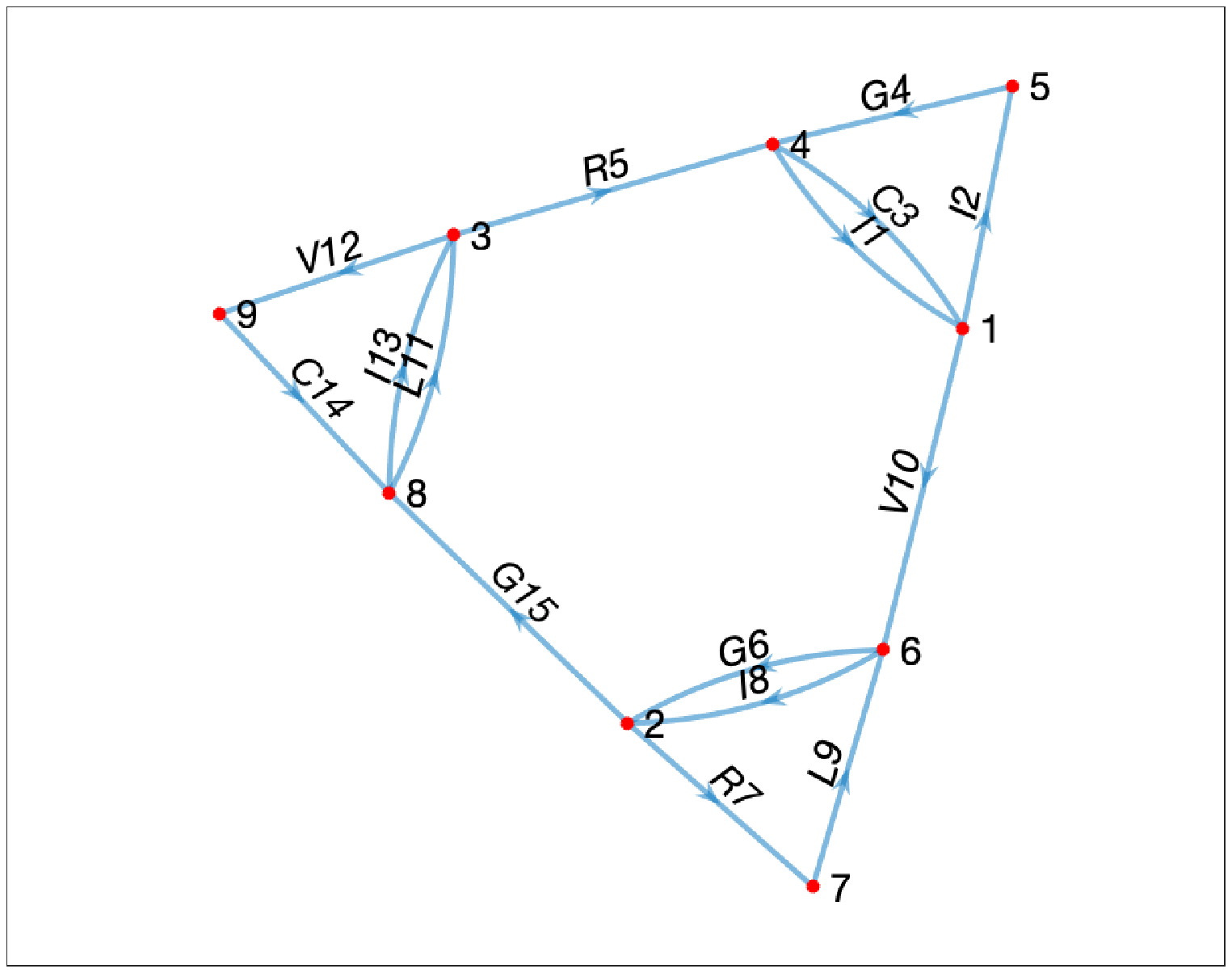}
\end{minipage}}
\hfill
\subfigure[]{
\begin{minipage}{.44\TW}
\begin{lstlisting}[basicstyle=\scriptsize\tt, keywords={}]

  T= 10 12  3 14  5  7  4  6
coT= ------------------------
 1 [ 0  0 -1  0  0  0  0  0 ] I|I1=@(t)sin(t)
 2 [ 0  0  1  0  0  0  1  0 ] I|I2=@(t)2*sin(3*t)
 8 [ 0  0  0  0  0  0  0 -1 ] I|I8=@(t)2*sin(3*t)
 9 [ 0  0  0  0  0  1  0  1 ] L|L9=0.99949
11 [ 0  1  0  1  0  0  0  0 ] L|L11=0.41452
13 [ 0  1  0  1  0  0  0  0 ] I|I13=@(t)2*sin(3*t)
15 [ 1 -1  1 -1  1  0  0  1 ] G|G15=0.10022
     V  V  C  C  R  R  G  G
     |  |  |  |  |  |  |  G6=0.98842
     |  |  |  |  |  |  G4=0.89849
     |  |  |  |  |  R7=0.53998
     |  |  |  |  R5=0.11816
     |  |  |  C14=0.8182
     |  |  C3=0.33083
     |  V12=@(t)cos(t)
     V10=@(t)-3*cos(2*t)
\end{lstlisting}
\end{minipage}}
\caption{Circuit with random elements.\label{fg:loop3ex}}
\end{figure}
An example with $b=15$ branches and $n=9$ nodes is in \fgref{loop3ex}(a);  \fgref{loop3ex}(b) displays the $F$ matrix, annotated with the specifications of the elements.
The $x$ solution components and the $y$ output components 
are plotted in \fgref{loop3xy}.

\begin{figure}[htb]
\nc\mpw{.24\TW}
\centering
\subfigure[Numerical solution]
{
%\nc\mpw{.32\TW}
\begin{minipage}{\TW}
\includegraphics[width=\mpw]{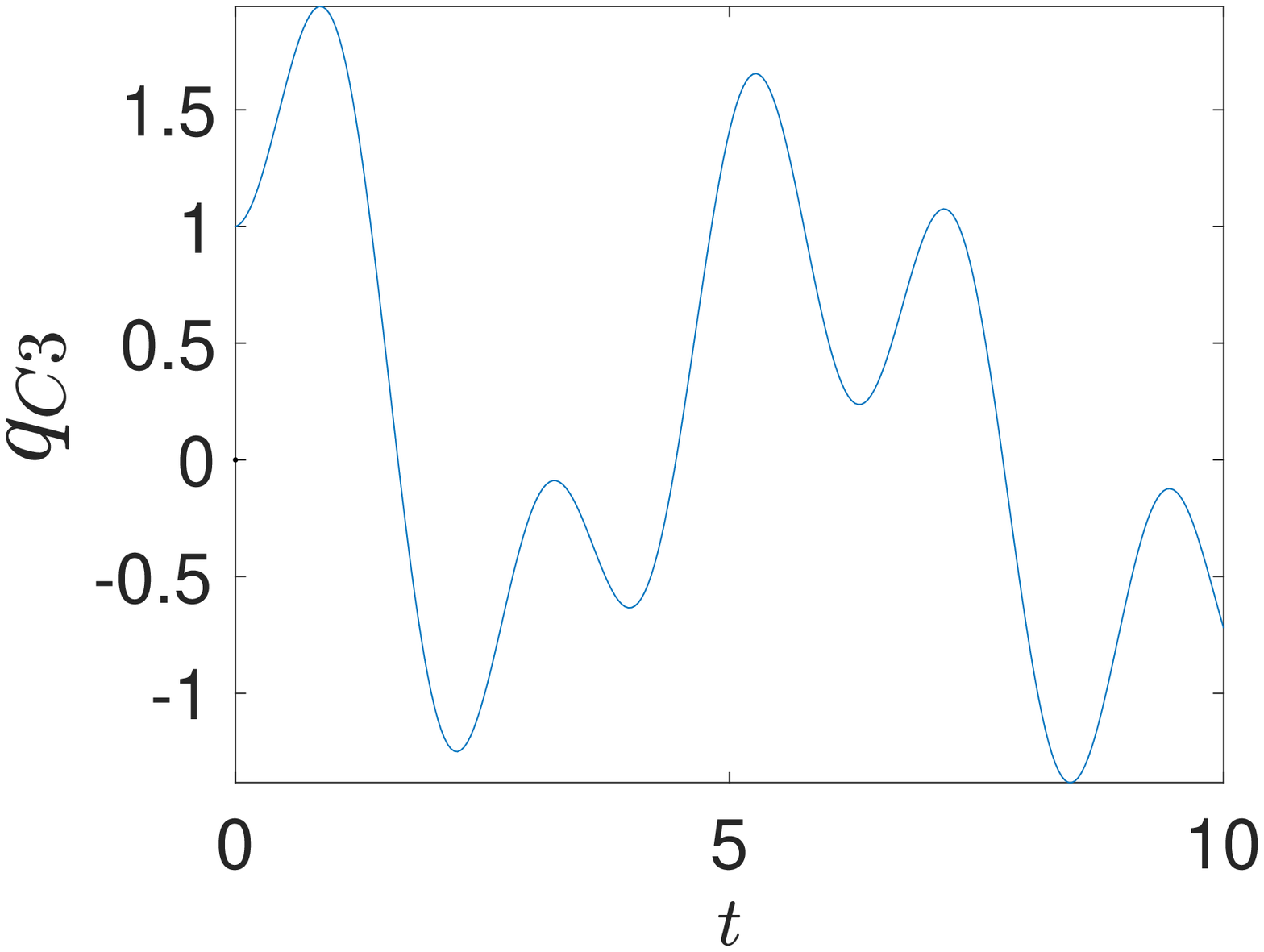}
\includegraphics[width=\mpw]{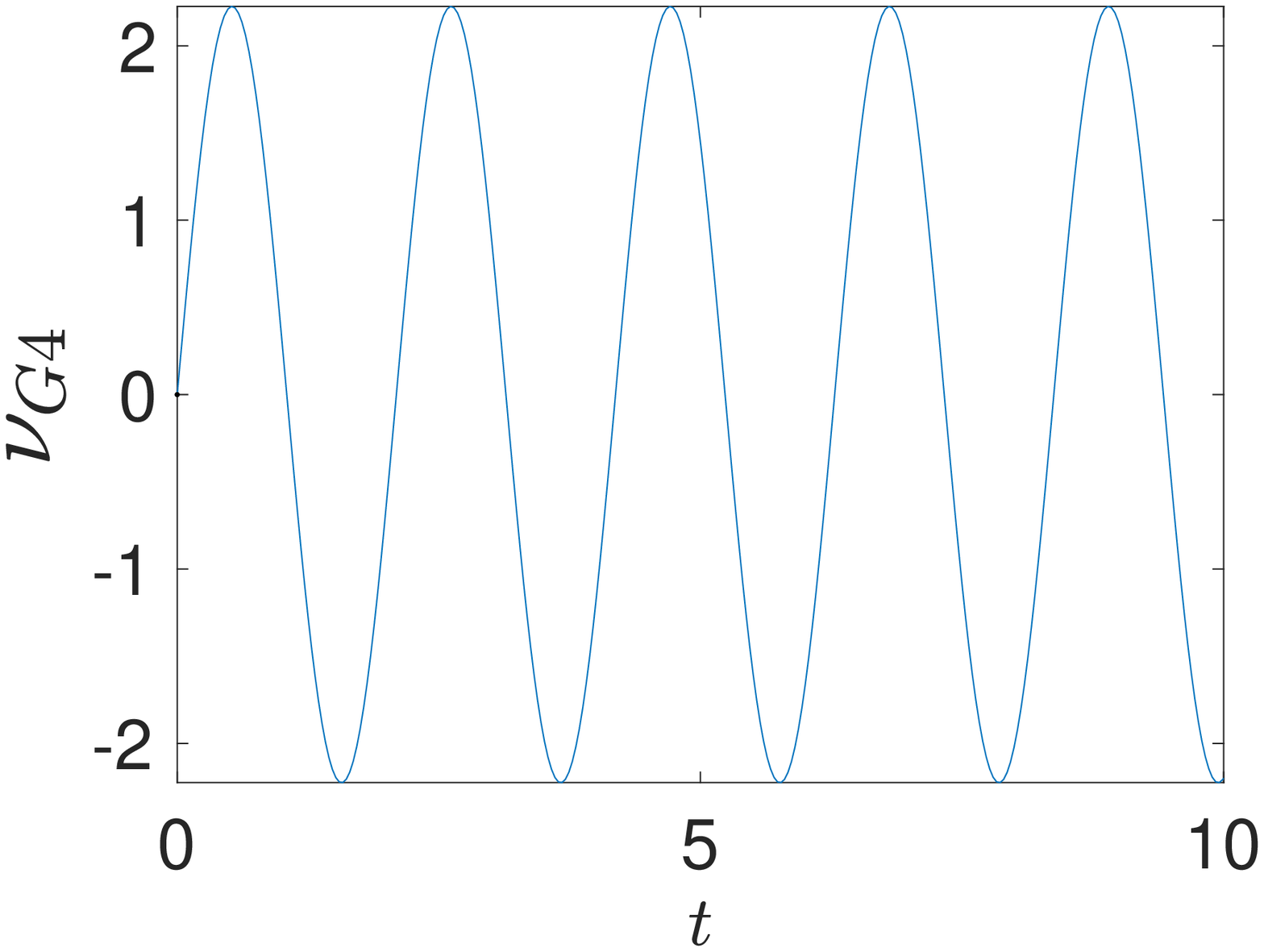}
\includegraphics[width=\mpw]{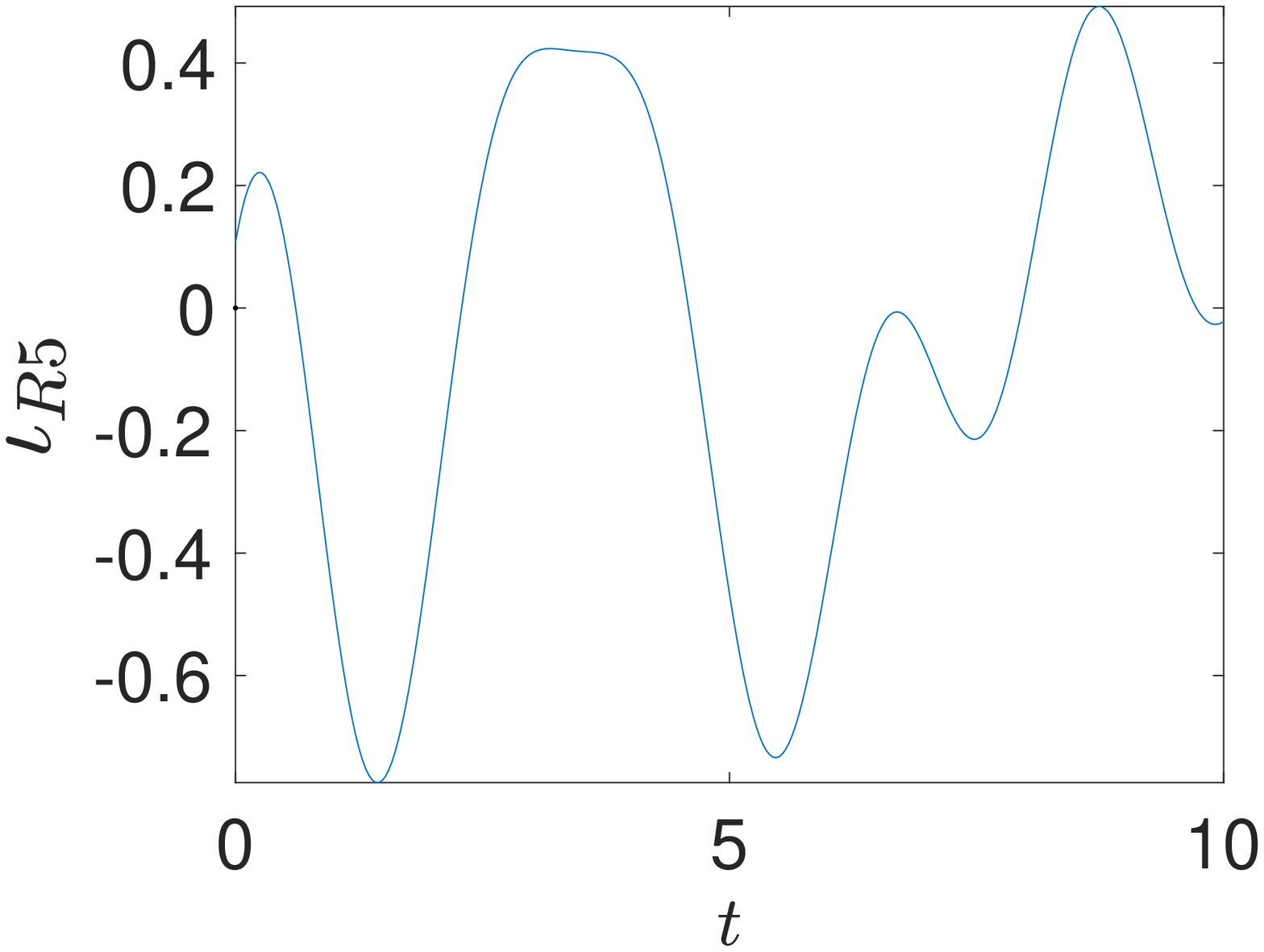}
\includegraphics[width=\mpw]{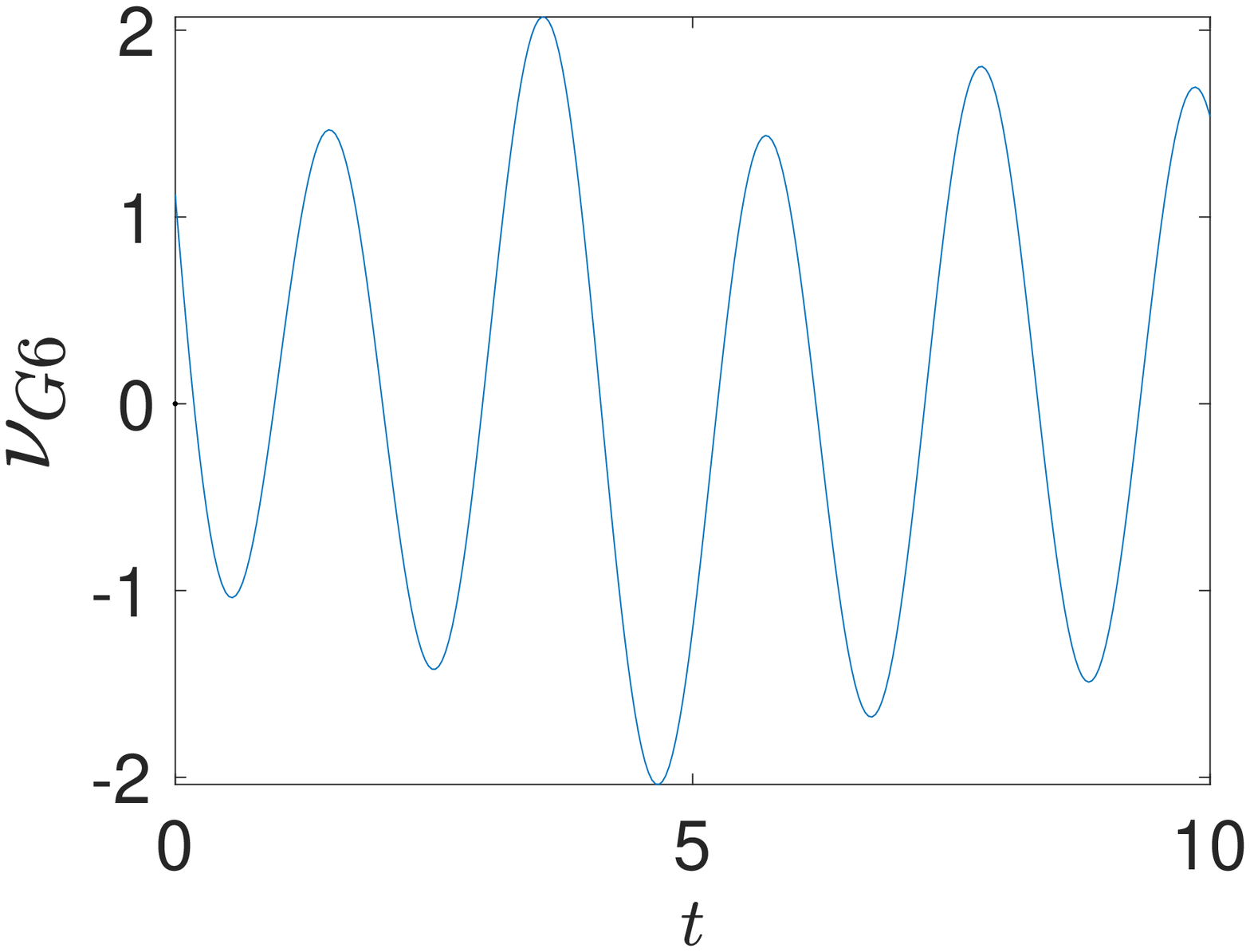}

\includegraphics[width=\mpw]{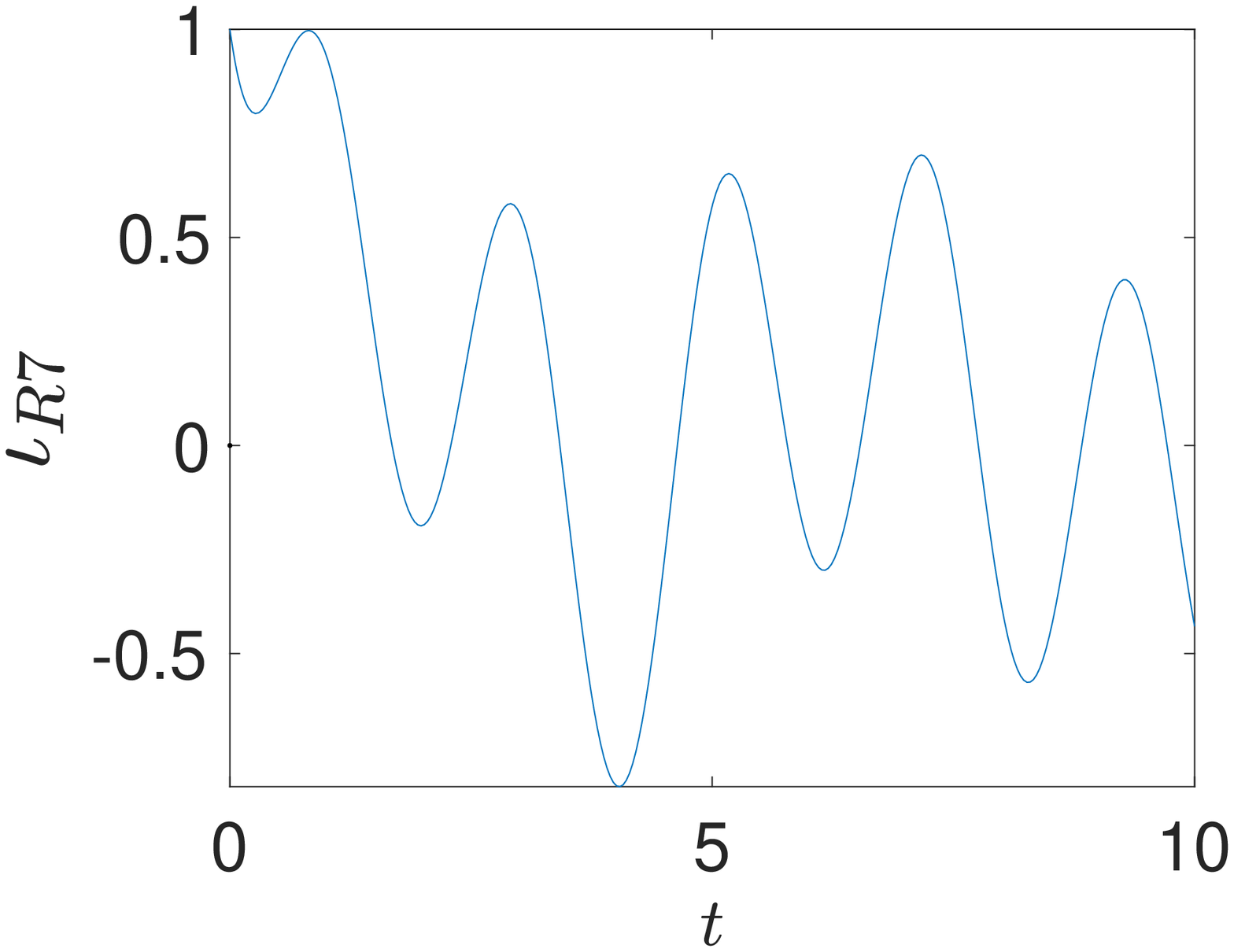}
\includegraphics[width=\mpw]{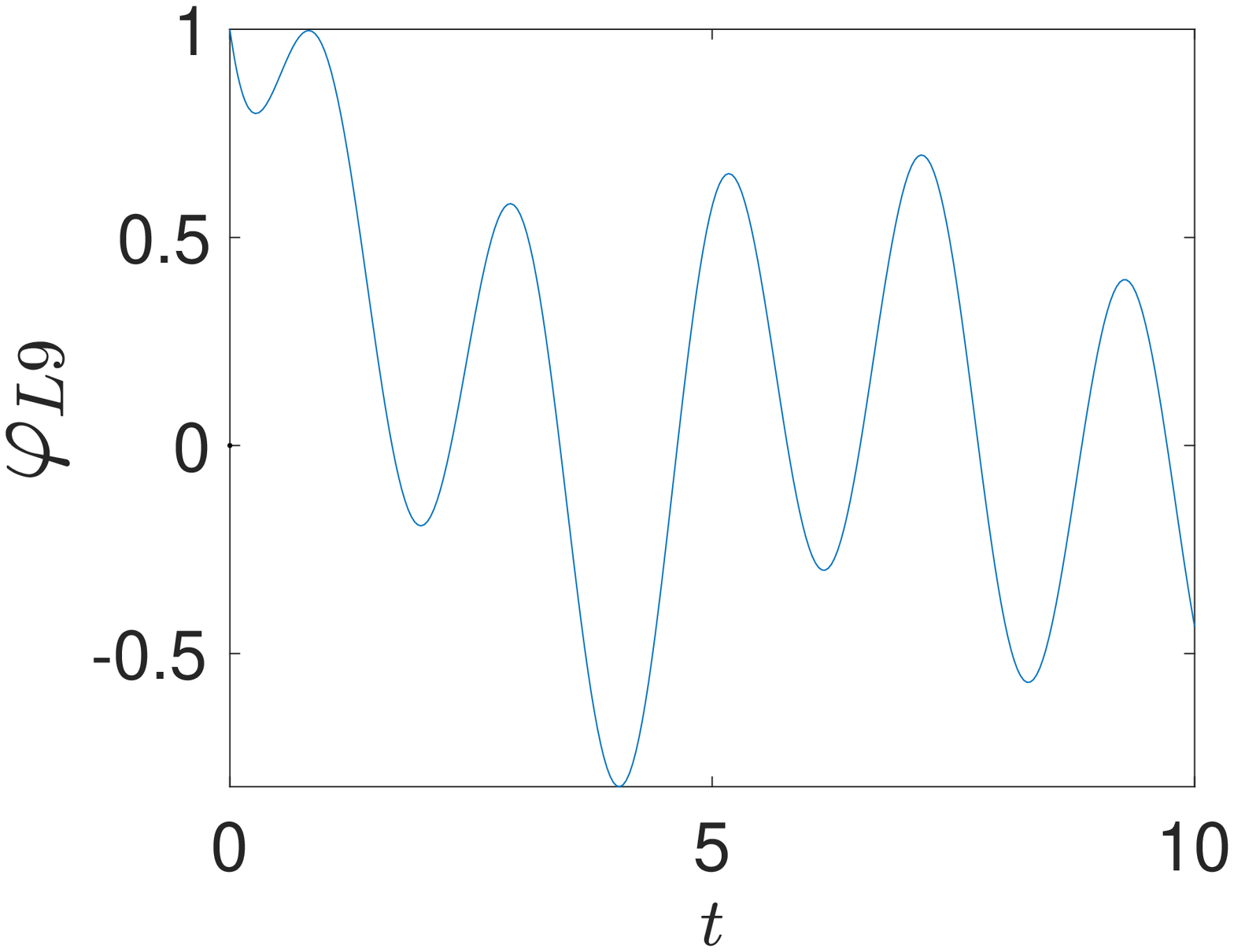}
\includegraphics[width=\mpw]{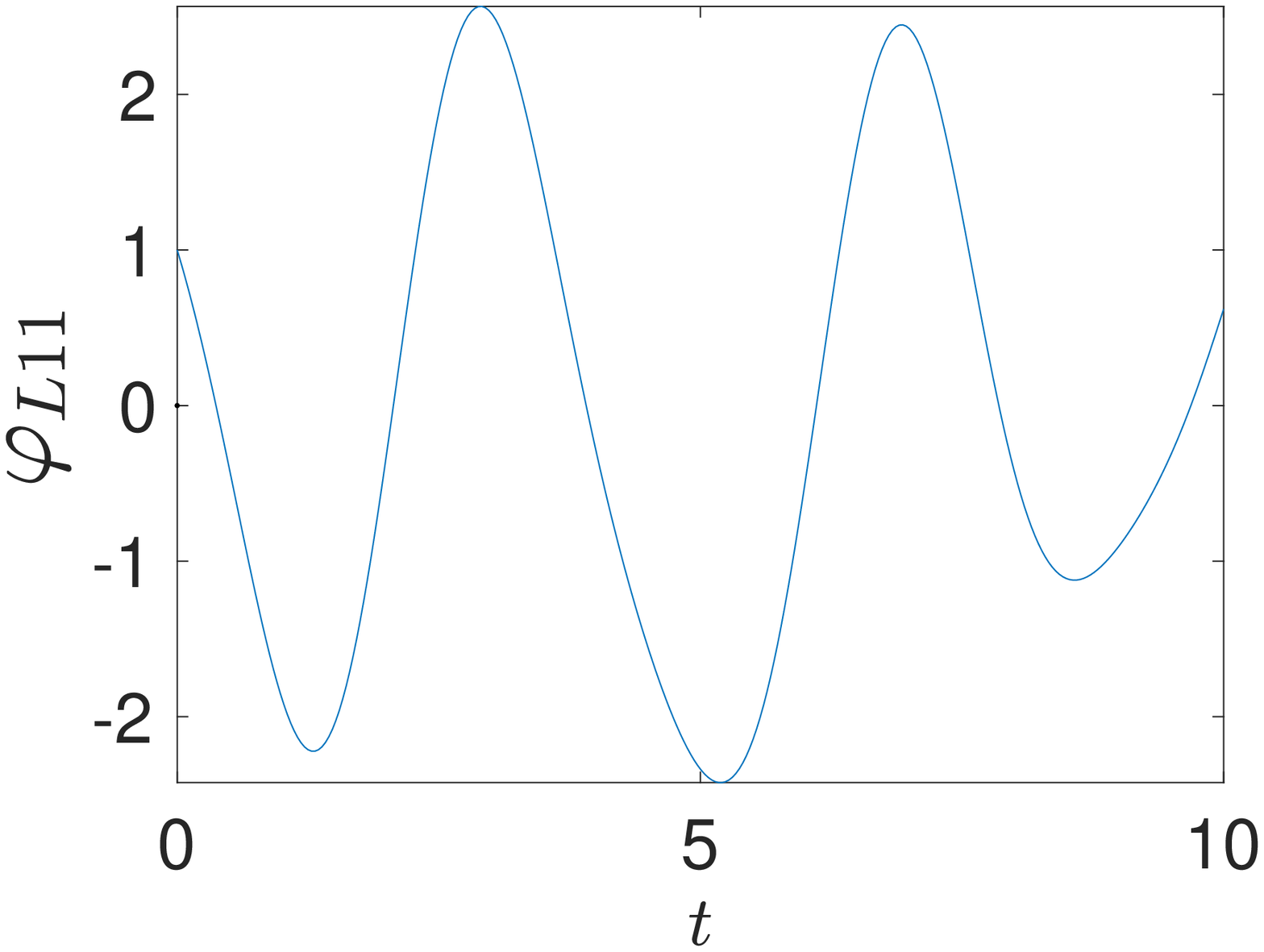}
\includegraphics[width=\mpw]{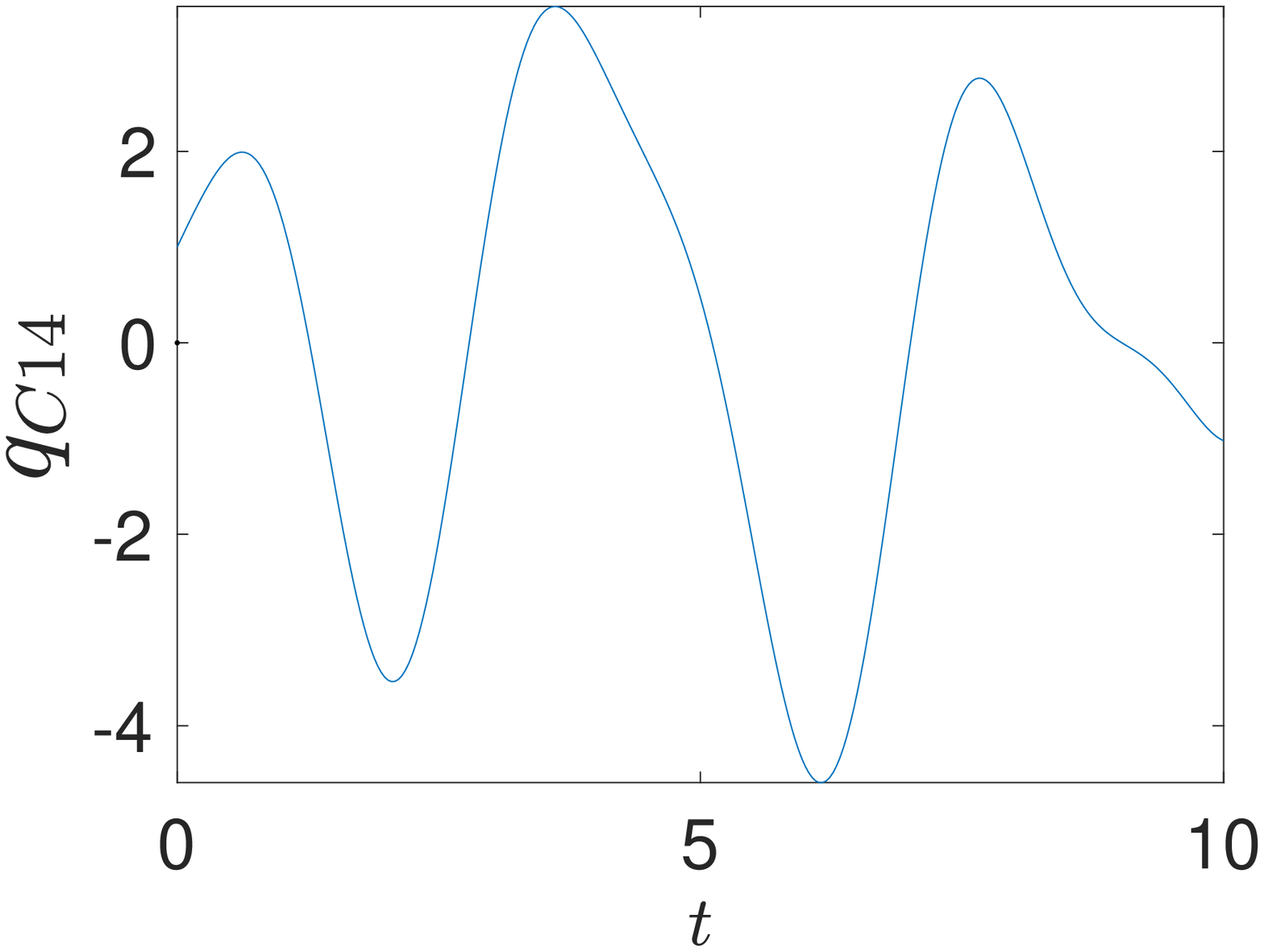}
\end{minipage}
}

\subfigure[Control output]{
\begin{minipage}{\TW}
\includegraphics[width=\mpw]{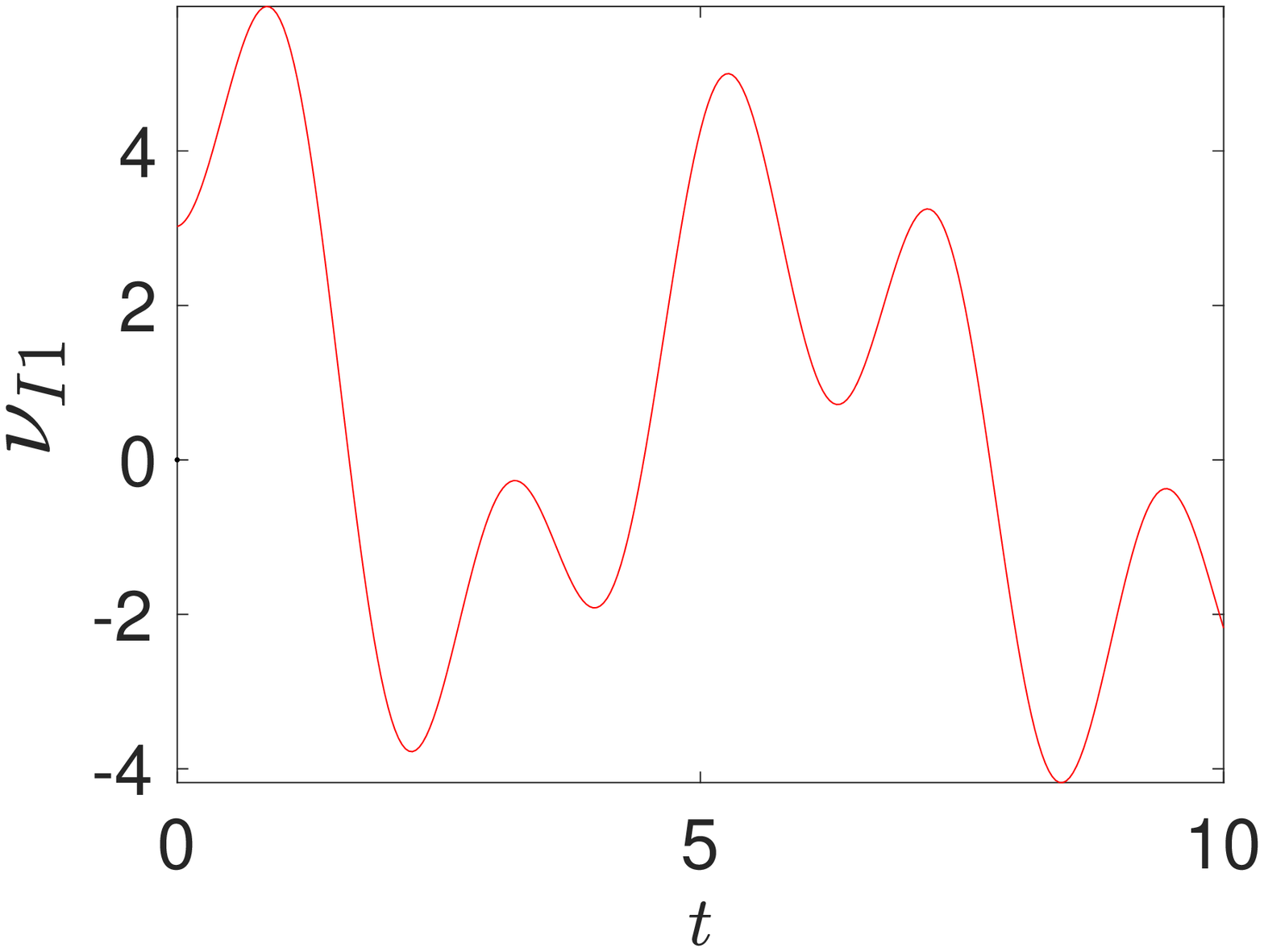}
\includegraphics[width=\mpw]{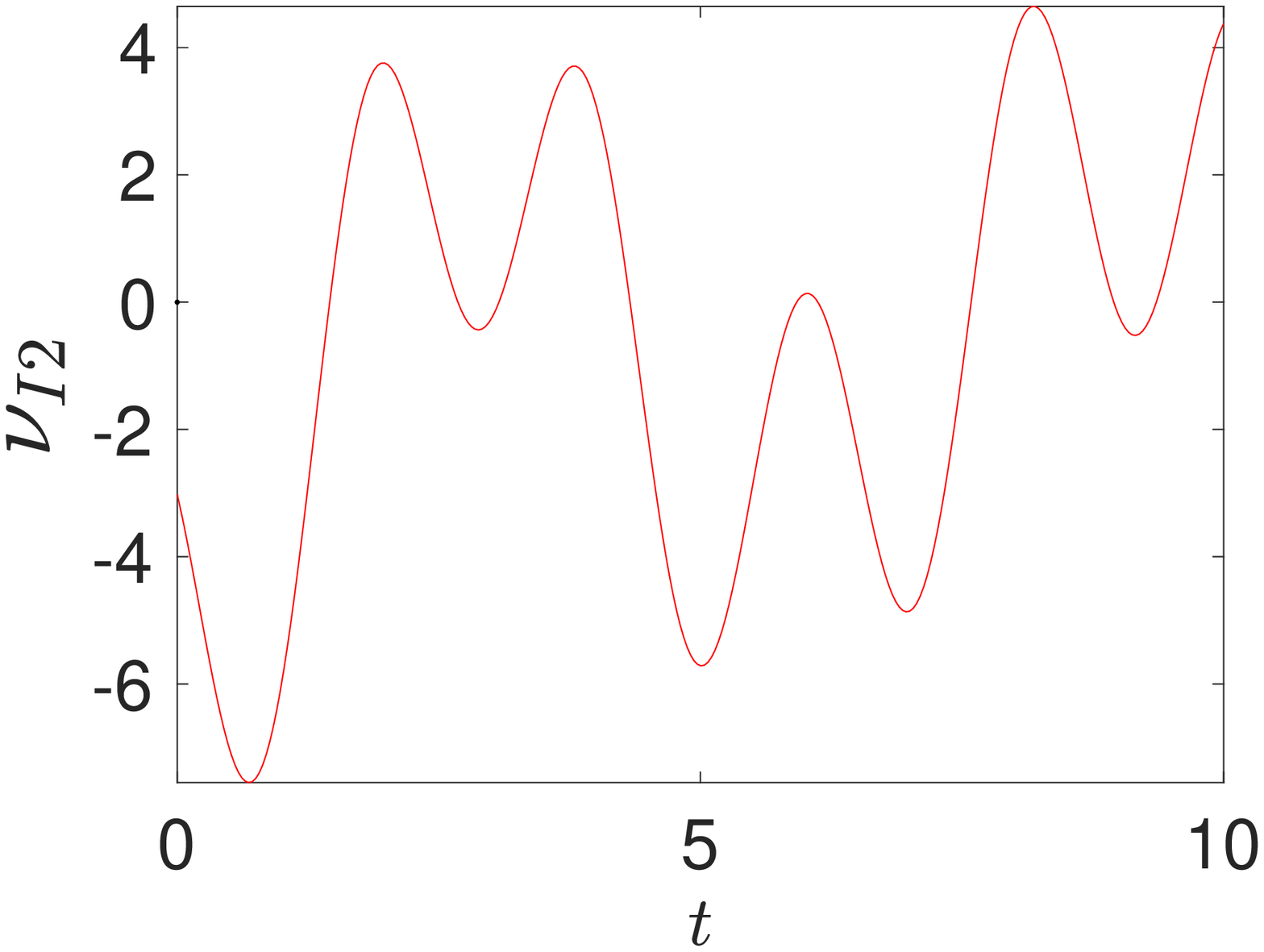}
\includegraphics[width=\mpw]{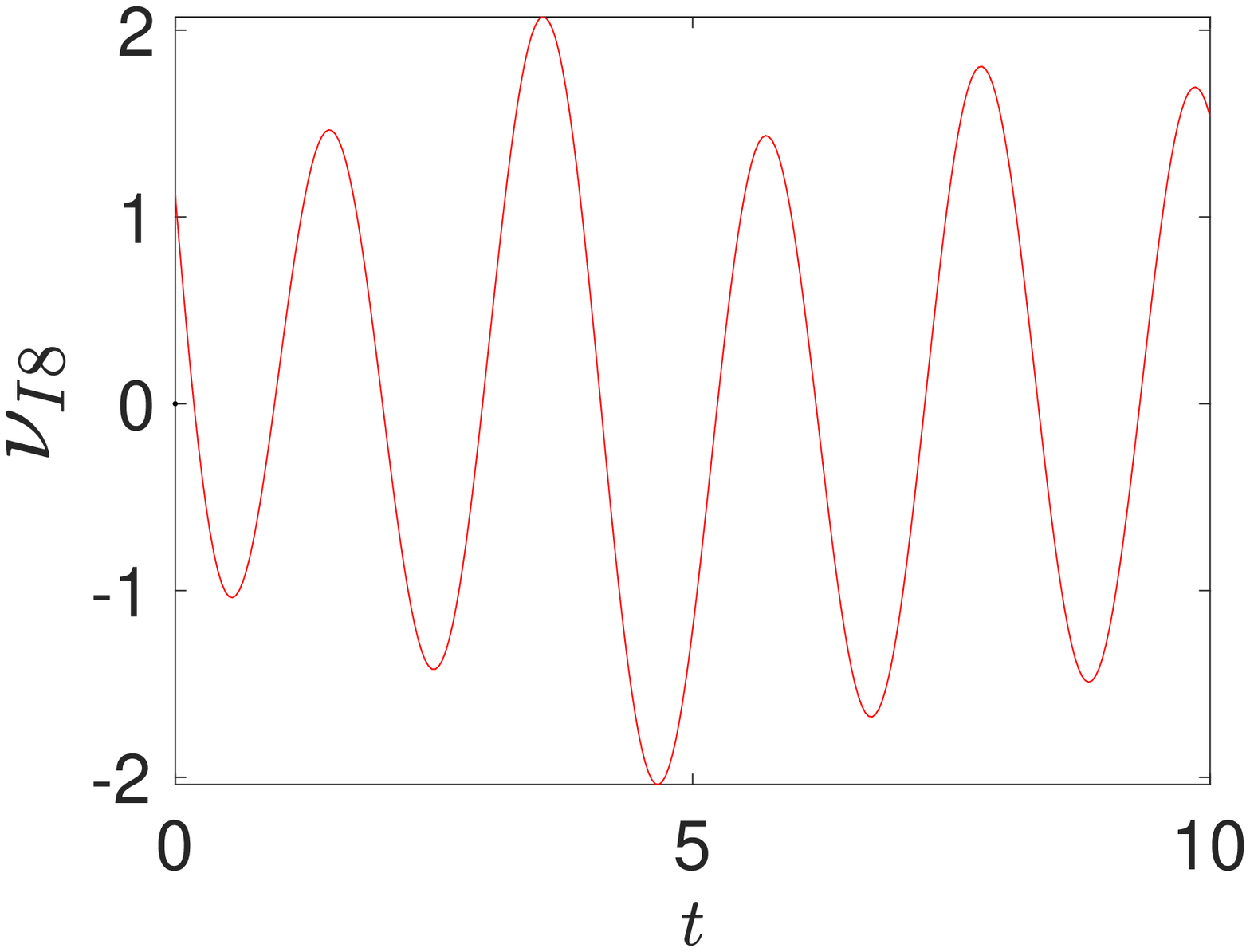}
\includegraphics[width=\mpw]{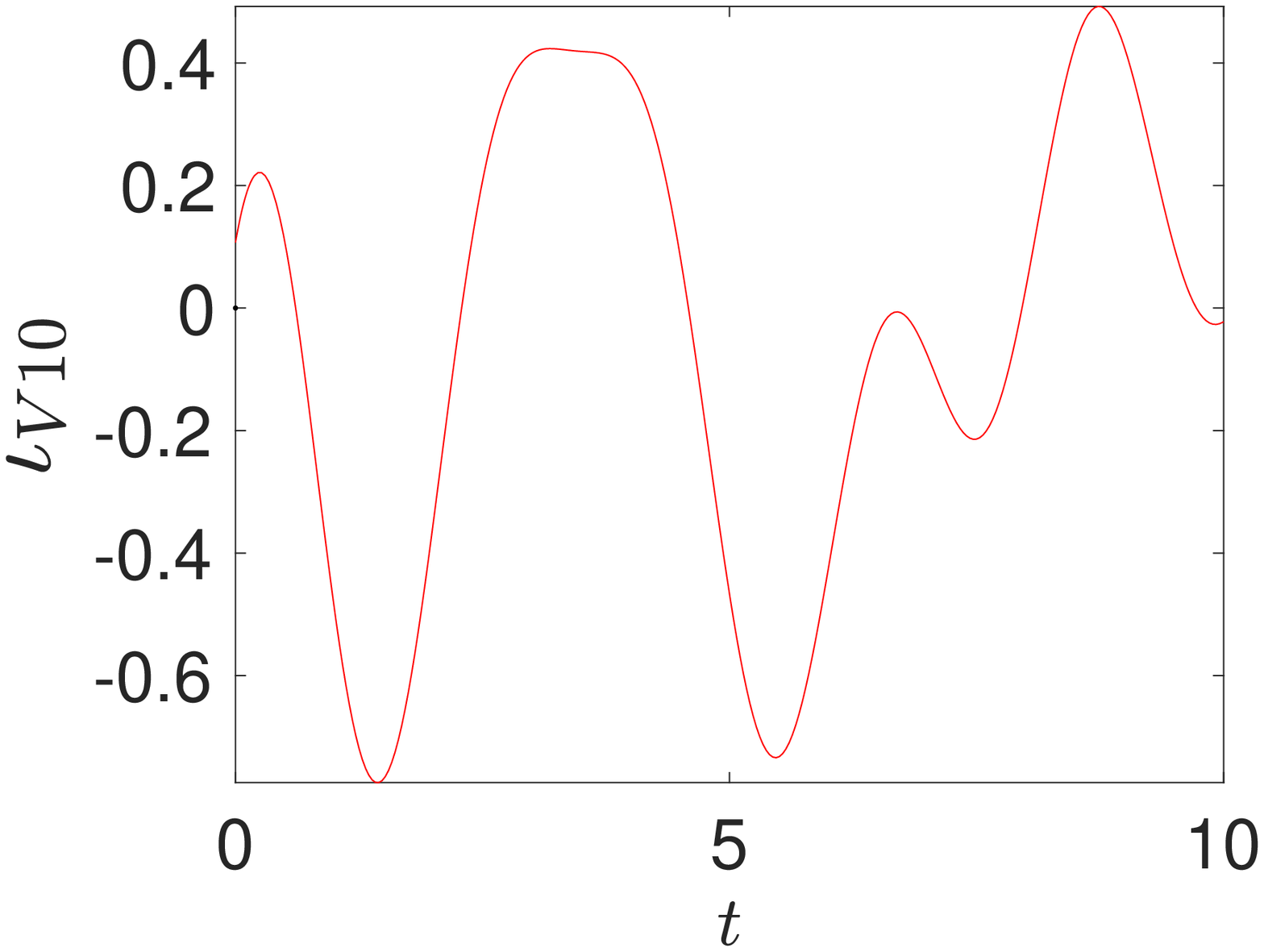}

\includegraphics[width=\mpw]{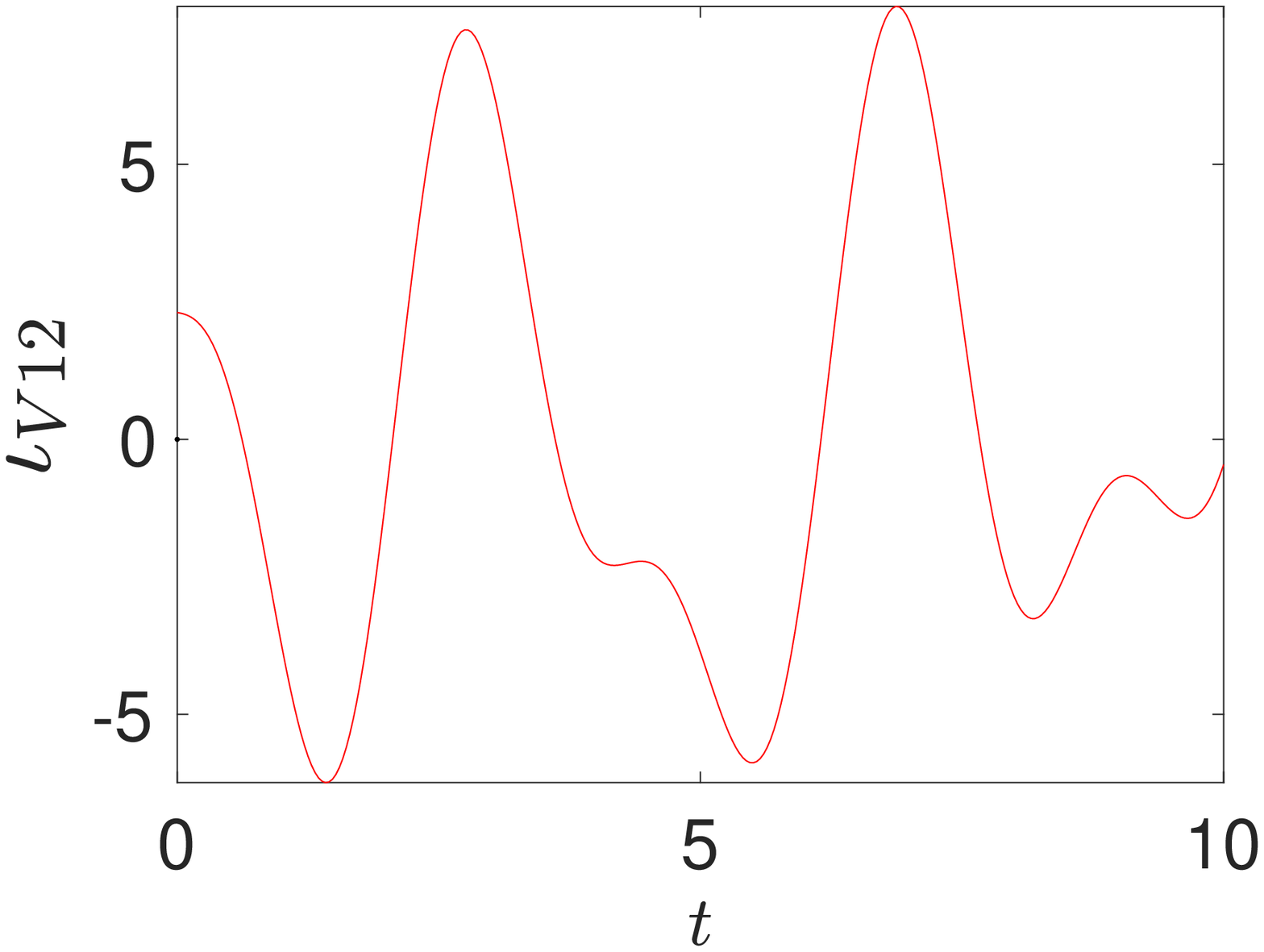}
\includegraphics[width=\mpw]{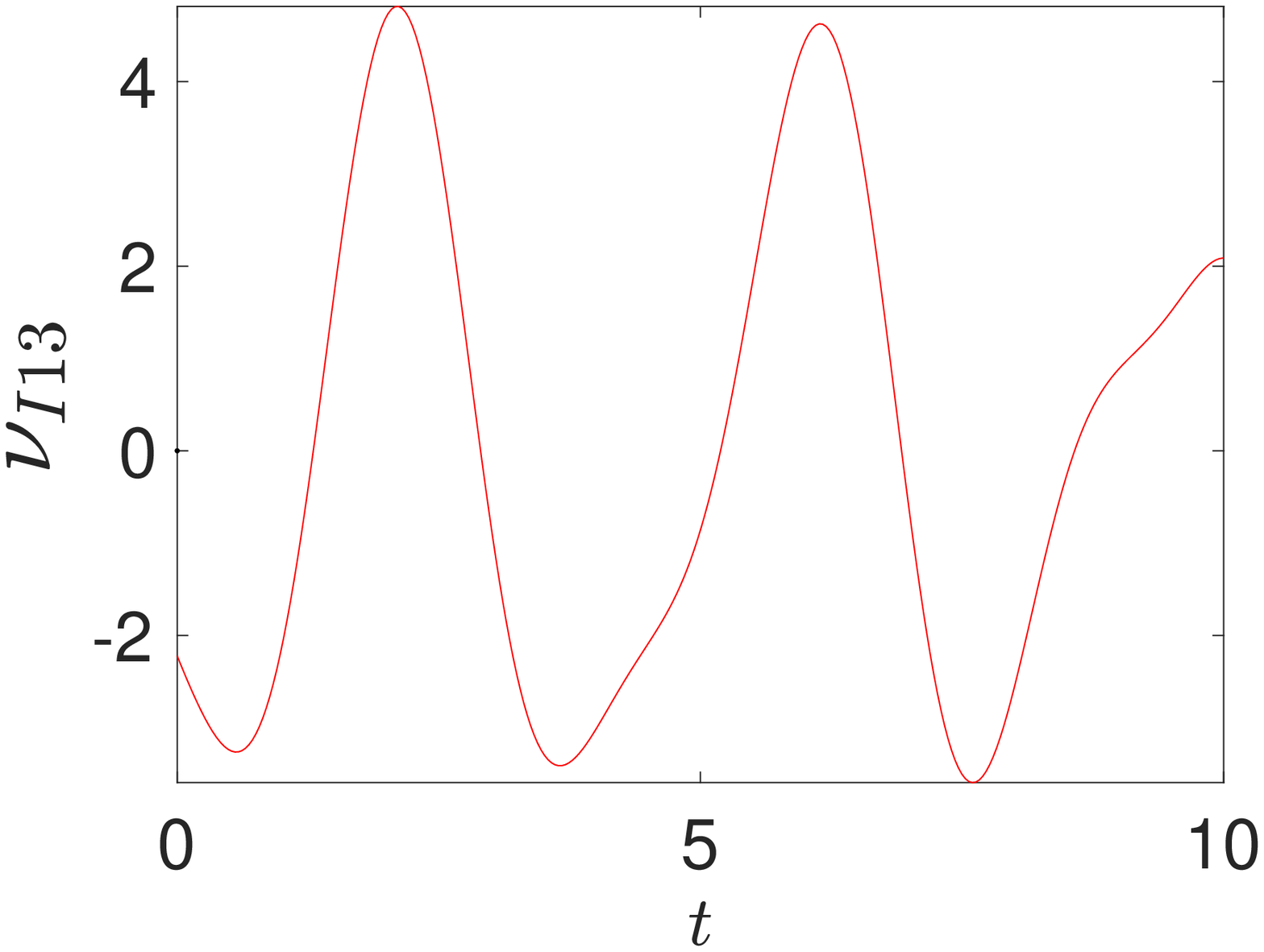}

\end{minipage}
}
\caption{Solution and output of the circuit from \fgref{loop3ex}.\label{fg:loop3xy}}

\end{figure}

\clearpage

\section{Conclusions}\label{sc:conclusions}

\paragraph{Old and new.}
The \pH concept is a powerful way to bring theory nearer physics, and suggest new computational approaches.
SA-amenability is a powerful concept in another way, making possible, or cheap to implement, numerical methods that without it are impossible or expensive.
Our contribution is to show that a particular mix of---mostly well known---ingredients gives a DAE that is both \pH and SA-amenable.
That it has index $\le1$ was an unexpected bonus, which we suspect comes from obeying the pH injunction that the Hamiltonian be an algebraic function of the state variables.
Writing Kirchhoff's laws in terms of the $F$ matrix goes back at least to \cite{Bashkow57,Bro63}, where also the zero blocks of $F$ are exploited though less transparently than here.
The existence of a normal tree is shown in \cite{Bro63}.
Our main proof is a distillation of that in \cite{Sch18} that the BOM circuit model's DAE is SA-amenable.
However, the BOM DAE is not pH and can have index 2.
Applying pH theory to electrical circuits is a topic of current research, see e.g.\ \cite{GerHRvdS21,GueBJR20}.
In \cite{GerHRvdS21}, selection of a ground node---not needed by CpH---plays a crucial role in defining the Dirac structure; and there is no discussion of the DAE index.
The aim of \cite{GueBJR20} is to put the MNA equations in port-Hamiltonian form; this can lead to a DAE of index 2.
Neither work discusses SA-amenability.
%\smallskip
To the best of our knowledge the CpH method construction, and SA-amenability and related proofs, are new.
So also appears to be our proof of \thref{normalT}, and the resulting direct counts of the number of linearly independent $\iCc\iVv$-cycles (including $\iCc$-cycles) $n_\iC=n_\iCc-r_{\iVv \iCc} + r_\iVv$, and similarly of $\iLl\iIi$-cutsets $n_\il=r_{\iVv \iCc \iDd \iLl} - r_{\iVv \iCc \iDd}$.
Our CICs method in \ssrf{getCICs} is new, but related to one in \cite{KroeMG92}.
Our method is easy to automate: as said in \ssrf{matlabimplem}, from a \li{CpH} object representing a circuit we generate code simulating that circuit, in the form of a \li{daefcn} function accepted by an index-1 solver such as \li{ode15i}.
Other advantages of being SA-amenable is the systematic reduction to an ODE \ssref{RedODE}, the easy eigenvalue analyis for LTI systems \ssref{LTI}, or the use of Mattsson-S\"oderlind dummy derivatives \cite{MatS93}.

\paragraph{Future work.}

A key feature of \pH theory is that a ``Dirac join'' of several pH systems---joining their Dirac structures in a way that makes a Dirac structure on the union of the systems---is again a pH system, whose Hamiltonian is the sum of the Hamiltonians of the parts.
Our \matlab implementation does not yet have this property and it is a priority to provide it.
This would follow in the steps of the \pH \modelica system recently reported in \cite{MarZY20}.
At a more modest level, our implementation currently only supports independent, possibly nonlinear, circuit elements.
We are working on upgrading it to handle the full generality of models in this paper, i.e.\ circuits containing dependent elements forming multi-terminal components (e.g.\ transistors) and other more complex circuit configurations (e.g.\ controlled sources).
In a future paper we aim to give proofs and algorithm details for the method in \ssrf{getCICs}.

\section*{Acknowledgments}
We thank Prof Volker Mehrmann and Prof Caren Tischendorf for valuable advice. 
We thank Jacek Kierzenka of the MathWorks for writing us an adapted version of the solver \li{ode15i} that supports writing a DAE in control input-output form \rf{ctrlIOform4}.

\appendix

\section{Terminology of circuit graphs vs other kinds of graph}\label{sc:graphterms}

We summarise where it is necessary or convenient for circuit graph terms to differ from what is ``often'' the case, i.e.\ common in other graph theory areas.
Let $\G=(V,\iE,\alpha)$ be a circuit graph.
\begin{compactitem}
  \item Riaza \cite[\S5.1.2]{Ria08} reserves the terms {\em branch}, {\em node} for elements of a circuit, and {\em edge}, {\em vertex} for the corresponding elements of the circuit graph. We keep this distinction.
  \item A circuit {\em graph} allows several edges between a given pair of vertices and forbids edges from a vertex to itself.
  Often the opposite of one or both of these is true.
  \item Often the {\em induced subgraph} of a subset of edges discards any vertex that is not the end of some edge in the subset.
  Here it is more useful to keep all vertices by default.
  \item Often a {\em spanning} subgraph of $\G=(V,\iE,\alpha)$ simply means one whose vertex set is all of $V$.
  Here, being connected is essential.
  \item A {\em link} is often called a chord; but both link and chord have multiple meanings in graph literature.
  \item Giving $\G$ an {\em orientation} is just a housekeeping aid for using Kirchoff's laws.
  It doesn't change the topology and is thus quite different from making $\G$ a {\em directed graph}.
  The circuit physics is the same whatever orientation is chosen.
\end{compactitem}

\section{The Signature Method \texorpdfstring{{\boldmath
(\Smethod)}}{(Sigma-method)} for DAEs}\label{sc:sigmethod}

Our preferred SA method for DAEs is the {\em Signature method} ({\em $\Sigma$-method}) \cite{Pry01}. 
The $\Sigma$-method can be applied to regular nonlinear DAEs of arbitrary high order $k$ of the form 
\begin{align*}%\label{eq:DAE_k}
f(t,x,\dot x,\dots,x^{(k)})=0,
\end{align*}
with $f:\Tivl\x\R^n\x\dots\x\R^n\to\R^n$ sufficiently smooth. 

The $\Sigma$-method consists of the following steps:
\begin{compactenum}
\item Building the {\em signature matrix} $\Sigma=[\sigma_{ij}]_{i,j=1,\dots,n}$ with
 $$\sigma_{ij}=\begin{cases}\text{highest order of derivative of $x_j$ in $f_i$,} \\ -\infty \text{ if $x_j$ does not occur in $f_i$.}\end{cases}$$

\item Finding a {\em highest value transversal (HVT)} of $\Sigma$, i.e.,
 a set 
 ${\cal T}=\{(1,j_1),(2,j_2),\dots,(n,j_n)\}$,
where  $(j_1,\dots,j_n)$ is a permutation of $(1,\dots,n)$,
with maximal  {\em value} $\text{Val}({\cal T}):=\sum_{(i,j)\in{\cal T}}\sigma_{ij}.$

\item Computing {\em offsets} $c=[c_i]_{i=1,\dots,n}$ and $d=[d_j]_{j=1,\dots,n}$ with $c_i \ge0$, $d_j \ge0$  such that
\begin{align}\label{eq:offsets}
  \text{$d_j-c_i\ge\sij$ for all $i,j=1,\dots,n$, with equality for $(i,j)\in\calT$}.
\end{align}

\item Forming the \sysJ $\_J = \bmx{\_J_{ij}}_{i,j=1}^n$ where (the two formulae are equivalent)
  \begin{align}\label{eq:sysJ}
    {\_J}_{ij}:=
    \begin{cases}
    \ds\tdbd{f_i}{x_j^{(\sigma_{ij})}} &\text{ if } d_j-c_i=\sigma_{ij}, \\
      0 &\text{ otherwise}
    \end{cases}
  \quad = \quad
  \begin{cases}
    \ds\tdbd{f_i}{x_j^{(d_j-c_i)}} &\text{if $d_j-c_i\ge0$}, \\
     0 &\text{otherwise}
  \end{cases}
  \end{align}
\item Form the reduced derivative array ${\cal F}(t,{\cal X})=0$ where
\[{\cal X}=[x_1,\dot x_1,\dots,x_1^{(d_1)}, \dots, x_n,\dot x_n,\dots,x_n^{(d_n)}]\+,\
  {\cal F}=[f_1,\dot f_1,\dots,f_1^{(c_1)}, \dots, f_n,\dot f_n,\dots,f_n^{(c_n)}]\+,
\]
\item  Success check: if ${\cal F}(t,{\cal X})=0$, considered locally as an algebraic system of $M=n+\sum_{i=1}^n c_i$ equations in $N=n+\sum_{j=1}^n d_j$ unknowns, has a solution $(t^*,{\cal X}^*)\in\Tivl\x\R^N$ where $\_J$ is nonsingular, then $(t^*,{\cal X}^*)$ is a consistent point, the method ``succeeds'', and the DAE is SA-amenable. 
\end{compactenum}

In case of success the {\em structural index}
\begin{align}\label{eq:sindex}
\nu_s:=\max_i c_i + (\text{$1$ if some $d_j=0$, $0$ otherwise})
\end{align}
is $\ge$ the differentiation index $\nu_d$.
Also $\text{Val}(\Sigma)$, defined as the value of the highest value transversal $\calT$ and also equal to $\sum_j d_j - \sum_i c_i$, equals the number of degrees of freedom of the system.

We have the {\em standard solution scheme (SSS)}: starting with $k = -\max d_j$, in order of increasing stage number $k$, 
\begin{align}\label{eq:SSS}
\begin{aligned}
&  \text{solve the equations} 
 & &\der{f_i}{k+c_i}=0
  &&\text{for those $i$ such that $k+c_i\ge0$} \\
&  \text{for the unknowns}
  &&\der{x_j}{k+d_j}
  &&\text{for those $j$ such that $k+d_j\ge0$}
\end{aligned}\biggr\}\;.
\end{align}
One interpretation of this is as a numerical scheme that first finds a consistent point and then expands the DAE solution through this point as a Taylor series; then, $k$ goes up to some Taylor order $k=k_{\max}>0$.
However the meaning used in this paper is that we go up to $k=0$: stages up to $k=-1$ are ``numerical'', and a consistent point is a set of values that solve these equations; while $k=0$ is ``symbolic''---a set of equations by which (using the Implicit Function Theorem) the highest derivatives $\der{x_j}{d_j}$ are proved to be functions of lower derivatives, thus defining an implicit ODE.
Let $\_J_k$ be the Jacobian of the equations in \rf{SSS} at stage $k$.
\begin{theorem}\label{th:Jrank}
(\cite[Proposition 4.1]{Pry01})
  For an SA-amenable DAE, each $\_J_k$ is of full row rank, and for all $k\ge0$ is equal to the \sysJ $\_J$ (hence is nonsingular).
\end{theorem}

\section{\PH dynamics}\label{sc:pHdyn}

\subsection{\PH system and Dirac structure}\label{ss:pHsystem}

A physical system is modelled by flows of energy within it, these flows being represented by {\em ports}: (flow, effort) pairs of variables whose product has the dimension of power.

The following paradigm is taken from \cite[\S2.6]{vdSchJ14}.
The $\Flo_\bullet$ are linear spaces and $\Eff_\bullet=\Flo_\bullet^*$ their duals.
In our context they are some $\R^d$, so the duality inner product $\<\flo,\eff>=\(\flo,\eff)$.
In general $\ZS$ is a manifold, with $T_z\ZS$ its {\em tangent space} at the state $z\in\ZS$, and $T_z^*\ZS$ the corresponding {\em co-tangent space}.
In our context $\ZS$ is some $\R^d$, so $T_z\ZS$ and $T_z^*\ZS$ can also be regarded as $\R^d$, again with duality inner product.
\medskip

A \pH system gives mathematical form to \fgref{pHgenl}. It is defined by
\begin{compactenum}[(pH1)]
  \item A state space $\ZS$ and a Hamiltonian $\Ham : \ZS \to \R$, defining energy-storage, corresponding to an energy-storing port $(\flo_\iStor,\eff_\iStor)\in T_z\ZS\x T_z^*\ZS$.
  \item A dissipative structure
  \begin{align}\label{eq:pHdyn2}
    \Res &\subset \Flo_\Res \x \Eff_\Res &\text{having $\<\flo,\eff>\le0$ for all $(\flo,\eff)\in\Res$},
  \end{align}
  corresponding to an energy-dissipating port $(\flo_\Res, \eff_\Res) \in \Flo_\Res \x \Eff_\Res$.
  \item An external port $(\flo_\iPort, \eff_\iPort) \in \Flo_\iPort\x\Eff_\iPort$.
  \item A Dirac structure linking (pH1)--(pH3),
  \begin{align}\label{eq:pHdyn3}
    \Dir \;&\subset\; T_z\ZS \x T_z^*\ZS \,\x\, \Flo_\Res\x\Eff_\Res \,\x\, \Flo_\iPort\x\Eff_\iPort.
  \end{align}
\end{compactenum}

\begin{definition}
\cite[\S2.2]{vdSchJ14}
Consider a finite-dimensional linear space $\Flo$ with dual $\Eff =\Flo^*$, 
and duality inner product $\<.,.>:\Flo \x \Eff\to\R$.
 A linear subspace $\Dir \subset \Flo \x \Eff$ is a {\em Dirac structure} if $\<\flo,\eff>= 0$, for all $(\flo,\eff) \in \Dir$, and ${\rm dim}\, \Dir = {\rm dim}\, \Flo$.
\end{definition}
Here, using the natural identification $(V\x V')^* = V^*\x V'^*$ for linear spaces $\iV,V'$, the space on the right of \rf{pHdyn3} is viewed as the product $\Flo\x\Eff$ of 
\begin{align}
  \Flo &= T_z\ZS \x \Flo_\Res \x \Flo_\iPort \label{eq:pHdyn4}
  \quad\text{and}\quad
  \Eff = \Flo^* = T_z^*\ZS \x \Eff_\Res \x \Eff_\iPort,
\inter{with the inner product}
  \<\flo,\eff> &= \<{(\flo_\iStor,\flo_\Res,\flo_\iPort)},{(\eff_\iStor,\eff_\Res,\eff_\iPort)}> =
  \({\flo_\iStor},\eff_\iStor) + \({\flo_\Res},\eff_\Res) + \({\flo_\iPort},\eff_\iPort) .\nonumber
\end{align}
The $\<\flo,\eff>= 0$ in the definition represents the {\em lossless} routing of energy between parts of the system.
In the circuit context the Dirac structure is essentially Kirchhoff's laws for a given circuit topology.
That they define lossless routing, is {\em Tellegen's theorem}.

\subsection{\PH system dynamics}

The dynamics of a system is specified by
\begin{align}
  \left(-\dot z(t), \dbd{\Ham}{z}(z(t)),\; \flo_\Res(t), \eff_\Res(t),\; \flo_\iPort(t), \eff_\iPort(t) \right) &\in \Dir(z(t)), \label{eq:pHdyn5a}\\
  (\flo_\Res(t), \eff_\Res(t)) &\in \Res(z(t)), \nonumber
\end{align}
for $t$ in some real interval.
A DAE that fits this structure---it might actually be an ODE---will be called a {\em \pH DAE (pHDAE)}.

A key point from the physical modelling angle is that $\Ham$ is an algebraic function of $z$---the time derivative $\dot z$ must not enter---and must include the whole set of variables on which system energy depends (and preferably nothing but them).
Hence we can use the chain rule in the form
\begin{align*}
  \dot\Ham(z(t)) = \dbd{\Ham}{z}\+\!(z(t))\;\dot z(t) = -\({\flo_\iStor},\eff_\iStor),
\end{align*}
the second equality because from \rf{pHdyn5a}, $-\dot z$ is identified with $\flo_\iStor$, and $\tdbd{\Ham}{z}$ with $\eff_\iStor$.
The Dirac structure enforces the {\em power balance equation}
$\({\flo_\iStor},\eff_\iStor) + \({\flo_\Res},\eff_\Res) + \({\flo_\iPort},\eff_\iPort) = 0$
and thus the DAE records energy flow: any solution of \rf{pHdyn5a} satisfies
\begin{align*}
  \dot\Ham = \({\flo_\Res},\eff_\Res) + \({\flo_\iPort},\eff_\iPort).
\end{align*}
In particular, energy is conserved if there is no dissipation or external power flow.

\medskip
Nothing is said about what variables occur in a pHDAE.
It will have the first order form $f(t,x(t),\dot x(t)) = 0$, where vector $x(t)$, the {\em DAE variables}, contains the {\em state variables} $z(t)\in\ZS$ as a subvector.
What the remaining $x_i$ are, depends on the scenarios chosen to specify the dissipative structure \rf{pHdyn2} and the external port behaviour.

\section{Positive definiteness and passivity results}

\subsection{Matrix positive definiteness}\label{ss:PD}

We cite without proof some results on {\em positive definite} matrices or matrix pairs that are standard, e.g.\ \cite{GolL13,SteS90}, and prove some that are special to this paper. The facts about matrix pairs are present, but not so directly stated, in \cite[\S6.2]{Ria08}.

\begin{definition}~\label{df:PD}
  A matrix $A\in\mathbb{R}^{n\times n}$ is said to be {\em positive definite (PD)} if $x\+Ax>0$ for all $x\in\mathbb{R}^n\setminus\{0\}$.
  A pair of matrices $(B,B')$, $B,B'\in\mathbb{R}^{n\times n}$ is called a {\em positive definite matrix pair (PDMP)} (resp. {\em negative definite matrix pair (NDMP)}) if for all $y,y'\in\mathbb{R}^n$ not both $0$, and $By+B'y'=0$, it holds that $y\+y'>0$ (resp. $y\+y'<0$).
\end{definition}
Note that positive definiteness of a matrix $A$ does not imply that $A$ is symmetric.

\begin{lemma}\label{lm:PDone}~
Let $A\in\mathbb{R}^{n\times n}$.
\begin{compactenum}[(i)]
  \item If $A$ is PD, then $A$ is nonsingular, and $A^{-1}$ is PD.
  \item $A$ is PD iff its symmetric part $(A+A\+)/2$ is PD.
  \item If $A$ is PD and $\tilde A=-\tilde A\+$ is skew-symmetric, then $A+\tilde A$ is PD.
\end{compactenum}
\end{lemma}

\begin{lemma}~\label{lm:PDMP1}
Let $(B,B')$ with $B,B'\in\mathbb{R}^{n\times n}$ be a PDMP (resp. NDMP). Then
\begin{compactenum}[(i)]
  \item $(B',B)$ is a PDMP (resp. NDMP);
  \item $B$ and $B'$ are nonsingular and $-B\`B'$, $-{B'}\`B$ are PD (resp. $B\`B'$, ${B'}\`B$ are PD);
  \item $(PBQ,PB'Q)$ is a PDMP (resp. NDMP) for $P,Q\in\mathbb{R}^{n\times n}$ with $P$ nonsingular, $Q$ orthogonal.
  \item $[C,C'] = [B,B']\;P$ is a PDMP (resp. NDMP), where $P$ is a permutation matrix as in \rf{mixed2}, i.e., $(C,C')$ results from swapping any given set of columns of $B$ with the corresponding columns of $B'$.
\end{compactenum}
\end{lemma}
\begin{proof}
We prove the results for PDMP $(B,B')$.
\begin{compactenum}[(i)]
 \item Follows immediately from ${y'}\+y=y\+y'$.
 \item Suppose $B$ is singular, so there exists $y\ne0$ with $By=0$. Setting $y'=0$ gives  $By+B'y'=0$ with $y,y'$ not both $0$, but $y\+y'=0$, contradicting the definition of a PDMP.
 For the second part let $A=-B\`B'$. For any $y'\ne0$, set $y=-B\`B'y'$.
 Then $y,y'$ are not both $0$ and $By+B'y'=0$ so $0<{y'}\+y = -{y'}\+ B\`B'y' = {y'}\+Ay'$. Hence $A$ is PD.
 \item Let $(C,C')=(PBQ,PB'Q)$ with $P,Q$ be as stated.
  For any $y,y'$, not both $0$, suppose that $0 = Cy+C'y' = P(Bz + B'z')$ where $z=Qy$, $z'=Qy'$. Since $P$ and $Q$ are nonsingular, $z,z'$ are not both $0$ and $Bz + B'z'=0$.
 As $(B,B')$ is a PDMP, we have $0<z\+z' = y\+Q\+Qy'=y\+y'$ since $Q$ is orthogonal.
 Since $y,y'$ were arbitrary this shows $(C,C')$ that is a PDMP.
 \item Let $y,y'\in\R^n$, not both 0, satisfy $Cy+Cy'=0$.
 Defining $[w;w']=P_\beta [y;y']$, we have $Bw+B'w' = Cy+C'y'=0$, with $w,w'$ not both $0$. So $w\+w'>0$ as $(B,B')$ is a PDMP.
 By \rf{mixed2} it holds that $w\+w'=y\+y'$ (as $w_jw_j'$ is replaced by $w_j'w_j$ for each $j$ where a swap was done, but the sum stays the same).
  Hence $y\+y'>0$, which proves $(C,C')$ is a PDMP.
  \end{compactenum}
\end{proof}

\begin{lemma}\label{lm:PDMP2}
Let $B,B'\in\mathbb{R}^{n\times n}$.
If $B$ is nonsingular and $B\`B'$ is PD then $(B,B')$ is a NDMP.  In particular $A\in\mathbb{R}^{n\times n}$ is PD if and only if $(\I,A)$ is a NDMP.
\end{lemma}
\begin{proof}
Follows directly from the previous lemma. 
\end{proof}

\begin{lemma}\label{lm:PD2}
If $M = \bmx{M_{11} &M_{12} \\M_{21} &M_{22}}\mx{\s n_1 \\\s n_2}$ is PD, then
\begin{align} P = \bmx{M_{11} - N\+M_{21}, &M_{12} - N\+M_{22} \\N, &\I} \label{eq:PD2}
\end{align}
is nonsingular for an arbitrary $n_2\x n_1$ matrix $N$.
\end{lemma}
\begin{proof}
  Take any $w = [u;v]$ and suppose $Pw=0$, so
  \begin{align*}
    0 &= (M_{11} - N\+M_{21})u + (M_{12} - N\+M_{22})v,\qquad
    0 = Nu + v 
  \inter{Premultiply the first equation by $u\+$ and substitute $-Nu = v$ given by the second equation:}
    0 &= u\+M_{11}u - u\+N\+M_{21}u + u\+M_{12}v - u\+N\+M_{22}v \\
      &= u\+M_{11}u  + v\+M_{21}u + u\+M_{12}v + v\+M_{22}v \qquad
      = w\+M w,
  \end{align*}
  so $w=0$ by the PD assumption, proving $P$ is nonsingular.
\end{proof}

\begin{lemma}\label{lm:HK}
Let $C,C',S\in\R^{n\x n}$. If $(C,C')$ is a PDMP and $S=-S\+$, then $M =\bmx{\I& S\\C& C'}$ is nonsingular.
\end{lemma}
\begin{proof}
Suppose $z = \bmx{x\\y}$ with $Mz = 0$.
So $x + Sy = 0$, hence $y\+x = y\+(-Sy)) = -y\+Sy = 0$ since $S$ is skew-symmetric.
Also $Cx + C'y = 0$, so as $(C,C')$ is a PDMP we have $x = y = 0$, i.e. $z = 0$. Hence $M$ is nonsingular.
\end{proof}

\subsection{Passivity results}%\label{ss:passivity}
There are various notions of passivity used in the literature. For a discussion of different definitions we refer to \cite{WyaCGGG81}.
Most generally, a system is passive if it can either absorb or store energy, but cannot create energy within itself.
Passivity can be characterized by a dissipation inequality \cite{Wil72}.
For \pH systems, passivity directly follows from the power-balance.
A \pH system is called {\em passive} (with respect to the port variables $\flo_\iPort,\eff_\iPort$ and storage function $\Ham$) if
\[\Ham(z(t_1))-\Ham(z(t_0))\leq \int_{t_0}^{t_1} \eff_\iPort\+\flo_\iPort\, dt\]
for all time instants $t_0\leq t_1$.
Hence, a passive system cannot store more energy in its storage function $\Ham$ than it is supplied with.
\medskip

We also need the concept of passivity for the dissipative structure $\Res$ described by a relation $r(\i_\iDd,\v_\iDd)=0$.

\begin{definition}\cite{ChuDK87}
The resistive structure $\Res$ is called {\em passive} [\/{\em strictly passive}\/] if $\(\i_\iDd,\v_\iDd) \geq 0$ [\/resp.\ $\(\i_\iDd,\v_\iDd) > 0$, except at the origin] for all consistent $(\i_\iDd, \v_\iDd)$. 
\end{definition}

\begin{definition}\label{def:slp}\cite{Chu80,Ria08}
Assume that there exist a disjoint union $\alpha\dcup\beta=\iDd$ such that
\begin{align*}
   \left.\left(\left[\frac{\partial r}{\partial \i_{\alpha}},\frac{\partial r}{\partial \v_{\beta}}\right]\right)\right|_{(\i_\iDd^*,\v_\iDd^*)}\qquad\text{is nonsingular}
\end{align*}
at a given consistent $(\i^*_\iDd,\v^*_\iDd)$. 
The resistive structure is called {\em strictly locally passive at $(\i^*_\iDd, \v^*_\iDd)$} if
\begin{align*}
D_{\cal R}:=-\bmx{\frac{\partial r}{\partial \v_{\alpha}}&\frac{\partial r}{\partial \i_{\beta}}}^{-1}\cdot \bmx{\frac{\partial r}{\partial \i_{\alpha}}&\frac{\partial r}{\partial \v_{\beta}}}.
\end{align*}
 is positive definite at $(\i^*_\iDd, \v^*_\iDd)$.
 The resistive structure $\Res$ is called {\em strictly
locally passive} if $D_{\cal R}$ is positive definite for all consistent $(\i^*_\iDd, \v^*_\iDd)$.
\end{definition}

Strict local passivity implies local passivity, the converse is not true.
Also passive components need not to be strictly locally passive as the following example shows.
\begin{example}
The ideal transformer can be described by the relation
 \[ \bmx{\v_1\\\i_2}=\bmx{0&k\\-k&0}\bmx{\i_1\\\v_2},\quad \text{where $k\in\R$ is the turns ratio}, \]
see e.g. \cite[p. 42ff]{BalB69}.
This resistive structure is passive, since
\[\(\v_\iDd,\i_\iDd)=\v_1\i_1+\v_2\i_2=k\v_2\i_1-k\v_2\i_1=0,\]
but not strictly locally passive as \Asmp5 is not satisfied.
\end{example}

From \cite[Proposition 6.1]{Ria08} it follows that the positive definiteness of $D_{\cal R}$ in Definition \ref{def:slp} does not depend on the specific subsets $\alpha,\beta$. In particular, in a strictly passive circuit one can choose arbitrarily the voltages or currents in every dissipative edge as a variable for a mixed description \rf{mixedForm}. 
Defining 
$$\bar E:=\frac{\partial r}{\partial \i_\iDd}(\i_\iDd^*, \v_\iDd^*), \text{ and } \bar D:=\frac{\partial r}{\partial \v_\iDd}(\i_\iDd^*, \v_\iDd^*)$$
strict local passivity implies that the incremental conductance and resistance matrices
\begin{align*}
 \bar G:=-\bar E^{-1}\bar D \quad \text{and}\quad \bar R:= -\bar D^{-1}\bar E
\end{align*}
are well-defined and positive definite (cf. \cite[p.296]{Ria08}).

In the linear setting  \Asmp5 reads
\begin{align*}
D \v_\iDd + E \i_\iDd = 0\qquad\text{for some }D,E\in\R^{n_\iDd\x n_\iDd}
\end{align*}
with positive definite pair $(D,E)$ so that $G=-E^{-1}D$ (or equivalently $R:= -D^{-1}E$) are positive definite.
Note that in the linear setting these properties are globally defined.
\medskip

Within the \pH setting, the dissipative structure is defined by $\Res \subset \Flo_\Res \x \Eff_\Res$ with $\(\eff_\Res,\flo_\Res)\leq0$.
Non-negativity of the dissipative power flow $\(\eff_\Res,\flo_\Res)$ respects the common convention that power flows {\em from} the ``boundary'' ports (i.e. ports related to $\iPort$ and $\Res$) {\em into} the system and {\em from} the internal network (the Dirac structure) {\em into} the energy storing elements $\iStor$, such that in the power balance equation all flows sum up to zero.
Here, it means the power flows from the energy-dissipating elements into the Dirac structure $\Dir$ (cf. \cite[p.\ 24]{vdSch04}).

By defining $\v_\iDd=\eff_\Res$ and $\i_\iDd=-\flo_\Res$ this setting corresponds to the sign convention in circuit analysis and in the definition of strict local passivity (e.g. the positive definiteness of conductance and resistance matrices) and 
\[\(\eff_\Res, \flo_\Res)\leq 0 \quad \Longleftrightarrow \quad \(\i_\iDd,\v_\iDd)\ge 0. \]

Using e.g.\ the voltage-controlled representation of resistors $\i_\iDd = \bar G \v_\iDd$ this gives
\[\(\i_\iDd,\v_\iDd) =\(\v_\iDd, \bar G \v_\iDd)\ge 0\quad \Longleftrightarrow \quad  \bar G\ge0,\]
i.e., the resistor is passive iff $\bar G$ is positive definite.

\section{\matlab code}\label{sc:matlabcode}

Below are the \matlab functions for the circuit in \exref{ex1} and for the Diode Clipper example generated from the CpH class, both in Model 2 form. They can be used by \matlab's \li{ode15i} solver.

\begin{lstlisting}[numbers=left]
function [f,y] = daefcnExample1_mm(t,x,xp)
%Function to specify circuit Example1 for ode15i.
% Formulation: CpH Model 2, independent, nonlinear, control input-output form.
% Tree [1, 2, 5, 6], cotree [3, 4, 7, 8]
% DAE variables x(1:6)=[qC1 qC2 vG iR φL1 φL2]
% Output variables y(1:2)=[iV vI]

%Parameters
C1 = 5e-6;
C2 = 5e-6;
G  = 1;
R  = 1;
L1 = 0.1;
L2 = 0.1;

%Constitutive relations
v1 = V(t);                i1 = 0;
v2 = x(1)/C1;             i2 = xp(1);
v3 = x(2)/C2;             i3 = xp(2);
v4 = x(3);                i4 = G*x(3);
v5 = R*x(4);              i5 = x(4);
v6 = xp(5);               i6 = x(5)/L1;
v7 = xp(6);               i7 = x(6)/L2;
v8 = 0;                   i8 = I(t);

% Dirac structure, CpH form
f = zeros(6,1); y = zeros(2,1);
y(1) = i1-i3+i4+i7;
f(1) = i2+i3-i7+i8;
f(2) = v3+v1-v2;
f(3) = v4-v1+v5;
f(4) = i5-i4-i7;
f(5) = i6-i7+i8;
f(6) = v7-v1+v2+v5+v6;
y(2) = v8-v2-v6;
y = -y; %to get sign of output right

  function val = V(t)
  val = 10*t.*sin(2*pi*100*t);
  end

  function val = I(t)
  val = 10*sin(10*t);
  end
end
\end{lstlisting}
\begin{lstlisting}[numbers=left]
function [f,y] = daefcnDiodeClipper_mm(t,x,xp)
%Function to specify circuit DiodeClipper for ode15i.
% Formulation: CpH Model 2, independent, nonlinear, control input-output form.
% Tree [1, 2], cotree [3, 4, 5]
% DAE variables x(1:3)=[iR vD1 vD2], 
% Output variables y(1:2)=[iV vI]

%Parameters
omega = 2*pi*1000;
tmax = 0.03; 
Vmax = 2;
Vfac = Vmax/tmax;
R  = 1000;
Is=1e-13;
VT=.025;

%Constitutive relations
v1 = V(t);                i1 = 0;
v2 = R*x(1);              i2 = x(1);
v3 = x(2);                i3 = diode(x(2));
v4 = x(3);                i4 = diode(x(3));
v5 = 0;                   i5 = I(t);

% Dirac structure, CpH form
f = zeros(3,1); y = zeros(2,1);
y(1) = i1+i3-i4+i5;
f(1) = i2+i3-i4+i5;
f(2) =  v3-v1-v2;
f(3) =  v4+v1+v2;
y(2) =  v5-v1-v2;
y = -y; %to get sign of output right

  function val = V(t)
  val = Vfac*t.*sin(omega*t);
  end

  function val = I(t)
  val = 0;
  end

  function val = diode(v)
  val = Is*(expm1(v/VT));
  end
end
\end{lstlisting}

\bibliographystyle{plain}
\bibliography{bibfile}

\end{document}